\documentclass[leqno,10pt]{amsart}
\usepackage{amsmath,amscd,amsthm,amsxtra}
\usepackage{epsfig,graphics,color,colortbl}
\usepackage{amssymb,latexsym}
\usepackage{mathabx}
\usepackage{mathrsfs,eucal,upgreek}
\usepackage[poly,all]{xy}
\usepackage{hyperref}
\usepackage{tikz-cd}
\usepackage{tikz}
\usetikzlibrary{fit,matrix,backgrounds}
\usetikzlibrary{decorations.pathreplacing,calligraphy,matrix}
\usepackage{arydshln}
\usepackage{ytableau}
\usepackage[draft]{todonotes}

\newtheorem{thm}{\bf Theorem}[section]
\newtheorem{df}[thm]{\bf Definition}
\newtheorem{prop}[thm]{\bf Proposition}
\newtheorem{cor}[thm]{\bf Corollary}
\newtheorem{lem}[thm]{\bf Lemma}
\newtheorem{rem}[thm]{\bf Remark}
\newtheorem{ex}[thm]{\bf Example}

\numberwithin{equation}{section}

\newcommand{\bs}{\boldsymbol}
\newcommand{\A}{\mathcal{A}}

\newcommand{\cB}{\mathcal{B}}

\newcommand{\cP}{\mathscr{P}}

\newcommand{\pf}{\noindent{\bfseries Proof. }}
\newcommand{\ov}{\overline}

\newcommand{\ba}{\bs{\rm a}}
\newcommand{\bb}{\bs{\rm b}}

\newcommand{\M}{{\mathcal{M}}}

\newcommand{\gl}{\mathfrak{gl}}
\newcommand{\Z}{\mathbb{Z}}

\newcommand{\te}{\widetilde{e}}
\newcommand{\tf}{\widetilde{f}}
\newcommand{\g}{\mathfrak{g}}
\newcommand{\td}{\widetilde}

\newcommand{\mc}{\mathcal}
\newcommand{\mf}{\mathfrak}
\newcommand{\La}{\Lambda}

\newcommand{\la}{\lambda}

\newcommand{\wh}{\widehat}
\newcommand{\de}{\delta}
\newcommand{\ot}{\otimes}

\newcommand{\blue}[1]{{\color{blue}#1}}
\newcommand{\red}[1]{{\color{red}#1}}

\ytableausetup {mathmode, boxsize=1.2em}

\begin{document}
\title[ ]
{Affine RSK correspondence and crystals of level zero extremal weight modules}

\author{JAE-HOON KWON}

\address{Department of Mathematical Sciences and RIM, Seoul National University, Seoul 08826, Korea}
\email{jaehoonkw@snu.ac.kr}

\author{HYUNSE LEE}

\address{Department of Mathematical Sciences, Seoul National University, Seoul 08826, Korea}
\email{cestlavie813@snu.ac.kr}


\thanks{This work is supported by the National Research Foundation of Korea(NRF) grant funded by the Korea government(MSIT) (No.\,2019R1A2C1084833 and 2020R1A5A1016126).}

\begin{abstract}
We give an affine analogue of the Robison-Schensted-Knuth (RSK) correspondence, which generalizes the affine Robinson-Schensted correspondence by Chmutov-Pylyavskyy-Yudovina. 
The affine RSK map sends a generalized affine permutation of period $(m,n)$ to a pair of tableaux $(P,Q)$ of the same shape, where $P$ belongs to a tensor product of level one perfect Kirillov-Reshetikhin crystals of type $A_{m-1}^{(1)}$, and $Q$ belongs to a crystal of extremal weight module of type $A_{n-1}^{(1)}$ when $m,n\ge 2$. 
We consider two affine crystal structures of types $A_{m-1}^{(1)}$ and $A_{n-1}^{(1)}$ on the set of generalized affine permutations, and show that the affine RSK map preserves the crystal equivalence. We also give a dual affine Robison-Schensted-Knuth correspondence.
\end{abstract}

\maketitle
\setcounter{tocdepth}{1}

\section{Introduction}
The Robinson-Schensted-Knuth (simply RSK) correspondence \cite{Kn} is a fundamental algorithm in algebraic combinatorics with rich applications and generalizations in various areas. In representation theory, we may regard it as a combinatorial aspect of Howe duality \cite{H89}, where the pair of groups $({\rm GL}_m,{\rm GL}_n)$ acts as mutual centralizers on the symmetric or exterior algebra generated by the tensor product of their natural representations. From {the viewpoint} of crystal base theory, the bijection can be further understood as an isomorphism of $(\gl_m,\gl_n)$-bicrystals \cite{La}.

The Robinson-Schensted (simply RS) correspondence, a special case of RSK correspondence restricted to $n\times n$ permutation matrices, gives a bijection from the symmetric group ${\mf S}_n$ on $n$ letters to the set of pairs of standard tableaux of the same shape {with size $n$}, and it provides a combinatorial tool to study the left and right cells of ${\mf S}_n$ in Kazhdan-Lusztig theory \cite{KL}. 

An affine analogue of the RS correspondence has been introduced by Shi \cite{Shi86,Shi91} in the study of affine Kazhdan-Lusztig cells, which maps an affine permutation $w$ to a pair of tabloids $(P(w),Q(w))$, but not necessarily in an injective way.
Recently in \cite{CPY}, Chmutov-Pylyavskyy-Yudovina construct a bijection from the set of extended affine permutations $w$ to the set of triples $(P(w),Q(w),\rho)$, which recovers the Shi's algorithm and characterizes each fiber of $(P(w),Q(w))$, where $\rho$ is an integral vector satisfying a condition {called dominance}. The main ingredient of the affine RS correspondence in \cite{CPY} is the {\em affine matrix-ball construction}, which generalizes the matrix-ball construction in the usual RSK algorithm \cite{Ful}.

In this paper we show that the affine RS correspondence in \cite{CPY} can be naturally extended to a bijection on the set of generalized affine permutations, which also admits an isomorphism of affine crystals.

Let $m$ and $n$ be positive integers greater than $1$. 
Let $\wh{\M}_{m\times n}$ be the set of matrices $A = (a_{ij})_{i,j \in \Z}$ with non-negative integral entries such that $a_{{i+m\, j+n}} = a_{ij}$  for all $i,j\in\Z$, and for each $j$, $a_{ij}=0$ except for finitely many $i$'s.
We show that there exists a bijection 
\begin{equation}\label{eq:affRSK}
\xymatrixcolsep{3pc}\xymatrixrowsep{0.5pc}\xymatrix{
\kappa:\ \ \wh{\M}_{m\times n} \ \ar@{->}[r]  & \ \
\displaystyle{\bigsqcup_{\la\in \cP_m\cap\cP_n} CSST_{[m]}(\la) \times \mc{B}_n(\la)} \\
\ \ A \  \ar@{|->}[r]  &\ \ (P_0,Q)}.
\end{equation}
Here the union is over the partitions $\la$ with length no more than $\min\{m,n\}$,  $CSST_{[m]}(\la)$ is the set of column semistandard tableaux of shape $\la$ with entries from $1$ to $m$, and $\mc{B}_n(\la)$ is the set of tableaux of shape $\la$ with entries in $\Z$, whose rectangular subtableaux are semistandard when $\la$ is decomposed into rectangles along its columns. The bijection is also well-defined when $m=1$ or $n=1$.

The bijection $\kappa$ recovers the affine RS correspondence in \cite{CPY} after a suitable standardization. A key observation is that the dominance condition on $\rho$ in \cite{CPY} fits nicely into description of the connected components in $\mc{B}_n(\la)$, and hence gives the image of $\kappa$ as in \eqref{eq:affRSK}.

We next consider an affine crystal structure related to \eqref{eq:affRSK}. We see that $CSST_{[m]}(\la)$ can be viewed as a finite affine crystal of type $A_{m-1}^{(1)}$ isomorphic to the tensor product of Kirillov-Reshetikhin crystals corresponding to fundamental representations (see, for example, \cite{FSO} and reference therein). Also, $\mc{B}_n(\la)$ has a structure of $A_{n-1}^{(1)}$-crystal isomorphic to the crystal of an extremal weight module associated to the {level-zero} extremal weight corresponding to $\la$ \cite{Kas94',Kas02}. It is shown in \cite{I} using the  theory of semi-infinite LS paths \cite{INS}, but an elementary and self-contained proof of this isomorphism for $\mc{B}(\la)$ is given in this paper. On the other hand, we consider two affine crystal structures of types $A_{m-1}^{(1)}$ and $A_{n-1}^{(1)}$ on $\wh{\M}_{m\times n}$, where the Kashiwara operators for $A_{m-1}^{(1)}$ commute with those for $A_{n-1}^{(1)}$. We expect that the affine crystal structures on $\wh{\M}_{m\times n}$ coincides with the ones in \cite{Lu99}. 

We show that the bijection $\kappa$ preserves the crystal equivalence for these crystal structures. Indeed, we show that $\kappa$
commutes with the Kashiwara operators for both $A_{m-1}$ and $A_{n-1}^{(1)}$, while the operators $\te_0$ and $\tf_0$ on $A$ for $A_{m-1}^{(1)}$ may change the second component $Q$ in \eqref{eq:affRSK} under $\kappa$. 
As a corollary, we have an isomorphism of $(A_{m-1},A_{n-1}^{(1)})$-bicrystals 
\begin{equation*}
 \wh{\M}_{m\times n} \cong \bigoplus_{\la\in \cP_m\cap\cP_n} 
CSST_{[m]}(\la) \times \mc{B}_n(\la),
\end{equation*}
where $\oplus$ means a disjoint union of crystals. 
We remark that the Knuth equivalence for the affine RS correspondence \cite{CLP} can be recovered using the crystal structure here.

We also have a dual affine RSK correspondence, where $\wh{\M}_{m\times n}$ is replaced by the set of binary matrices and $CSST_{[m]}(\la)$ is replaced by the set of row semistandard tableaux of the conjugate shape $\la'$. This can be viewed as a level-zero analogue of the decomposition of the crystal \cite{G,GL} associated to the higher-level $q$-deformed Fock space \cite{U}.

Finally, we remark that an affine RSK correspondence is also given by Imamura- Mucciconi-Sasamoto \cite{IMS}, where the algorithm is given by using dynamics of Sagan-Stanley's skew RSK correspondence \cite{SS}. The algorithm looks different from the one by affine matrix-ball construction in this paper. It would be interesting to compare these two algorithms. A representation theoretic interpretation of the identity corresponding to the bijection in \cite{IMS} is also recently given using representations of current Lie algebras \cite{FKM}. 

The paper is organized as follows: In Section \ref{sec:prel}, we introduce necessary background and notions for the affine RSK algorithm, which are analogous to those introduced in \cite{CPY}. In Section \ref{sec:aff RSK}, we construct the affine RSK correspondence $\kappa$ in \eqref{eq:affRSK} (Theorem \ref{thm:main-1}), where the proof of its well-definedness and bijectivity is given in Section \ref{sec:Proof of affine RSK}.
In Section \ref{sec:aff bicrystal}, we discuss affine crystal structure on both sides of \eqref{eq:affRSK}, and its compatibility with $\kappa$ (Theorem \ref{thm:main-2}). The proof is given in Section \ref{sec:proof of main-2}. In Section \ref{sec:dual aff RSK}, we give a dual affine RSK correspondence.\vskip 2mm

{\bf Acknowledgement} The authors would like to thank Dongkwan Kim for kind explanation on the works \cite{CFKLY,CPY} and the referees for very careful reading and helpful comments.

\section{Preliminaries}\label{sec:prel}
Let us introduce some necessary notations and terminologies following \cite{CPY} with slight modification and generalization.

\subsection{Generalized affine permutations}
Throughout the paper, let $m$ and $n$ denote the positive integers. Let $\Z_{\geq 0}$ denote the set of non-negative integers and $[k]=\{\,1,2,\dots,k\,\}$ for $k\ge 1$ with the usual linear order.

Let
\begin{equation*}
\wh{\M}_{m\times n}=
\left\{\,
A = (a_{ij})_{i,j \in \Z}\,\,\Bigg\vert\,\, 
\begin{array}{l}
(1)\ \text{$a_{ij} \in \Z_{{\geq 0}}$ and $a_{i+m\, j+n} = a_{ij}$ for all $i,j\in\Z$},\\
(2)\ \text{for each $j$, $a_{ij}=0$ except for finitely many $i$'s}
\end{array}
\,\right\}.
\end{equation*}

An {\em extended affine permutation} of $n$ is a bijection $w: \Z \longrightarrow \Z$ satisfying $w(i+n)=w(i)$ for each $i \in \Z$.
The set of extended affine permutations of $n$ is the extended Weyl group of type $A_{n-1}^{(1)}$ if $n > 1$.
Suppose that $m=n$ and a matrix $A \in \wh{\M}_{m\times n}$ satisfies $\sum_{j\in\Z}a_{ij}=\sum_{i\in\Z}a_{ij}=1$ for each $i$ and $j$.
Then $A$ can be viewed as an extended affine permutation $w_A$, which is given by $w_A(i)=j$ for a pair $(i,j)$ such that $a_{ij}=1$.
In this sense, we call an element of $\wh{\M}_{m\times n}$ a {\em generalized affine permutation} of period $(m, n)$.
Following \cite{CPY}, we call $A \in \wh{\M}_{n\times n}$ a {\em partial (extended affine) permutation} if $\sum_{j\in\Z}a_{ij}=\sum_{i\in\Z}a_{ij}\le 1$ for each $i$ and $j$. 
For $A\in \wh{\M}_{m\times n}$, let 
\begin{equation*}
{\rm supp}(A)=\{\,c=(i,j)\in\Z^2\,|\,a_{ij}>0\,\}    
\end{equation*}
be the support of $A\in \wh{\M}_{m\times n}$.
We denote by ${\mathbb O}$ the zero matrix.

An element $c=(i,j)$ in the lattice $\Z^2$ will be often called a {\em cell}. 
We regard each cell as a position of the entry $a_{ij}$ in $A=(a_{ij})\in \wh{\M}_{m\times n}$, where the row index $i$ (resp. the column index $j$) is increasing from top to bottom (resp. from left to right). With this convention, we define partial orders $>_{\tt NW}$, $\ge_{\tt nw}$ and $\le_{\tt ne}$ on $\Z^2$ as follows:
\begin{itemize}
    \item[(1)] $c_1 >_{\tt NW} c_2$ if and only if $i_1 < i_2$ and $j_1 < j_2$,
    \item[(2)] $c_1 \ge_{\tt nw} c_2$ if and only if $i_1 \le i_2$ and $j_1 \le j_2$,
    \item[(3)] $c_1 \le_{\tt ne} c_2$ if and only if $i_1 \ge i_2$ and $j_1 \le j_2$,
\end{itemize}
for $c_1 = (i_1, j_1), c_2 = (i_2, j_2)\in \Z^2$. 
By convention, we use ${\tt N}$ (or ${\tt E}, {\tt W}, {\tt S}$) to emphasize strict inequality, while ${\tt n}$ (or ${\tt e}, {\tt w}, {\tt s}$) allows equality of the co-ordinates of cells. 

Let $\tau=\tau_{m,n}$ denote {the} bijection on $\Z^2$ given by 
\begin{equation*}
\tau(i,j)=(i+m,j+n)\quad ((i,j)\in\Z^2).
\end{equation*} 

\begin{ex}\label{ex:A}{\rm
The following is an example of a generalized affine permutation in $\wh{\M}_{4 \times 5}$.\vskip 2mm


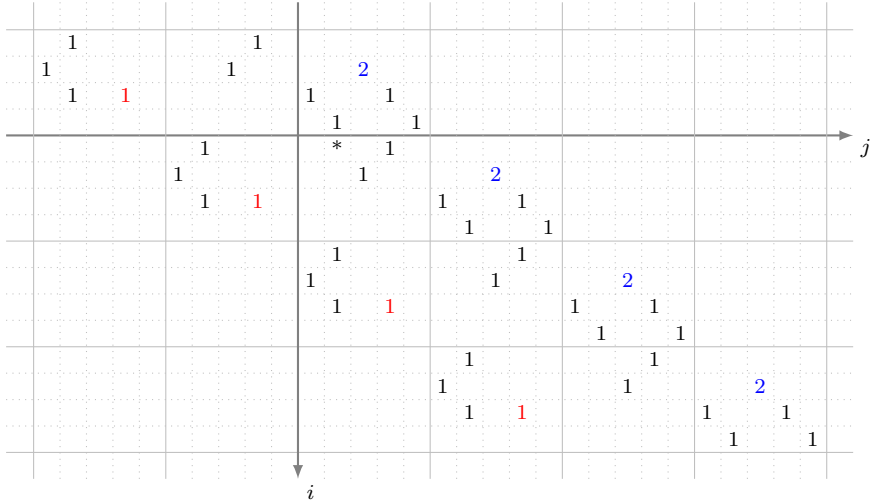
\begin{center}
\begin{tikzpicture}[every node/.style={font=\footnotesize,scale=1}]
\matrix (M)[matrix of math nodes,nodes in empty cells,nodes={rectangle,minimum height=1.0em,minimum width=1.0em,inner sep=0pt,anchor=center,align=center}]
{
& &&&&& &&&&& &&&&& &&&&& &&&&& &&&&& &\\
& &1&&&& &&&1&& &&&&& &&&&& &&&&& &&&&& &\\
& 1&&&&& &&1&&& &&\blue{2}&&& &&&&& &&&&& &&&&& &\\
& &1&&\red{1}&& &&&&& 1&&&1&& &&&&& &&&&& &&&&& &\\
& &&&&& &&&&& &1&&&1& &&&&& &&&&& &&&&& &\\
& &&&&& &1&&&& &*&&1&& &&&&& &&&&& &&&&& &\\
& &&&&& 1&&&&& &&1&&& &&\blue{2}&&& &&&&& &&&&& &\\
& &&&&& &1&&\red{1}&& &&&&& 1&&&1&& &&&&& &&&&& &\\
& &&&&& &&&&& &&&&& &1&&&1& &&&&& &&&&& &\\
& &&&&& &&&&& &1&&&& &&&1&& &&&&& &&&&& &\\
& &&&&& &&&&& 1&&&&& &&1&&& &&\blue{2}&&& &&&&& &\\
& &&&&& &&&&& &1&&\red{1}&& &&&&& 1&&&1&& &&&&& &\\
& &&&&& &&&&& &&&&& &&&&& &1&&&1& &&&&& &\\
& &&&&& &&&&& &&&&& &1&&&& &&&1&& &&&&& &\\
& &&&&& &&&&& &&&&& 1&&&&& &&1&&& &&\blue{2}&&& &\\
& &&&&& &&&&& &&&&& &1&&\red{1}&& &&&&& 1&&&1&& &\\
& &&&&& &&&&& &&&&& &&&&& &&&&& &1&&&1& &\\
& &&&&& &&&&& &&&&& &&&&& &&&&& &&&&& &\\
& &&&&& &&&&& &&&&& &&&&& &&&&& &&&&& &\\
};
\draw[solid,lightgray](M-1-1.south west)--(M-1-32.south east);
\draw[dotted,lightgray](M-2-1.south west)--(M-2-32.south east);
\draw[dotted,lightgray](M-3-1.south west)--(M-3-32.south east);
\draw[dotted,lightgray](M-4-1.south west)--(M-4-32.south east);
\draw[-latex,thick,gray](M-5-1.south west)--(M-5-32.south east);
\node (j) at (M-6-33) {$j$};
\draw[dotted,lightgray](M-6-1.south west)--(M-6-32.south east);
\draw[dotted,lightgray](M-7-1.south west)--(M-7-32.south east);
\draw[dotted,lightgray](M-8-1.south west)--(M-8-32.south east);
\draw[solid,lightgray](M-9-1.south west)--(M-9-32.south east);
\draw[dotted,lightgray](M-10-1.south west)--(M-10-32.south east);
\draw[dotted,lightgray](M-11-1.south west)--(M-11-32.south east);
\draw[dotted,lightgray](M-12-1.south west)--(M-12-32.south east);
\draw[solid,lightgray](M-13-1.south west)--(M-13-32.south east);
\draw[dotted,lightgray](M-14-1.south west)--(M-14-32.south east);
\draw[dotted,lightgray](M-15-1.south west)--(M-15-32.south east);
\draw[dotted,lightgray](M-16-1.south west)--(M-16-32.south east);
\draw[solid,lightgray](M-17-1.south west)--(M-17-32.south east);
\draw[solid,lightgray](M-1-1.north east)--(M-18-1.south east);
\draw[dotted,lightgray](M-1-2.north east)--(M-18-2.south east);
\draw[dotted,lightgray](M-1-3.north east)--(M-18-3.south east);
\draw[dotted,lightgray](M-1-4.north east)--(M-18-4.south east);
\draw[dotted,lightgray](M-1-5.north east)--(M-18-5.south east);
\draw[solid,lightgray](M-1-6.north east)--(M-18-6.south east);
\draw[dotted,lightgray](M-1-7.north east)--(M-18-7.south east);
\draw[dotted,lightgray](M-1-8.north east)--(M-18-8.south east);
\draw[dotted,lightgray](M-1-9.north east)--(M-18-9.south east);
\draw[dotted,lightgray](M-1-10.north east)--(M-18-10.south east);
\draw[-latex,thick,gray](M-1-11.north east)--(M-18-11.south east);
\node (i) at (M-19-12) {$i$};
\draw[dotted,lightgray](M-1-12.north east)--(M-18-12.south east);
\draw[dotted,lightgray](M-1-13.north east)--(M-18-13.south east);
\draw[dotted,lightgray](M-1-14.north east)--(M-18-14.south east);
\draw[dotted,lightgray](M-1-15.north east)--(M-18-15.south east);
\draw[solid,lightgray](M-1-16.north east)--(M-18-16.south east);
\draw[dotted,lightgray](M-1-17.north east)--(M-18-17.south east);
\draw[dotted,lightgray](M-1-18.north east)--(M-18-18.south east);
\draw[dotted,lightgray](M-1-19.north east)--(M-18-19.south east);
\draw[dotted,lightgray](M-1-20.north east)--(M-18-20.south east);
\draw[solid,lightgray](M-1-21.north east)--(M-18-21.south east);
\draw[dotted,lightgray](M-1-22.north east)--(M-18-22.south east);
\draw[dotted,lightgray](M-1-23.north east)--(M-18-23.south east);
\draw[dotted,lightgray](M-1-24.north east)--(M-18-24.south east);
\draw[dotted,lightgray](M-1-25.north east)--(M-18-25.south east);
\draw[solid,lightgray](M-1-26.north east)--(M-18-26.south east);
\draw[dotted,lightgray](M-1-27.north east)--(M-18-27.south east);
\draw[dotted,lightgray](M-1-28.north east)--(M-18-28.south east);
\draw[dotted,lightgray](M-1-29.north east)--(M-18-29.south east);
\draw[dotted,lightgray](M-1-30.north east)--(M-18-30.south east);
\draw[solid,lightgray](M-1-31.north east)--(M-18-31.south east);
\end{tikzpicture}
\end{center}

The co-ordinate of $*$ is {$(1,2)$} and hence the co-ordinates of the cells with red entries are $(-1,-6)$, $(3,-1)$, $(7,4)$, $(11,9)$ from northwest to southeast, while those with blue entries are
$(-2,3)$, $(2,8)$, $(6,13)$, $(10,18)$.


}
\end{ex}

\subsection{Standardization of generalized affine permtations}\label{subsec:standardization}

Let $a=(a_j)_{j \in \Z}$ be a single row matrix with $a_j \in \Z_{\ge 0}$. If $r = \sum_{j \in \Z} a_j < \infty$, we define $[r] \times \infty$ matrix $a^{\circ}=(a^\circ_{ij})_{i\in [r],j\in \Z}$ by
\begin{equation}\label{eq:row std}
a^\circ_{ij}=
\begin{cases}
1 & \text{$\sum_{s=j+1}^{\infty}a_s < i \le \sum_{s=j}^{\infty}a_s$},\\
0 & \text{otherwise}
\end{cases}
\end{equation}
for each $j \in \Z$. For example, if $a=(\dots,0,1,0,3,2,0,\dots)$ with $\sum_{j \in \Z} a_j = 6$, then
\begin{equation*}
a^\circ= 
\begin{pmatrix}
\cdots & 0 & 0 & 0 & 0 & 1 & 0 & \cdots \\
\cdots & 0 & 0 & 0 & 0 & 1 & 0 & \cdots \\
\cdots & 0 & 0 & 0 & 1 & 0 & 0 & \cdots \\
\cdots & 0 & 0 & 0 & 1 & 0 & 0 & \cdots \\
\cdots & 0 & 0 & 0 & 1 & 0 & 0 & \cdots \\
\cdots & 0 & 1 & 0 & 0 & 0 & 0 & \cdots \\
\end{pmatrix}.
\end{equation*} Note that each row of $a^\circ$ has exactly one non-zero entry 1. If $\sum_{j \in \Z} a_j = 0$, we regard $a^\circ$ as an empty matrix with no row (or removing the matrix $a$). 

For a generalized affine permutation $A=(a_{ij})_{i,j \in \Z} \in \wh{\M}_{m\times n}$, we define $A^\circ$ to be the matrix obtained from $A$ by replacing each row $A_i = (a_{ij})_{j \in \Z}$ with $A^\circ_i$ for $i \in \Z$ as in \eqref{eq:row std}.
Similary, we define $A^{\circ '}$ with respect to the columns of $A$, that is, $A^{\circ '}=((A^t)^{\circ})^t$, where $A^t$ denotes the transpose of $A$.
We define the {\em standardization} of $A$ to be
\begin{equation}\label{eq:standardization of A}
    A^{\tt st} = (A^{\circ})^{\circ'}.
\end{equation}
It is straightforward to see that $A^{\tt st}=(A^{\circ})^{\circ'} = (A^{\circ'})^{\circ}$. Note that $A^{\tt st}$ is an extended affine permutation of $K$ if $A$ is non-zero and $K = \sum_{i=1}^{m}\sum_{j\in\Z} a_{ij}=\sum_{j=1}^{n}\sum_{i \in \Z}a_{ij}$.
We regard $\mathbb{O}^{\tt st}$ as an empty matrix with no row and no column.

Let us describe more explicitly the index sets for the rows and columns in $A^{\tt st}$.
Let 
\begin{equation}\label{eq:row and col of A}
    {\rm row}(A)=\left(\ \sum_{j \in \Z}a_{1j},\ \dots,\ \sum_{j \in \Z}a_{mj}\right),\quad
    {\rm col}(A)=\left(\ \sum_{i \in \Z}a_{i1},\ \dots,\ \sum_{i \in \Z}a_{in}\right)
\end{equation}
be the row and column contents of $A$, respectively and write ${\rm row}(A)=(\alpha_1, \dots, \alpha_m)$, ${\rm col}(A)=(\beta_1, \dots, \beta_n)$.
Let $K=\alpha_1 + \cdots + \alpha_m = \beta_1 + \cdots + \beta_n$ and assume that $K \ge 1$, i.e., $A$ is not the zero matrix.
For $i \in [m]$ and $j \in [j]$, let
\begin{equation*}
\begin{split}
    & I_i = \{\, k \in [K]\,|\, \alpha_1 + \dots + \alpha_{i-1} < k \le \alpha_1 + \dots + \alpha_{i-1} + \alpha_{i}\,\},\\
    & J_j = \{\, l \in [K]\,|\, \beta_1 + \dots + \beta_{j-1} < l \le \beta_1 + \dots + \beta_{j-1} + \beta_{j}\,\},
\end{split}
\end{equation*}
where we understand the empty sum is $0$, and let 
\begin{equation*}
    I_{i+sm} = I_i + sK,\quad J_{j+sn} = J_j + sK \quad (s\in\Z).
\end{equation*}
Then we have
\begin{equation}\label{eq:row and column for A perm}
\begin{split}
[K]&=\bigsqcup_{i\in [m]}I_i
     =\bigsqcup_{j\in [n]}J_j,\\
\Z &=\bigsqcup_{i\in Z}I_i
     =\bigsqcup_{j\in Z}J_j.
\end{split}
\end{equation}

\begin{ex}{\rm
Let $m=3$, $n=4$. If $A \in \wh{\M}_{3 \times 4}$ is a generalized affine permutation given by
 \vskip 2mm

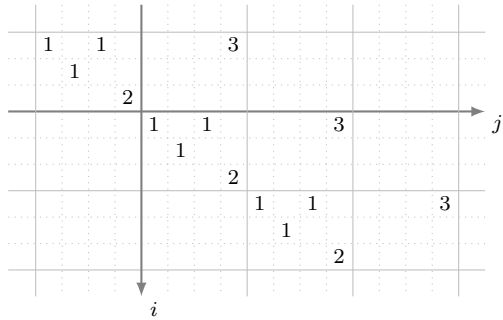
\begin{center}
\begin{tikzpicture}[every node/.style={font=\footnotesize}]
\matrix (M)[matrix of math nodes,nodes in empty cells,nodes={rectangle,minimum height=1.0em,minimum width=1.0em,inner sep=0pt,anchor=center,align=center}]
{
& &&&& &&&& &&&& &&&& &\\
& 1&&1&& &&&3& &&&& &&&& &\\
& &1&&& &&&& &&&& &&&& &\\
& &&&2& &&&& &&&& &&&& &\\
& &&&& 1&&1&& &&&3& &&&& &\\
& &&&& &1&&& &&&& &&&& &\\
& &&&& &&&2& &&&& &&&& &\\
& &&&& &&&& 1&&1&& &&&3& &\\
& &&&& &&&& &1&&& &&&& &\\
& &&&& &&&& &&&2& &&&& &\\
& &&&& &&&& &&&& &&&& &\\
& &&&& &&&& &&&& &&&& &\\
};
\draw[solid,lightgray](M-1-1.south west)--(M-1-18.south east);
\draw[dotted,lightgray](M-2-1.south west)--(M-2-18.south east);
\draw[dotted,lightgray](M-3-1.south west)--(M-3-18.south east);
\draw[-latex,thick,gray](M-4-1.south west)--(M-4-18.south east);
\node (j) at (M-5-19) {$j$};
\draw[dotted,lightgray](M-5-1.south west)--(M-5-18.south east);
\draw[dotted,lightgray](M-6-1.south west)--(M-6-18.south east);
\draw[solid,lightgray](M-7-1.south west)--(M-7-18.south east);
\draw[dotted,lightgray](M-8-1.south west)--(M-8-18.south east);
\draw[dotted,lightgray](M-9-1.south west)--(M-9-18.south east);
\draw[solid,lightgray](M-10-1.south west)--(M-10-18.south east);
\draw[solid,lightgray](M-1-1.north east)--(M-11-1.south east);
\draw[dotted,lightgray](M-1-2.north east)--(M-11-2.south east);
\draw[dotted,lightgray](M-1-3.north east)--(M-11-3.south east);
\draw[dotted,lightgray](M-1-4.north east)--(M-11-4.south east);
\draw[-latex,thick,gray](M-1-5.north east)--(M-11-5.south east);
\node (i) at (M-12-6) {$i$};
\draw[dotted,lightgray](M-1-6.north east)--(M-11-6.south east);
\draw[dotted,lightgray](M-1-7.north east)--(M-11-7.south east);
\draw[dotted,lightgray](M-1-8.north east)--(M-11-8.south east);
\draw[solid,lightgray](M-1-9.north east)--(M-11-9.south east);
\draw[dotted,lightgray](M-1-10.north east)--(M-11-10.south east);
\draw[dotted,lightgray](M-1-11.north east)--(M-11-11.south east);
\draw[dotted,lightgray](M-1-12.north east)--(M-11-12.south east);
\draw[solid,lightgray](M-1-13.north east)--(M-11-13.south east);
\draw[dotted,lightgray](M-1-14.north east)--(M-11-14.south east);
\draw[dotted,lightgray](M-1-15.north east)--(M-11-15.south east);
\draw[dotted,lightgray](M-1-16.north east)--(M-11-16.south east);
\draw[solid,lightgray](M-1-17.north east)--(M-11-17.south east);
\end{tikzpicture}
\end{center}
{then we have
\begin{gather*}
I_1 = \{\,1, 2, 3, 4, 5\,\}, \quad I_2 = \{\, 6\,\}, \quad I_3 = \{\, 7, 8\,\}, \\
J_1 = \{\,1 \,\}, \quad J_2 = \{\, 2\,\}, \quad J_3 = \{\, 3\,\}, \quad J_4 = \{\, 4, 5, 6, 7, 8\,\},
\end{gather*}} and $A^{\tt st}$ is

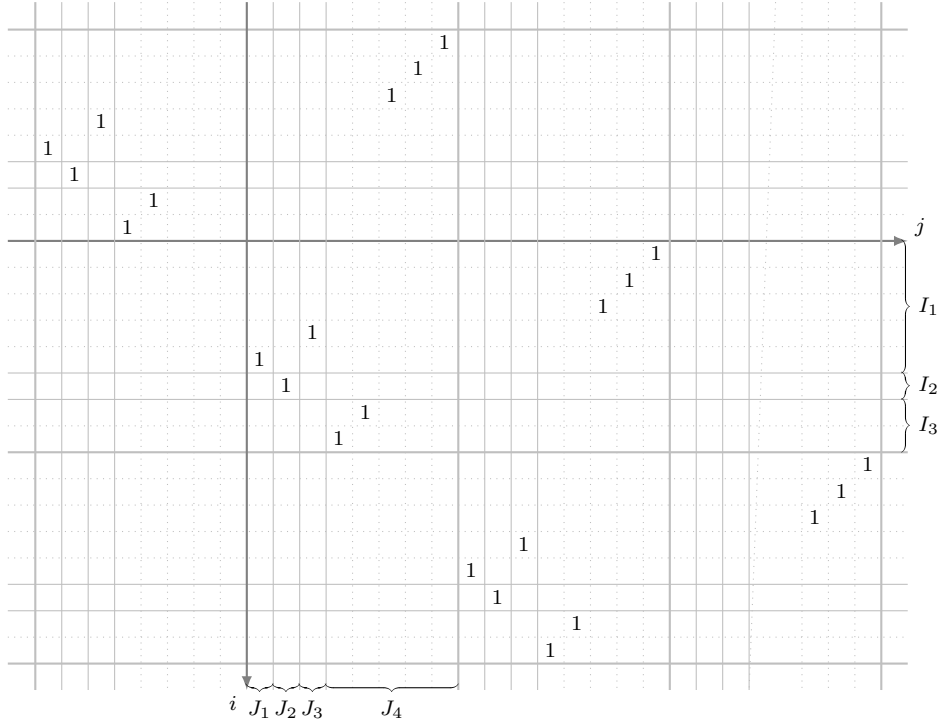
\begin{center}
\begin{tikzpicture}[every node/.style={font=\footnotesize,scale=1},BC/.style = {decorate, decoration={calligraphic brace, amplitude=3pt, raise=1mm},
        pen colour={black}
            }]
\matrix (M)[matrix of math nodes,nodes in empty cells,nodes={rectangle,minimum height=1.0em,minimum width=1.0em,inner sep=0pt,anchor=center,align=center}]
{
& &&&&&&&& &&&&&&&& &&&&&&&& &&&&&&&& &\\
& &&&&&&&& &&&&&&&1& &&&&&&&& &&&&&&&& &\\
& &&&&&&&& &&&&&&1&& &&&&&&&& &&&&&&&& &\\
& &&&&&&&& &&&&&1&&& &&&&&&&& &&&&&&&& &\\
& &&1&&&&&& &&&&&&&& &&&&&&&& &&&&&&&& &\\
& 1&&&&&&&& &&&&&&&& &&&&&&&& &&&&&&&& &\\
& &1&&&&&&& &&&&&&&& &&&&&&&& &&&&&&&& &\\
& &&&&1&&&& &&&&&&&& &&&&&&&& &&&&&&&& &\\
& &&&1&&&&& &&&&&&&& &&&&&&&& &&&&&&&& &\\
& &&&&&&&& &&&&&&&& &&&&&&&1& &&&&&&&& &\\
& &&&&&&&& &&&&&&&& &&&&&&1&& &&&&&&&& &\\
& &&&&&&&& &&&&&&&& &&&&&1&&& &&&&&&&& &\\
& &&&&&&&& &&1&&&&&& &&&&&&&& &&&&&&&& &\\
& &&&&&&&& 1&&&&&&&& &&&&&&&& &&&&&&&& &\\
& &&&&&&&& &1&&&&&&& &&&&&&&& &&&&&&&& &\\
& &&&&&&&& &&&&1&&&& &&&&&&&& &&&&&&&& &\\
& &&&&&&&& &&&1&&&&& &&&&&&&& &&&&&&&& &\\
& &&&&&&&& &&&&&&&& &&&&&&&& &&&&&&&1& &\\
& &&&&&&&& &&&&&&&& &&&&&&&& &&&&&&1&& &\\
& &&&&&&&& &&&&&&&& &&&&&&&& &&&&&1&&& &\\
& &&&&&&&& &&&&&&&& &&1&&&&&& &&&&&&&& &\\
& &&&&&&&& &&&&&&&& 1&&&&&&&& &&&&&&&& &\\
& &&&&&&&& &&&&&&&& &1&&&&&&& &&&&&&&& &\\
& &&&&&&&& &&&&&&&& &&&&1&&&& &&&&&&&& &\\
& &&&&&&&& &&&&&&&& &&&1&&&&& &&&&&&&& &\\
& &&&&&&&& &&&&&&&& &&&&&&&& &&&&&&&& &\\
& &&&&&&&& &&&&&&&& &&&&&&&& &&&&&&&& &\\
};
\draw[thick,lightgray](M-1-1.south west)--(M-1-34.south east);
\draw[dotted,lightgray](M-2-1.south west)--(M-2-34.south east);
\draw[dotted,lightgray](M-3-1.south west)--(M-3-34.south east);
\draw[dotted,lightgray](M-4-1.south west)--(M-4-34.south east);
\draw[dotted,lightgray](M-5-1.south west)--(M-5-34.south east);
\draw[solid,lightgray](M-6-1.south west)--(M-6-34.south east);
\draw[solid,lightgray](M-7-1.south west)--(M-7-34.south east);
\draw[dotted,lightgray](M-8-1.south west)--(M-8-34.south east);
\draw[-latex,thick,gray](M-9-1.south west)--(M-9-34.south east);
\node (j) at (M-9-35) {$j$};
\draw[dotted,lightgray](M-10-1.south west)--(M-10-34.south east);
\draw[dotted,lightgray](M-11-1.south west)--(M-11-34.south east);
\draw[dotted,lightgray](M-12-1.south west)--(M-12-34.south east);
\draw[dotted,lightgray](M-13-1.south west)--(M-13-34.south east);
\draw[solid,lightgray](M-14-1.south west)--(M-14-34.south east);
\draw[solid,lightgray](M-15-1.south west)--(M-15-34.south east);
\draw[dotted,lightgray](M-16-1.south west)--(M-16-34.south east);
\draw[thick,lightgray](M-17-1.south west)--(M-17-34.south east);
\draw[dotted,lightgray](M-18-1.south west)--(M-18-34.south east);
\draw[dotted,lightgray](M-19-1.south west)--(M-19-34.south east);
\draw[dotted,lightgray](M-20-1.south west)--(M-20-34.south east);
\draw[dotted,lightgray](M-21-1.south west)--(M-21-34.south east);
\draw[solid,lightgray](M-22-1.south west)--(M-22-34.south east);
\draw[solid,lightgray](M-23-1.south west)--(M-23-34.south east);
\draw[dotted,lightgray](M-24-1.south west)--(M-24-34.south east);
\draw[thick,lightgray](M-25-1.south west)--(M-25-34.south east);
\draw[thick,lightgray](M-1-1.north east)--(M-26-1.south east);
\draw[solid,lightgray](M-1-2.north east)--(M-26-2.south east);
\draw[solid,lightgray](M-1-3.north east)--(M-26-3.south east);
\draw[solid,lightgray](M-1-4.north east)--(M-26-4.south east);
\draw[dotted,lightgray](M-1-5.north east)--(M-26-5.south east);
\draw[dotted,lightgray](M-1-6.north east)--(M-26-6.south east);
\draw[dotted,lightgray](M-1-7.north east)--(M-26-7.south east);
\draw[dotted,lightgray](M-1-8.north east)--(M-26-8.south east);
\draw[-latex,thick,gray](M-1-9.north east)--(M-26-9.south east);
\node (i) at (M-27-9) {$i$};
\draw[solid,lightgray](M-1-10.north east)--(M-26-10.south east);
\draw[solid,lightgray](M-1-11.north east)--(M-26-11.south east);
\draw[solid,lightgray](M-1-12.north east)--(M-26-12.south east);
\draw[dotted,lightgray](M-1-13.north east)--(M-26-13.south east);
\draw[dotted,lightgray](M-1-14.north east)--(M-26-14.south east);
\draw[dotted,lightgray](M-1-15.north east)--(M-26-15.south east);
\draw[dotted,lightgray](M-1-16.north east)--(M-26-16.south east);
\draw[thick,lightgray](M-1-17.north east)--(M-26-17.south east);
\draw[solid,lightgray](M-1-18.north east)--(M-26-18.south east);
\draw[solid,lightgray](M-1-19.north east)--(M-26-19.south east);
\draw[solid,lightgray](M-1-20.north east)--(M-26-20.south east);
\draw[dotted,lightgray](M-1-21.north east)--(M-26-21.south east);
\draw[dotted,lightgray](M-1-22.north east)--(M-26-22.south east);
\draw[dotted,lightgray](M-1-23.north east)--(M-26-23.south east);
\draw[dotted,lightgray](M-1-24.north east)--(M-26-24.south east);
\draw[thick,lightgray](M-1-25.north east)--(M-26-25.south east);
\draw[solid,lightgray](M-1-26.north east)--(M-26-26.south east);
\draw[solid,lightgray](M-1-27.north east)--(M-26-27.south east);
\draw[solid,lightgray](M-1-28.north east)--(M-26-28.south east);
\draw[dotted,lightgray](M-1-29.north east)--(M-26-28.south east);
\draw[dotted,lightgray](M-1-30.north east)--(M-26-30.south east);
\draw[dotted,lightgray](M-1-31.north east)--(M-26-31.south east);
\draw[dotted,lightgray](M-1-32.north east)--(M-26-32.south east);
\draw[thick,lightgray](M-1-33.north east)--(M-26-33.south east);
\draw[BC] (M-10-34.north) -- node[right=2mm] {$I_1$} (M-14-34.south);
\draw[BC] (M-15-34.north) -- node[right=2mm] {$I_2$} (M-15-34.south);
\draw[BC] (M-16-34.north) -- node[right=2mm] {$I_3$} (M-17-34.south);
\draw[BC] (M-26-10.east) -- node[below=2mm] {$J_1$} (M-26-10.west);
\draw[BC] (M-26-11.east) -- node[below=2mm] {$J_2$} (M-26-11.west);
\draw[BC] (M-26-12.east) -- node[below=2mm] {$J_3$} (M-26-12.west);
\draw[BC] (M-26-17.east) -- node[below=2mm] {$J_4$} (M-26-13.west);

\end{tikzpicture}
\end{center}

}
\end{ex}

\begin{rem}\label{rem:block submatrix}
{\rm 
For $c=(i,j) \in {\rm supp}(A)$, we denote by $A^{\tt st}_c$ the matrix in $\M_{\Z\times\Z}$, which is equal to $A^{\tt st}$ at the positions of $(k,l)\in I_i \times J_j$ and has zero entries elsewhere. 
Then $A^{\tt st}_c$ has an $a_{ij}\times a_{ij}$ block submatrix at $I_i\times J_j$ with $1$ on the antidiagonal, and zero entries elsewhere.
} 
\end{rem}

The following lemmas can be checked easily.
\begin{lem}\label{lem: descending Ast}
Let $c_1, c_2\in {\rm supp}(A^{\tt st})$ be given with $c_1=(i_1,j_1)$ and $c_2=(i_2,j_2)$.
\begin{itemize}
\item[(1)] If $i_1<i_2$ and $c_1,c_2\in I_i \times \Z$ for some $i \in \Z$, then $c_2 \le_{\tt ne} c_1$.

\item[(2)] If $j_1<j_2$ and $c_1,c_2\in \Z \times J_j$ for some $j \in \Z$, then $c_1 \le_{\tt ne} c_2$.
\end{itemize}
\end{lem}

\begin{lem}\label{lem:nw order on Ast}
Let $c_1, c_2\in {\rm supp}(A^{\tt st})$ be given with $c_i\in {\rm supp}\left(A^{\tt st}_{c'_i}\right)$ for some $c'_i\in {\rm supp}(A)$ $(i=1,2)$. Then we have
\begin{itemize}
\item[(1)] $c_2 >_{\tt NW} c_1$ if and only if $c'_2 >_{\tt NW} c'_1$,    

\item[(2)] $c_2 = \tau_{K,K}(c_1)$ implies $c'_2 = \tau_{m,n}(c'_1)$.

\end{itemize}
\end{lem}

\begin{rem}\label{rem:st other versions}
{\rm
Let $a=(a_j)_{j \in \Z}$ with $a_j \in \Z_{\ge 0}$ and assume that $r = \sum_{j \in \Z} a_j < \infty$. We define $[r]\times \Z$ matrix $a^{\bullet}=(a^\bullet_{ij})_{i\in [r],j\in \Z}$ by
\begin{equation*}
a^\bullet_{ij}=
\begin{cases}
1 & \text{$\sum_{s=-\infty}^{j-1}a_s < i \le \sum_{s=-\infty}^{j}a_s$},\\
0 & \text{otherwise}.
\end{cases}
\end{equation*} Then $a^\bullet$ is given by rearranging the rows of $a^\circ$ in a reverse way. For example, if $a=(\dots,0,1,0,3,2,0,\dots)$, we have
\begin{equation*}
a^\bullet= 
\begin{pmatrix}
\cdots & 0 & 1 & 0 & 0 & 0 & 0 & \cdots \\
\cdots & 0 & 0 & 0 & 1 & 0 & 0 & \cdots \\
\cdots & 0 & 0 & 0 & 1 & 0 & 0 & \cdots \\
\cdots & 0 & 0 & 0 & 1 & 0 & 0 & \cdots \\
\cdots & 0 & 0 & 0 & 0 & 1 & 0 & \cdots \\
\cdots & 0 & 0 & 0 & 0 & 1 & 0 & \cdots \\
\end{pmatrix},
\end{equation*}
In general, for $A\in \wh{\M}_{m\times n}$, we define $A^\bullet$ and $A^{\bullet'}=((A^t)^{\bullet})^t$. We may define {other versions} of standardization of $A$ using $\circ$, $\circ'$, $\bullet$ and $\bullet'$.

}\end{rem}

\subsection{Proper numberings}
Suppose that a non-zero $A\in \wh{\M}_{m\times n}$ is given.

\begin{df}[cf.~{{\cite[Definition 3.1]{CPY}}}]\label{def:proper}
{\rm 
A {\em numbering on $A$} is a function $d: {\rm supp}(A)\longrightarrow \Z$. It is called {\em proper} if
\begin{itemize}
\item[(1)] $d(c_2) < d(c_1)$ for $c_1, c_2\in {\rm supp}(A)$ with $c_2 >_{\tt NW} c_1$,

\item[(2)] for $c_1\in {\rm supp}(A)$, there exists $c_2\in {\rm supp}(A)$ with $ c_2>_{\tt NW} c_1$ and $d(c_2) = d(c_1) - 1$.

\end{itemize}
The conditions (1) and (2) are called {\em monotone} and {\em continuous}, respectively.}
\end{df}
The notion of proper numbering on a partial permutation was introduced in \cite{CPY}. It also plays an important role when we describe an affine analogue for the matrix-ball construction of the RSK correspondence in this paper.

Let $d$ be a numbering  on $A$. Let $d^{\tt st}$ be the numbering on ${\rm supp}(A^{\tt st})$ given by
\begin{equation*}
d^{\tt st}(c)=
d(c')\quad \text{if $c\in {\rm supp}(A^{\tt st}_{c'})$ for some $c'\in {\rm supp}(A)$}.
\end{equation*}

\begin{lem}\label{lem:d and dst}
We have the following:
\begin{itemize}
\item[(1)] $d$ is a proper numbering on $A$ if and only if $d^{\tt st}$ is a proper numbering on $A^{\tt st}$, 

\item[(2)] any proper numbering on $A^{\tt st}$ is given by $d^{\tt st}$ for a unique proper numbering $d$ on $A$.
\end{itemize}
\end{lem}
\pf
(1) Since no two cells corresponding to non-zero entries in $A^{\tt st}_{c}$ ($c\in {\rm supp}(A)$) are comparable with respect to $>_{\tt NW}$ (see Remark \ref{rem:block submatrix}), it follows from Lemma \ref{lem:nw order on Ast}(1) that $d$ satisfies the conditions Definition \ref{def:proper}(1) and (2) if and only if $d^{\tt st}$ does so.

(2) Let $d'$ be a proper numbering on $A^{\tt st}$. We claim that $d'=d^{\tt st}$ for some proper numbering $d$ on $A$. By Lemma \ref{lem:nw order on Ast}(1), it suffices to show that $d'$ is constant on $A^{\tt st}_{c}$ for each $c\in {\rm supp}(A)$. Suppose that it does not hold. Then there exist $c_1, c_2\in A^{\tt st}_{c}$ for some $c\in {\rm supp}(A)$ such that $d'(c_1)<d'(c_2)$. By Definition \ref{def:proper}(1), there exists $c_3\in  A^{\tt st}_{c'}$ for some $c' \in {\rm supp}(A)$ such that $c'>_{\tt NW} c$ and $d'(c_3)=d'(c_1)$. Since $c_3>_{\tt NW}c_1$ by Lemma \ref{lem:nw order on Ast}, it is a contradiction. This proves the claim.
\qed\newline

\begin{lem}[cf.~{{\cite[Proposition3.4]{CPY}}}]\label{lem:period is unique}
Under the above hypothesis, the following hold.
\begin{itemize}
\item[(1)] For any proper numbering $d$ on $A$, there exists a positive integer $\ell$, which we call {the {\em period} of $d$}, such that
$d(\tau(c)) = d(c) + \ell$ for $c\in {\rm supp}(A)$.
\item[(2)] If $\ell$ and $\ell'$ are the periods of any two proper numberings $d$ and $d'$ on $A$, respectively, then we have $\ell=\ell'$, which we call the {\em width} of $A$.
\end{itemize}
\end{lem}
\pf The assertions are proved in \cite[Proposition 3.4]{CPY}, when $A$ is a partial permutation. Therefore for arbitrary $A$, (1) and (2) follow from Lemmas \ref{lem:nw order on Ast}(2) and \ref{lem:d and dst}.
\qed
\newline

\subsection{Streams and channels}

Let us verify the existence of a proper numbering on $A\in \wh{\M}_{m\times n}$. This can be done by using standardization and the result in case of partial permutations in \cite{CPY}.

\begin{df}[cf.~{{\cite[Definition 3.20]{CPY}}}]\label{def:stream}
{\rm
A {\em stream} is an infinite collection of cells ${\bf s}=\{c_i\}_{i\in \Z}$, which is invariant under $\tau$ and forms a chain with respect to $>_{\tt NW}$, that is, $c_{i}>_{\tt NW}c_{i+1}$ for $i \in \Z$. 
A {\em flow of a stream ${\bf s}$} is {the} positive integer $l$ such that $\tau(c_i)=c_{i+l}$ for $i \in \Z$. 
A {\em defining data of a stream ${\bf s}=\{c_i=(a_i, b_i)\}_{i\in \Z}$} of flow $l$ is the triple  $(\ba,\bb,r)$, where
\begin{itemize}
\item[(1)] $\ba = (a_{1+r_1},\dots,a_{l+r_1})$ with $1\le a_{1+r_1}< \dots < a_{l+r_1}\le m$,  

\item[(2)] $\bb = (b_{1+r_2},\dots,b_{l+r_2})$ with $1\le b_{1+r_2}< \dots < b_{l+r_2}\le n$, 

\item[(3)] $r=r_1-r_2$.

\end{itemize}
}
\end{df}

Suppose that a non-zero $A\in \wh{\M}_{m\times n}$ is given.
\begin{df}[cf.~{{\cite[Definition 3.6]{CPY}}}]\label{def:channel}
{\rm 
A stream ${\bf s}$ is called a {\em stream of $A$} if ${\bf s}\subset {\rm supp}(A)$. A stream ${\bf s}$ of $A$ is called a {\em channel of $A$} if its flow is maximal among the streams of $A$.}
\end{df}

Let $C=\{c_i\}_{i \in \Z}$ be a channel of $A$ and let $c \in {\rm supp}(A)$. Let $k$ be the maximal integer such that $c_k >_{\tt NW} c$. The maximal property of channel forces that $c \ngtr_{\tt NW} c_{k+1}$. 
This implies either $c_{k+1} \le_{\tt ne} c$ or $c \le_{\tt ne} c_{k+1}$. 
In other words, we have 
\begin{equation*}
 {\rm supp}(A) = C_{\tt ne} \cup C_{\tt sw},
\end{equation*}
where
\begin{equation}\label{eq:C_ne C_sw}
\begin{split}
C_{\tt ne} = \{\,c \in {\rm supp}(A)\,  \vert\, \text{$c' \le_{\tt ne} c$ for some $c' \in C$}\,\},\\
C_{\tt sw} = \{\,c \in {\rm supp}(A)\,  \vert\, \text{$c \le_{\tt ne} c'$ for some $c' \in C$}\,\},
\end{split} 
\end{equation}
and $C_{\tt ne} \cap C_{\tt sw} = C$.
Let $\mc{C}_A$ denote the set of channels of $A$, which is finite.
For $C_1, C_2\in \mc{C}_A$, we define 
\begin{equation}\label{eq:SW order on channels}
\text{$C_1 \succcurlyeq_{\tt sw} C_2$ if and only if $C_1 \subset (C_2)_{\tt sw}$. }
\end{equation}
Then $\succcurlyeq_{\tt sw}$ is a partial order on $\mc{C}_A$.

Let ${\bf s}=\{c_i\}_{i\in\Z}$ be a stream of $A$. Let ${\bf s}^{\tt st}=\{c^{\tt st}_i\}_{i\in\Z}$ be the collection of cells 
such that $c^{\tt st}_i$ is the cell corresponding to the most southwest $1$ in the block submatrix $A^{\tt st}_{c_i}$ for $i \in \Z$.

On the other hand, let ${\bf t}=\{c'_i\}_{i\in\Z}$ be a stream of $A^{\tt st}$. Let ${\bf t}^{{}^{{}_{\swarrow}}}$ be the collection of cells obtained from ${\bf t}$ by replacing each $c'_i$ by the cell corresponding to the most southwest non-zero entry in $A^{\tt st}_{c}$ if $c'_i\in {\rm supp}\left( A^{\tt st}_{c} \right)$ for $c\in {\rm supp}(A)$. Also let $\ov{\bf t}=\{\ov{c}_i\}_{i\in \Z}$ be the collection of cells in ${\rm supp}(A)$ such that $c_i\in {\rm supp}\left(A^{\tt st}_{\ov{c}_i}\right)$ for $i\in\Z$.

\begin{lem}\label{lem:streams comparison}
Under the above hypothesis, we have the following:
\begin{itemize}
\item[(1)] ${\bf s}^{\tt st}$ is a stream of $A^{\tt st}$ with the same flow as ${\bf s}$,

\item[(2)] ${\bf t}^{{}^{{}_{\swarrow}}}$ and $\ov{\bf t}$ are streams of $A^{\tt st}$ and $A$ with the same flow as ${\bf t}$, respectively,

\item[(3)] ${\bf t}$ is a channel if and only if ${\bf t}^{{}^{{}_{\swarrow}}}$  is a channel, in which case we have ${\bf t}^{{}^{{}_{\swarrow}}}\succcurlyeq_{\tt sw}{\bf t}$,

\item[(4)] ${\bf s}$ is a channel if and only if ${\bf s}^{\tt st}$ is a channel, 

\item[(5)] for any channel ${\bf s}'$ of $A$, ${\bf s}\succcurlyeq_{\tt sw}{\bf s}'$ if and only if ${\bf s}^{\tt st}\succcurlyeq_{\tt sw}{{\bf s}'}^{\tt st}$.
\end{itemize}
\end{lem}
\pf It is straightforward to check (1) and (2). The assertion (3) follows from the fact that no two cells of non-zero entries in $A^{\tt st}_{c}$ for $c\in {\rm supp}(A)$ belong to a stream of $A^{\tt st}$.
The assertion (4) follows from the fact that for any stream ${\bf t}$ of $A^{\tt st}$, there exists a unique stream ${\bf s}$ of $A$ such that ${\bf s}^{\tt st}={\bf t}^{{}^{{}_{\swarrow}}}$. This implies (5) immediately.
\qed

\begin{prop}[cf.~{{\cite[Proposition 3.13]{CPY}}}]\label{prop:lattice property}
The set $\mc{C}_A$ has a greatest element with respect to $\succcurlyeq_{\tt sw}$, which we denote by $C_A^{\tt sw}$. 
\end{prop}
\pf As in the proof of Lemma \ref{lem:period is unique}, we reduce the proof to the case when $A$ is a partial permutation.
Let $C_1$ and $C_2$ be two channels of $A$. By Lemma \ref{lem:streams comparison}, $C_1^{\tt st}$ and $C_2^{\tt st}$ are {two channels} of $A^{\tt st}$. By \cite[Proposition 3.13]{CPY}, there exists a channel $C$ of $A^{\tt st}$ such that $C\succcurlyeq_{\tt sw} C_1^{\tt st}, C_2^{\tt st}$. Let $C_0$ be a channel of $A$ such that $C_0^{\tt st} = {C}^{{}^{{}_{\swarrow}}}$. Hence we have $C_0^{\tt st} = {C}^{{}^{{}_{\swarrow}}}\succcurlyeq_{\tt sw} C\succcurlyeq_{\tt sw} C_1^{\tt st}, C_2^{\tt st}$, which implies that $C_0\succcurlyeq_{\tt sw} C_1, C_2$ by Lemma \ref{lem:streams comparison}. This completes the proof.
\qed\newline

\begin{df}[cf.~{{\cite[Definition 3.14]{CPY}}}]
{\rm
We call $C_A^{\tt sw}$ the \textit{southwest channel of $A$}.}
\end{df}

\begin{cor}
Under the above hypothesis, we have 
$$\left(C_A^{\tt sw}\right)^{\tt st}=C_{A^{\tt st}}^{\tt sw}.$$
\end{cor}
\pf It follows from Lemma \ref{lem:streams comparison}(4) and (5).
\qed\newline

Let $C = \{c_i\}_{i\in \Z}$ be a channel of $A$ and let ${d}_0$ be a numbering on $C$ given by ${d}_0(c_i)=i$ for $i\in\Z$. For $c\in {\rm supp}(A)$, we define
\begin{equation}\label{eq:channel numbering}
d^{C}_{A}(c) = 
{\rm sup}
\left\{\,{d}_0(c'_k) + k \ \Bigg\vert\
\begin{array}{l}
\text{$c'_k>_{\tt NW} \dots >_{\tt NW} c'_0$ is a chain in ${\rm supp}(A)$}\\
\ \text{ ($k\ge 0$) such that $c'_0 = c$ and $c'_k \in C$}
\end{array}
\,\right\}.
\end{equation}

\begin{prop}[cf.~{{\cite[Proposition 3.10]{CPY}}}]\label{prop:ch numbering}
The numbering $d^{C}_{A}$ on $A$ is a well-defined proper numbering. Moreover, we have $d^{C}_A(c)={d}_0(c)$ for $c \in C$.
\end{prop}
\pf Put $d=d^{C}_A$. 
We claim that 
\begin{equation}\label{eq:sw proper A and Ast}
d^{\tt st}=d^{C^{\tt st}}_{A^{\tt st}},
\end{equation} 
where $d^{C^{\tt st}}_{A^{\tt st}}$ is defined in the same way as in \eqref{eq:channel numbering} with respect to $C^{\tt st}$.

Let $c\in {\rm supp}(A^{\tt st})$ be given. 
Let $\ov{c}\in {\rm supp}(A)$ such that $c\in {\rm supp}(A^{\tt st}_{\ov{c}})$. 
Suppose that we are given a chain $c_k>_{\tt NW} \dots >_{\tt NW} c_0$ in ${\rm supp}(A)$ from $c_0 = \ov{c}$ to $c_k \in C$. Let $c'_i$ be a cell in ${\rm supp}(A^{\tt st}_{c_i})$ ($1 \le i \le k$) such that $c'_k \in C^{\tt st}$. Then we have a chain $c'_k>_{\tt NW} \dots >_{\tt NW} c'_0$ in ${\rm supp}(A^{\tt st})$ from $c'_0 = c$ to $c_k \in C^{\tt st}$. Conversely, if we have a chain of length $k+1$ in ${\rm supp}(A^{\tt st})$ from $c$ to $C^{\tt st}$, then we can associate a chain of length $k+1$ in ${\rm supp}(A)$ from $\ov{c}$ to $C$ by taking representatives of the block submatrix to which each cell in the chain belongs to. This implies that $d^{\tt st}(c)=d^{C^{\tt st}}_{A^{\tt st}}(c)$, which proves the claim. 
Hence, $d$ is also proper by Lemma \ref{lem:d and dst}(1) since $d^{C^{\tt st}}_{A^{\tt st}}$ is proper by \cite[Proposition 3.10]{CPY}. It follows from $d^{\tt st}=d^{C^{\tt st}}_{A^{\tt st}}$ and \cite[Proposition 11.3]{CPY} that
$d(c)={d}_0(c)$ $(c \in C)$.
\qed

\begin{rem}\label{rem:period=flow=width}
{{\rm
Note that the period of the channel numbering $d_A^C$ or, equivalently, the width of $A$ is equal to the flow of a channel $C$ of $A$.
}}
\end{rem}

\begin{df}{\rm
When $C = C^{\tt sw}_A$, we write $d^{\tt sw}_A = d^{C^{\tt sw}_A}_A$ for short, and call it the {\em southwest channel numbering on $A$}.

}
\end{df}

\begin{ex}\label{ex:swc numbering of A}{\rm
Let $A$ be the generalized affine permutation in Example \ref{ex:A}. The southwest channel of $A$ is given by
$$ C^{\tt sw}_A = \left\{\, \cdots >_{\tt NW} \tau^{-1}(4,7) >_{\tt NW} (2,3) >_{\tt NW} (3,6) >_{\tt NW} (4,7) >_{\tt NW} \tau(2,3) >_{\tt NW} \cdots\, \right\}, $$ which has flow 3. Then the southwest channel numbering $d^{\tt sw}_A$ is represented as follows:\vskip 2mm

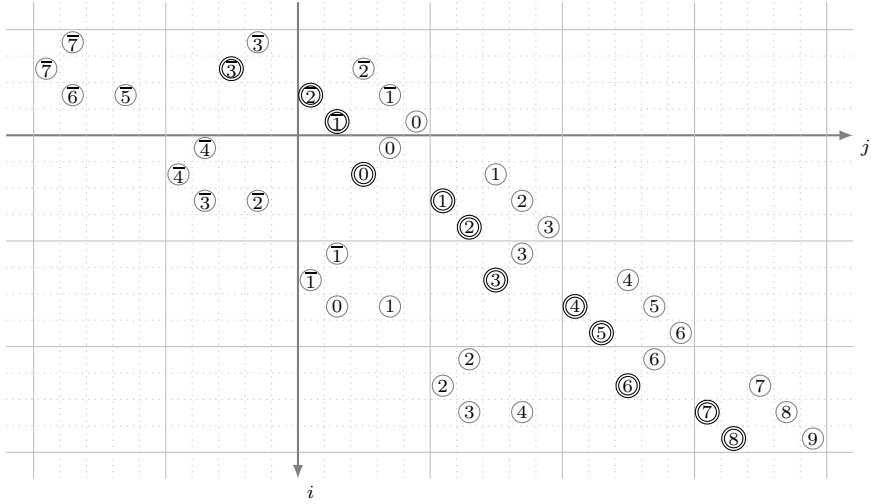
\begin{center}
\begin{tikzpicture}[every node/.style={font=\scriptsize}]
\matrix (M)[matrix of math nodes,nodes in empty cells,nodes={rectangle,minimum height=1.0em,minimum width=1.0em,inner sep=0pt,anchor=center,align=center}]
{
& &&&&& &&&&& &&&&& &&&&& &&&&& &&&&& &\\
& &$\ov{7}$&&&& &&&$\ov{3}$&& &&&&& &&&&& &&&&& &&&&& &\\
& $\ov{7}$&&&&& &&$\ov{3}$&&& &&$\ov{2}$&&& &&&&& &&&&& &&&&& &\\
& &$\ov{6}$&&$\ov{5}$&& &&&&& $\ov{2}$&&&$\ov{1}$&& &&&&& &&&&& &&&&& &\\
& &&&&& &&&&& &$\ov{1}$&&&0& &&&&& &&&&& &&&&& &\\
& &&&&& &$\ov{4}$&&&& &&&0&& &&&&& &&&&& &&&&& &\\
& &&&&& $\ov{4}$&&&&& &&0&&& &&1&&& &&&&& &&&&& &\\
& &&&&& &$\ov{3}$&&$\ov{2}$&& &&&&& 1&&&2&& &&&&& &&&&& &\\
& &&&&& &&&&& &&&&& &2&&&3& &&&&& &&&&& &\\
& &&&&& &&&&& &$\ov{1}$&&&& &&&3&& &&&&& &&&&& &\\
& &&&&& &&&&& $\ov{1}$&&&&& &&3&&& &&4&&& &&&&& &\\
& &&&&& &&&&& &0&&1&& &&&&& 4&&&5&& &&&&& &\\
& &&&&& &&&&& &&&&& &&&&& &5&&&6& &&&&& &\\
& &&&&& &&&&& &&&&& &2&&&& &&&6&& &&&&& &\\
& &&&&& &&&&& &&&&& 2&&&&& &&6&&& &&7&&& &\\
& &&&&& &&&&& &&&&& &3&&4&& &&&&& 7&&&8&& &\\
& &&&&& &&&&& &&&&& &&&&& &&&&& &8&&&9& &\\
& &&&&& &&&&& &&&&& &&&&& &&&&& &&&&& &\\
& &&&&& &&&&& &&&&& &&&&& &&&&& &&&&& &\\
};
\draw[solid,lightgray](M-1-1.south west)--(M-1-32.south east);
\draw[dotted,lightgray](M-2-1.south west)--(M-2-32.south east);
\draw[dotted,lightgray](M-3-1.south west)--(M-3-32.south east);
\draw[dotted,lightgray](M-4-1.south west)--(M-4-32.south east);
\draw[-latex,thick,gray](M-5-1.south west)--(M-5-32.south east);
\node (j) at (M-6-33) {$j$};
\draw[dotted,lightgray](M-6-1.south west)--(M-6-32.south east);
\draw[dotted,lightgray](M-7-1.south west)--(M-7-32.south east);
\draw[dotted,lightgray](M-8-1.south west)--(M-8-32.south east);
\draw[solid,lightgray](M-9-1.south west)--(M-9-32.south east);
\draw[dotted,lightgray](M-10-1.south west)--(M-10-32.south east);
\draw[dotted,lightgray](M-11-1.south west)--(M-11-32.south east);
\draw[dotted,lightgray](M-12-1.south west)--(M-12-32.south east);
\draw[solid,lightgray](M-13-1.south west)--(M-13-32.south east);
\draw[dotted,lightgray](M-14-1.south west)--(M-14-32.south east);
\draw[dotted,lightgray](M-15-1.south west)--(M-15-32.south east);
\draw[dotted,lightgray](M-16-1.south west)--(M-16-32.south east);
\draw[solid,lightgray](M-17-1.south west)--(M-17-32.south east);
\draw[solid,lightgray](M-1-1.north east)--(M-18-1.south east);
\draw[dotted,lightgray](M-1-2.north east)--(M-18-2.south east);
\draw[dotted,lightgray](M-1-3.north east)--(M-18-3.south east);
\draw[dotted,lightgray](M-1-4.north east)--(M-18-4.south east);
\draw[dotted,lightgray](M-1-5.north east)--(M-18-5.south east);
\draw[solid,lightgray](M-1-6.north east)--(M-18-6.south east);
\draw[dotted,lightgray](M-1-7.north east)--(M-18-7.south east);
\draw[dotted,lightgray](M-1-8.north east)--(M-18-8.south east);
\draw[dotted,lightgray](M-1-9.north east)--(M-18-9.south east);
\draw[dotted,lightgray](M-1-10.north east)--(M-18-10.south east);
\draw[-latex,thick,gray](M-1-11.north east)--(M-18-11.south east);
\node (i) at (M-19-12) {$i$};
\draw[dotted,lightgray](M-1-12.north east)--(M-18-12.south east);
\draw[dotted,lightgray](M-1-13.north east)--(M-18-13.south east);
\draw[dotted,lightgray](M-1-14.north east)--(M-18-14.south east);
\draw[dotted,lightgray](M-1-15.north east)--(M-18-15.south east);
\draw[solid,lightgray](M-1-16.north east)--(M-18-16.south east);
\draw[dotted,lightgray](M-1-17.north east)--(M-18-17.south east);
\draw[dotted,lightgray](M-1-18.north east)--(M-18-18.south east);
\draw[dotted,lightgray](M-1-19.north east)--(M-18-19.south east);
\draw[dotted,lightgray](M-1-20.north east)--(M-18-20.south east);
\draw[solid,lightgray](M-1-21.north east)--(M-18-21.south east);
\draw[dotted,lightgray](M-1-22.north east)--(M-18-22.south east);
\draw[dotted,lightgray](M-1-23.north east)--(M-18-23.south east);
\draw[dotted,lightgray](M-1-24.north east)--(M-18-24.south east);
\draw[dotted,lightgray](M-1-25.north east)--(M-18-25.south east);
\draw[solid,lightgray](M-1-26.north east)--(M-18-26.south east);
\draw[dotted,lightgray](M-1-27.north east)--(M-18-27.south east);
\draw[dotted,lightgray](M-1-28.north east)--(M-18-28.south east);
\draw[dotted,lightgray](M-1-29.north east)--(M-18-29.south east);
\draw[dotted,lightgray](M-1-30.north east)--(M-18-30.south east);
\draw[solid,lightgray](M-1-31.north east)--(M-18-31.south east);
\draw[gray](M-2-3)circle(0.4em);
\draw[gray](M-2-10)circle(0.4em);
\draw[gray](M-3-2)circle(0.4em);
\draw[double](M-3-9)circle(0.4em);
\draw[gray](M-3-14)circle(0.4em);
\draw[gray](M-4-3)circle(0.4em);
\draw[gray](M-4-5)circle(0.4em);
\draw[double](M-4-12)circle(0.4em);
\draw[gray](M-4-15)circle(0.4em);
\draw[double](M-5-13)circle(0.4em);
\draw[gray](M-5-16)circle(0.4em);
\draw[gray](M-6-8)circle(0.4em);
\draw[gray](M-6-15)circle(0.4em);
\draw[gray](M-7-7)circle(0.4em);
\draw[double](M-7-14)circle(0.4em);
\draw[gray](M-7-19)circle(0.4em);
\draw[gray](M-8-8)circle(0.4em);
\draw[gray](M-8-10)circle(0.4em);
\draw[double](M-8-17)circle(0.4em);
\draw[gray](M-8-20)circle(0.4em);
\draw[double](M-9-18)circle(0.4em);
\draw[gray](M-9-21)circle(0.4em);
\draw[gray](M-10-13)circle(0.4em);
\draw[gray](M-10-20)circle(0.4em);
\draw[gray](M-11-12)circle(0.4em);
\draw[double](M-11-19)circle(0.4em);
\draw[gray](M-11-24)circle(0.4em);
\draw[gray](M-12-13)circle(0.4em);
\draw[gray](M-12-15)circle(0.4em);
\draw[double](M-12-22)circle(0.4em);
\draw[gray](M-12-25)circle(0.4em);
\draw[double](M-13-23)circle(0.4em);
\draw[gray](M-13-26)circle(0.4em);
\draw[gray](M-14-18)circle(0.4em);
\draw[gray](M-14-25)circle(0.4em);
\draw[gray](M-15-17)circle(0.4em);
\draw[double](M-15-24)circle(0.4em);
\draw[gray](M-15-29)circle(0.4em);
\draw[gray](M-16-18)circle(0.4em);
\draw[gray](M-16-20)circle(0.4em);
\draw[double](M-16-27)circle(0.4em);
\draw[gray](M-16-30)circle(0.4em);
\draw[double](M-17-28)circle(0.4em);
\draw[gray](M-17-31)circle(0.4em);
\end{tikzpicture}
\end{center}
Here each circle represents a cell in ${\rm supp}(A)$ and the double circles forms the southwest channel $C^{\tt sw}_A$. The number in each circle is the numbering $d^{\tt sw}_A$ of that cell, where $\ov{k}$ means $-k$. For example, $d^{\tt sw}_A(2,8) = d^{\tt sw}_A(3,6) = d^{\tt sw}_A(7,4) = 1$.
}
\end{ex}

\begin{rem}\label{rem:ne ch}
{\rm
Consider the southwest channel $C^{\tt sw}_{A^t}$ of $A^t$. Then
\begin{equation*}
\begin{split}
C^{\tt ne}_{A} &= \{\, (i, j) \, \vert \, (j, i) \in C^{\tt sw}_{A^t} \,\} \\
\end{split}
\end{equation*}
is the minimal element in $\mc{C}_A$ with respect to $\ge_{\tt sw}$, which is called the {\em northeast channel of $A$}.
Let $d^{\tt ne}_A$ be the channel numbering corresponding to $C^{\tt ne}_A$. Then it follows from definition that
$d^{\tt sw}_{A^t}(j, i) = d^{\tt ne}_{A}(i, j)$ for $(j,i) \in {\rm supp}(A^t)$.
}
\end{rem}

The following lemma gives a characterization of a channel numbering (cf.~\cite[Remark 11.8]{CPY}).

\begin{lem}\label{lem:characterization of ch numbering}
Let $C$ be a channel of $A$. 
Let $d$ be a proper numbering on $A$ such that $d(c) = d^{C}_{A}(c)$ for $c \in C$. Then the following are equivalent:
\begin{itemize}
\item[(1)] $d=d^{C}_{A}$,
\item[(2)] for $c \in {\rm supp}(A)$, there exists a chain $c_k>_{\tt NW} \dots >_{\tt NW} c_0$ in ${\rm supp}(A)$ such that $c_0 = c$, $c_k \in C$ and $d(c_i) = d(c) - i$ for $0\le i \le k$,
\item[(3)] if $d'$ is a proper numbering such that $d'(c)=d^C_A(c)$ for $c \in C$, then we have $d(c) \le d'(c)$ for every $c \in {\rm supp}(A)$.
\end{itemize} 
\end{lem}

\pf Suppose that (1) holds. Let $c_k>_{\tt NW} \dots >_{\tt NW} c_0=c$ be a chain which gives the maximum value $d_0(c_k) + k$ in \eqref{eq:channel numbering}. Since $d$ is monotone, we have
\begin{equation}\label{eq:inequality}
d(c_k) + k \le d(c_{k-1}) + k - 1 \le \dots \le d(c)=d_0(c_k) + k. 
\end{equation}
Since $c_k \in C$, we have $d(c_k) = d_0(c_k)$ by Proposition \ref{prop:ch numbering}. Thus all the inequalities in \eqref{eq:inequality} are in fact equalities and hence, $d(c_i) = d(c) - i$ for $0\le i \le k$. This implies (2).

Suppose that (2) holds. 
For $c \in {\rm supp}(A)$, let $c_k>_{\tt NW} \dots >_{\tt NW} c_0=c$ be a chain satisfying the condition in (2). 
Let $d'$ be a proper numbering such that $d'=d^C_A$ on $C$.
Along this chain, we have
\begin{equation*}
d'(c_k) + k \le d'(c_{k-1}) + k - 1 \le \dots \le d'(c).
\end{equation*}
from the monotonicity of $d'$. Since $d'(c_k) = d(c_k)$ we conclude that $d(c) = d(c_k) + k = d'(c_k) + k \le d'(c)$. This implies (3).

Suppose that (3) holds. Then, in particular, we have $d(c) \le d^C_A(c)$ for $c \in {\rm supp}(A)$ by letting $d'=d^C_A$. Let $c_k>_{\tt NW} \dots >_{\tt NW} c_0=c$ be a chain which gives the maximal value $d_A^{C}(c)$. We have
$d(c_k) + k \le d(c_{k-1}) + k - 1 \le \dots \le d(c)$ from the monotonicity of $d$. Then we see that $d^C_A(c) = d^C_A(c_k) + k = d(c_k) + k \le d(c)$.
Hence $d(c) = d^C_A(c)$.
\qed

\begin{rem}{\rm
By Lemma \ref{lem:characterization of ch numbering}, we may understand that the channel numbering $d^C_A$ is the proper numbering with minimal values among the proper numberings $d$ which coincide with $d_0$ on $C$.
} 
\end{rem}

\begin{cor}\label{cor:two ch numbering with intersection}
Let $C_1$, $C_2$ be channels of $A$ such that $C_1\cap C_2\neq \emptyset$. If $d^{C_1}_A=d^{C_2}_A$ on $C_1 \cap C_2$, then we have $d^{C_1}_A=d^{C_2}_A$.
\end{cor}
\pf
Let 
$C_1 = \{ \cdots >_{\tt NW} c_0 >_{\tt NW} \cdots >_{\tt NW} c_l = \tau(c_0) >_{\tt NW} \cdots \}$,
where $l$ is the common flow of $C_1$ and $C_2$. We may assume that $c_0 \in C_1\cap C_2$. 
Along $C_1$, we have
\begin{equation*}
    d^{C_2}_A(c_0) + l \le d^{C_2}_A(c_1) + (l-1) \le \cdots \le d^{C_2}_A(c_l).
\end{equation*}
By Lemma \ref{lem:period is unique}, we see that $d^{C_2}_A(c_l) = d^{C_2}_A(c_0) + l$. Hence all the inequalities above are equalities, and for $0 \le i \le l$
$$d^{C_2}_A(c_i) = d^{C_2}_A(c_0) + i = d^{C_1}_A(c_0) + i = d^{C_1}_A(c_i).$$ 
By Lemma \ref{lem:characterization of ch numbering}(3), we have $d^{C_1}_A(c) \le d^{C_2}_A(c)$ for $c \in {\rm supp}(A)$.
The reverse inequality can be obtained in the same way. Hence we obtain $d^{C_1}_A=d^{C_2}_A$.
\qed
\newline

\subsection{Tableaux}\label{subsec:tableaux}
Let $\cP$ be the set of partitions $\la=(\la_1,\la_2,\cdots)$. Following \cite{Ful}, we regard a partition $\lambda$ as its Young diagram.
Let $\la'=(\la'_1,\la'_2,\cdots)$ be the conjugate of $\la$, $\ell(\la)$ the length of $\la$, and $|\la|=\sum_{i\ge 1}\la_i$. For $n\geq 1$, let $\cP_{n}=\{\,\la\in\cP\,|\,\ell(\la)\leq n\,\}$. 

{Let $\mc{A}$ be $\Z$ or $[k]$ ($k\ge 1$) with the usual linear orders.}
For $\la\in\cP$, let $SST_\A(\la)$ be the set of semistandard (or $\A$-semistandard) tableaux of shape $\la$, that is, tableaux with entries in $\A$ such that 
(1) the entries in each row are weakly increasing from left to right, 
(2) the entries in each column are strictly increasing from top to bottom.

We denote by $CSST_\A(\la)$ the set of column semistandard tableaux of shape $\la$, that is, tableaux of shape $\la$ with entries in $\A$ which satisfy the condition (2). We also let $CST_\A(\la)$ be the set of $T\in CSST_\A(\la)$ such that all the entries in $T$ are distinct.

{Similarly, we let $RSST_\A(\la)$ be the set of row semistandard tableaux of shape $\la$, that is, tableaux of shape $\la$ with entries in $\mc{A}$ which satisfy the condition (1), and let $RST_\A(\la)$ be the set of $T\in RSST_\A(\la)$ such that all the entries in $T$ are distinct.}

\section{Affine RSK correspondence}\label{sec:aff RSK}

\subsection{Proper numberings and zig-zags}\label{subsec:zig-zags}

\begin{df}\label{df:zig-zag}
{\rm Let ${\bf z}=\{c_i\}_{i\in \Z}$ be an infinite collection of cells such that $c_{i+1}$ is adjacent to $c_i$ for $i \in \Z$.
\begin{itemize}
\item[(1)] We say that ${\bf z}$ is a {\em zig-zag} when it satisfies the following:

(i) $c_{i+1}$ is located to the east or north of $c_i$ for all $i$,

(ii) $c_{i+1}$ is located to the east of $c_i$ for $i\gg 0$,

(iii) $c_{i-1}$ is located to the south of $c_i$ for $i\ll 0$.
 
\item[(2)] A cell $c_i$ in a zig-zag ${\bf z}$ is called an {\em inner corner} if $c_{i-1}$ is located to the south and $c_{i+1}$ is located to the east of $c_i$, and called an {\em outer corner} if $c_{i-1}$ is located to the west and $c_{i+1}$ is located to the north of $c_i$. Note that a zig-zag ${\bf z}$ has at least one inner corner.
 
\item[(3)] A cell $c=(i_l,j_k)$ is called the {\em back-post corner} of a zig-zag ${\bf z}$  when $c_k=(i_k,j_k)$ and $c_l=(i_l,j_l)$ are the leftmost and the rightmost inner corners of ${\bf z}$, respectively.
 
\end{itemize}}
\end{df}

Let $d$ be a proper numbering on a non-zero $A\in \wh{\M}_{m\times n}$. For each $k\in\Z$, the level set $d^{-1}(k)$ forms a chain with respect to $\le_{\tt ne}$, by monotonicity of $d$ (Definition \ref{def:proper}). 

We associate a set of zig-zags $Z_d=\{{\bf z}_k\}_{k \in \Z}$ to $d$, where the inner corners of ${\bf z}_k$ are the set of maximal elements in $d^{-1}(k)$ with respect to $\ge_{\tt nw}$. It is straightforward to see that $\{ {\bf z}_k \}_{k \in \Z}$ satisfies
\begin{itemize}
    \item[(z.1)] the inner corners of each ${\bf z}_k$ are contained in ${\rm supp}(A)$,
    \item[(z.2)] ${\bf z}_k$'s are mutually disjoint and ${\rm supp}(A) \subseteq \bigsqcup_{k \in \Z} {\bf z}_k$, 
    \item[(z.3)] ${\bf z}_k$ is located to the southeast of ${\bf z}_{k-1}$ for $k\in\Z$ in the sense that
    \begin{equation}\label{eq:NW chain on zig-zag}
	\text{
	for each $c_1 \in {\bf z}_k$, there exists $c_2 \in {\bf z}_{k-1}$ such that $c_2 >_{\tt NW} c_1$.}
	\end{equation}
\end{itemize}
Conversely, a set of zig-zags $Z=\{{\bf z}_k\}_{k \in \Z}$ satisfying (z.1)-(z.3) determines a unique proper numbering $d^{Z}$ on $A$ given by
\begin{equation}\label{eq:numbering defined by zig-zags}
d^{Z}(c) = k \quad \text{if $c \in {\rm supp}(A)\cap {\bf z}_k$},
\end{equation}
whose associated {set of} zig-zags is $Z$.

\subsection{Matrix-ball construction}\label{subsec:mbc}
Now suppose that a non-zero $A\in \wh{\M}_{m\times n}$ is given and let
\begin{itemize}
\item[$\bullet$] $\{{\bf z}_k\}_{k\in\Z}$ : the set of zig-zags associated to $d_A^{\tt sw}$.
\end{itemize}
Note that $\tau({\bf z}_k)={\bf z}_{k+\ell}$ ($k\in \Z$) if $\ell$ is the period of $d_A^{\tt sw}$ (cf.~Lemma \ref{lem:period is unique}).

\begin{ex}\label{ex:zigzag of A}{\rm
Let $A$ and $d^{\tt sw}_A$ be as in Example \ref{ex:swc numbering of A}. The period of $d^{\tt sw}_A$ is 3. The zig-zags ${\bf z}_k$ corresponding to the level sets $(d_A^{\tt sw})^{-1}(k)$ for $k=1, 2, 3$ are given as red lines below.\vskip 2mm

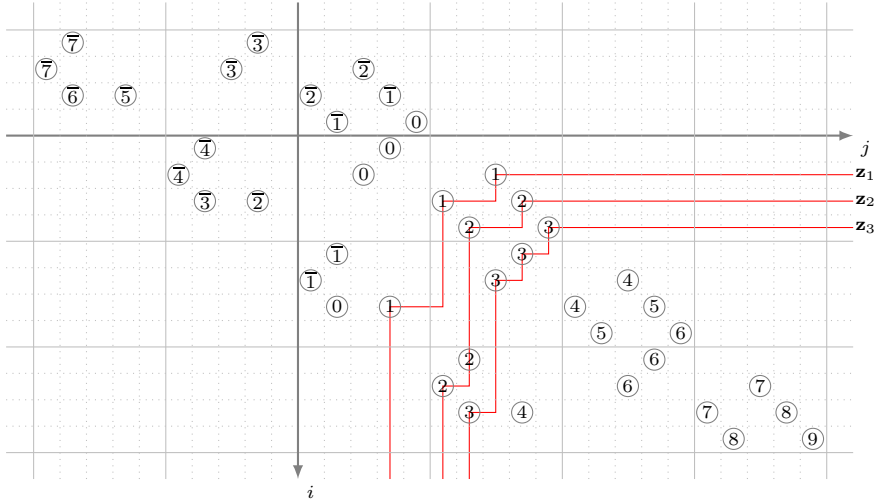
\begin{center}
\begin{tikzpicture}[every node/.style={font=\scriptsize}]
\matrix (M)[matrix of math nodes,nodes in empty cells,nodes={rectangle,minimum height=1.0em,minimum width=1.0em,inner sep=0pt,anchor=center,align=center}]
{
& &&&&& &&&&& &&&&& &&&&& &&&&& &&&&& &\\
& &$\ov{7}$&&&& &&&$\ov{3}$&& &&&&& &&&&& &&&&& &&&&& &\\
& $\ov{7}$&&&&& &&$\ov{3}$&&& &&$\ov{2}$&&& &&&&& &&&&& &&&&& &\\
& &$\ov{6}$&&$\ov{5}$&& &&&&& $\ov{2}$&&&$\ov{1}$&& &&&&& &&&&& &&&&& &\\
& &&&&& &&&&& &$\ov{1}$&&&0& &&&&& &&&&& &&&&& &\\
& &&&&& &$\ov{4}$&&&& &&&0&& &&&&& &&&&& &&&&& &\\
& &&&&& $\ov{4}$&&&&& &&0&&& &&1&&& &&&&& &&&&& &\\
& &&&&& &$\ov{3}$&&$\ov{2}$&& &&&&& 1&&&2&& &&&&& &&&&& &\\
& &&&&& &&&&& &&&&& &2&&&3& &&&&& &&&&& &\\
& &&&&& &&&&& &$\ov{1}$&&&& &&&3&& &&&&& &&&&& &\\
& &&&&& &&&&& $\ov{1}$&&&&& &&3&&& &&4&&& &&&&& &\\
& &&&&& &&&&& &0&&1&& &&&&& 4&&&5&& &&&&& &\\
& &&&&& &&&&& &&&&& &&&&& &5&&&6& &&&&& &\\
& &&&&& &&&&& &&&&& &2&&&& &&&6&& &&&&& &\\
& &&&&& &&&&& &&&&& 2&&&&& &&6&&& &&7&&& &\\
& &&&&& &&&&& &&&&& &3&&4&& &&&&& 7&&&8&& &\\
& &&&&& &&&&& &&&&& &&&&& &&&&& &8&&&9& &\\
& &&&&& &&&&& &&&&& &&&&& &&&&& &&&&& &\\
& &&&&& &&&&& &&&&& &&&&& &&&&& &&&&& &\\
};
\draw[solid,lightgray](M-1-1.south west)--(M-1-32.south east);
\draw[dotted,lightgray](M-2-1.south west)--(M-2-32.south east);
\draw[dotted,lightgray](M-3-1.south west)--(M-3-32.south east);
\draw[dotted,lightgray](M-4-1.south west)--(M-4-32.south east);
\draw[-latex,thick,gray](M-5-1.south west)--(M-5-32.south east);
\node (j) at (M-6-33) {$j$};
\node (j) at (M-7-33) {${\bf z}_1$};
\node (j) at (M-8-33) {${\bf z}_2$};
\node (j) at (M-9-33) {${\bf z}_3$};
\draw[dotted,lightgray](M-6-1.south west)--(M-6-32.south east);
\draw[dotted,lightgray](M-7-1.south west)--(M-7-32.south east);
\draw[dotted,lightgray](M-8-1.south west)--(M-8-32.south east);
\draw[solid,lightgray](M-9-1.south west)--(M-9-32.south east);
\draw[dotted,lightgray](M-10-1.south west)--(M-10-32.south east);
\draw[dotted,lightgray](M-11-1.south west)--(M-11-32.south east);
\draw[dotted,lightgray](M-12-1.south west)--(M-12-32.south east);
\draw[solid,lightgray](M-13-1.south west)--(M-13-32.south east);
\draw[dotted,lightgray](M-14-1.south west)--(M-14-32.south east);
\draw[dotted,lightgray](M-15-1.south west)--(M-15-32.south east);
\draw[dotted,lightgray](M-16-1.south west)--(M-16-32.south east);
\draw[solid,lightgray](M-17-1.south west)--(M-17-32.south east);
\draw[solid,lightgray](M-1-1.north east)--(M-18-1.south east);
\draw[dotted,lightgray](M-1-2.north east)--(M-18-2.south east);
\draw[dotted,lightgray](M-1-3.north east)--(M-18-3.south east);
\draw[dotted,lightgray](M-1-4.north east)--(M-18-4.south east);
\draw[dotted,lightgray](M-1-5.north east)--(M-18-5.south east);
\draw[solid,lightgray](M-1-6.north east)--(M-18-6.south east);
\draw[dotted,lightgray](M-1-7.north east)--(M-18-7.south east);
\draw[dotted,lightgray](M-1-8.north east)--(M-18-8.south east);
\draw[dotted,lightgray](M-1-9.north east)--(M-18-9.south east);
\draw[dotted,lightgray](M-1-10.north east)--(M-18-10.south east);
\draw[-latex,thick,gray](M-1-11.north east)--(M-18-11.south east);
\node (i) at (M-19-12) {$i$};
\draw[dotted,lightgray](M-1-12.north east)--(M-18-12.south east);
\draw[dotted,lightgray](M-1-13.north east)--(M-18-13.south east);
\draw[dotted,lightgray](M-1-14.north east)--(M-18-14.south east);
\draw[dotted,lightgray](M-1-15.north east)--(M-18-15.south east);
\draw[solid,lightgray](M-1-16.north east)--(M-18-16.south east);
\draw[dotted,lightgray](M-1-17.north east)--(M-18-17.south east);
\draw[dotted,lightgray](M-1-18.north east)--(M-18-18.south east);
\draw[dotted,lightgray](M-1-19.north east)--(M-18-19.south east);
\draw[dotted,lightgray](M-1-20.north east)--(M-18-20.south east);
\draw[solid,lightgray](M-1-21.north east)--(M-18-21.south east);
\draw[dotted,lightgray](M-1-22.north east)--(M-18-22.south east);
\draw[dotted,lightgray](M-1-23.north east)--(M-18-23.south east);
\draw[dotted,lightgray](M-1-24.north east)--(M-18-24.south east);
\draw[dotted,lightgray](M-1-25.north east)--(M-18-25.south east);
\draw[solid,lightgray](M-1-26.north east)--(M-18-26.south east);
\draw[dotted,lightgray](M-1-27.north east)--(M-18-27.south east);
\draw[dotted,lightgray](M-1-28.north east)--(M-18-28.south east);
\draw[dotted,lightgray](M-1-29.north east)--(M-18-29.south east);
\draw[dotted,lightgray](M-1-30.north east)--(M-18-30.south east);
\draw[solid,lightgray](M-1-31.north east)--(M-18-31.south east);
\begin{scope}[on background layer]
\draw[thin,red](M-18-15.south) -- (M-12-15.center) -| (M-8-17.center) -| (M-7-19.center)  -- (M-7-32.east);
\draw[thin,red](M-18-17.south) -- (M-15-17.center) -| (M-9-18.center) -| (M-8-20.center)  -- (M-8-32.east);
\draw[thin,red](M-18-18.south) -- (M-16-18.center) -| (M-11-19.center) -| (M-10-20.center) -| (M-9-21.center)  -- (M-9-32.east);
\end{scope}
\draw[gray](M-2-3)circle(0.4em);
\draw[gray](M-2-10)circle(0.4em);
\draw[gray](M-3-2)circle(0.4em);
\draw[gray](M-3-9)circle(0.4em);
\draw[gray](M-3-14)circle(0.4em);
\draw[gray](M-4-3)circle(0.4em);
\draw[gray](M-4-5)circle(0.4em);
\draw[gray](M-4-12)circle(0.4em);
\draw[gray](M-4-15)circle(0.4em);
\draw[gray](M-5-13)circle(0.4em);
\draw[gray](M-5-16)circle(0.4em);
\draw[gray](M-6-8)circle(0.4em);
\draw[gray](M-6-15)circle(0.4em);
\draw[gray](M-7-7)circle(0.4em);
\draw[gray](M-7-14)circle(0.4em);
\draw[gray](M-7-19)circle(0.4em);
\draw[gray](M-8-8)circle(0.4em);
\draw[gray](M-8-10)circle(0.4em);
\draw[gray](M-8-17)circle(0.4em);
\draw[gray](M-8-20)circle(0.4em);
\draw[gray](M-9-18)circle(0.4em);
\draw[gray](M-9-21)circle(0.4em);
\draw[gray](M-10-13)circle(0.4em);
\draw[gray](M-10-20)circle(0.4em);
\draw[gray](M-11-12)circle(0.4em);
\draw[gray](M-11-19)circle(0.4em);
\draw[gray](M-11-24)circle(0.4em);
\draw[gray](M-12-13)circle(0.4em);
\draw[gray](M-12-15)circle(0.4em);
\draw[gray](M-12-22)circle(0.4em);
\draw[gray](M-12-25)circle(0.4em);
\draw[gray](M-13-23)circle(0.4em);
\draw[gray](M-13-26)circle(0.4em);
\draw[gray](M-14-18)circle(0.4em);
\draw[gray](M-14-25)circle(0.4em);
\draw[gray](M-15-17)circle(0.4em);
\draw[gray](M-15-24)circle(0.4em);
\draw[gray](M-15-29)circle(0.4em);
\draw[gray](M-16-18)circle(0.4em);
\draw[gray](M-16-20)circle(0.4em);
\draw[gray](M-16-27)circle(0.4em);
\draw[gray](M-16-30)circle(0.4em);
\draw[gray](M-17-28)circle(0.4em);
\draw[gray](M-17-31)circle(0.4em);
\end{tikzpicture}
\end{center}

}
\end{ex}

Let
\begin{itemize}
\item[$\bullet$] $A^\flat$ : the matrix in $\wh{\M}_{m\times n}$ obtained from $A$ by 

\subitem (i) subtracting one at the inner corners in ${\bf z}_k$ ($k\in\Z$),     

\subitem (ii) adding one at the outer corners in ${\bf z}_k$ ($k\in\Z$),

\item[$\bullet$] $A^{(t)}$ ($t\ge 0$) : the matrices in $\wh{\M}_{m\times n}$ defined inductively as follows:
\begin{equation*}
A^{(0)}=A,\quad A^{(t)}=\left( A^{(t-1)} \right)^\flat.
\end{equation*}
\end{itemize}
We have $A^{(t)}={\mathbb O}$ for a sufficiently large $t$, since the number of outer corners in each zig-zag ${\bf z}_k$ is strictly smaller than the number of inner corners. Let $s \geq 1$ be the integer such that
\begin{equation*}
A^{(s-1)}\neq {\mathbb O},\quad A^{(s)} ={\mathbb O}
\end{equation*}
For $1\le t\le s$, we consider 
\begin{itemize}
\item[$\bullet$] $\{{\bf z}^{(t)}_k\}_{k\in \Z}$ : the set of zig-zags associated to $d^{\tt sw}_{A^{(t-1)}}$,

\item[$\bullet$] ${\bf s}^{(t)}$ : the stream consisting of the back-post corners of ${\bf z}^{(t)}_k$  $(k\in \Z)$,

\item[$\bullet$] $\mu_t$ : the flow of ${\bf s}^{(t)}$, equivalently, the width of $A^{(t-1)}$,

\item[$\bullet$] $(\ba_t,\bb_t,\rho_t)$ : the defining data of ${\bf s}^{(t)}$ $(1\le t\le s)$.
\end{itemize}

\begin{lem}[cf.~{{\cite[Lemma 14.9]{CPY}}}]\label{lem:width decreasing}
The flows of ${\bf s}^{(1)},\dots,{\bf s}^{(s)}$ {are wealy decreasing positive integers}, that is, $\mu_1\ge \dots \ge \mu_{s}>0$. 
\end{lem}
\pf Let $\ell$ and $\ell'$ be the widths of $A$ and $A^\flat$, respectively. 
Consider the zig-zags ${\bf z}_k$  $(k\in \Z)$ associated to $d^{\tt sw}_A$. 
By definition, ${\rm supp}(A^\flat)$ is contained in the union of the ${\bf z}_k$'s. Suppose that $\ell < \ell'$. 
Let $C=\{c_i \}_{i \in \Z}$ be a channel of $A^\flat$. 
Note that each cell $c_i$ belongs to a different zig-zag ${\bf z}_k$.
But the width of $A$ is $\ell$, which implies that at least two cells among $\ell'$ cells $c_1 >_{\tt NW} c_2 >_{\tt NW} \dots >_{\tt NW} c_{l'} = \tau(c_0)$ {belong} to the same zig-zag ${\bf z}_k$ for some $k$. This is a contradiction.
\qed
\newline

Now we let
\begin{itemize}
\item[$\bullet$] $\la=\mu'$ : the conjugate partition of $\mu=(\mu_1,\dots,\mu_s)\in \cP_{s}$,

\item[$\bullet$] $P_0$ : the tableau of shape $\la$, whose $t$-th column from the left is ${\bf a}_{t}$ ($1\le t\le s$),

\item[$\bullet$] $Q_0$ : the tableau of shape $\la$, whose $t$-th column from the left is  ${\bf b}_{t}$ ($1\le t\le s$), 

\item[$\bullet$] $\rho=(\rho_1,\dots,\rho_s)\in \Z^s$.
 
\end{itemize}
Here we understand that the entries of the $t$-th column in $P_0$ and $Q_0$ are given by those in $\ba_t$ and $\bb_t$ respectively. Note that $\la\in \cP_m\cap\cP_n$, and $P_0, Q_0$ are column semistandard tableaux.

Summarizing, we have a map

\begin{equation}\label{eq:Aff RSK 0}
\xymatrixcolsep{3pc}\xymatrixrowsep{0.5pc}\xymatrix{
\kappa_0 :\ \wh{\M}_{m\times n} \ \ar@{->}[r]  & \ \
\displaystyle\bigsqcup_{\la\in \cP_m\cap\cP_n} CSST_{[m]}(\la)\times CSST_{[n]}(\la) \times \Z^{\la_1}  \\
\quad A \  \ar@{|->}[r]  &\ \  (P_0, Q_0, \rho)},
\end{equation}
where we assume $\kappa_0({\mathbb O}) = (\emptyset, \emptyset, (0, \dots, 0))$ with $\emptyset$ the empty tableau.
\begin{ex}\label{ex:kappa_0 of A}{\rm
Let $A$ be the generalized affine permutation in Example \ref{ex:A}, {i.e.,} \vskip 2mm

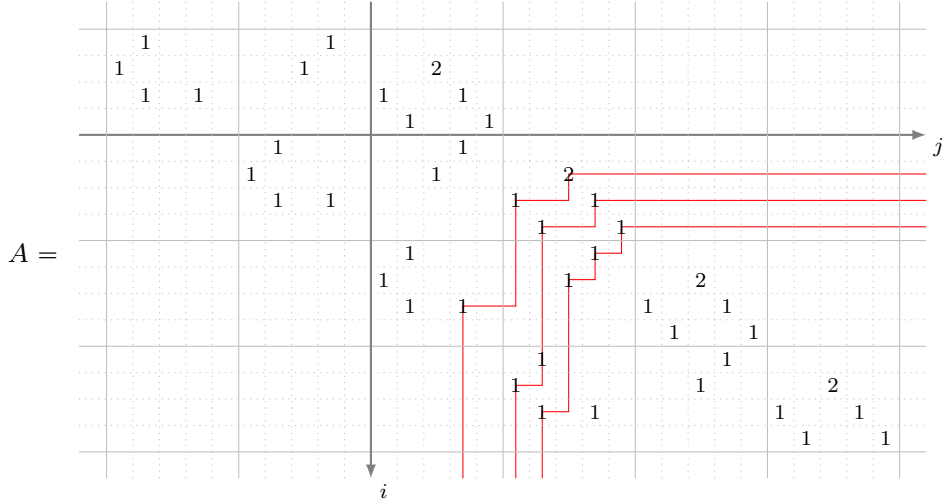
\begin{center}
\begin{tikzpicture}[every node/.style={font=\footnotesize,scale=1}]
\matrix (M)[matrix of math nodes,nodes in empty cells,nodes={rectangle,minimum height=1.0em,minimum width=1.0em,inner sep=0pt,anchor=center,align=center}]
{
& &&&&& &&&&& &&&&& &&&&& &&&&& &&&&& &\\
& &1&&&& &&&1&& &&&&& &&&&& &&&&& &&&&& &\\
& 1&&&&& &&1&&& &&2&&& &&&&& &&&&& &&&&& &\\
& &1&&1&& &&&&& 1&&&1&& &&&&& &&&&& &&&&& &\\
& &&&&& &&&&& &1&&&1& &&&&& &&&&& &&&&& &\\
& &&&&& &1&&&& &&&1&& &&&&& &&&&& &&&&& &\\
& &&&&& 1&&&&& &&1&&& &&2&&& &&&&& &&&&& &\\
& &&&&& &1&&1&& &&&&& 1&&&1&& &&&&& &&&&& &\\
& &&&&& &&&&& &&&&& &1&&&1& &&&&& &&&&& &\\
& &&&&& &&&&& &1&&&& &&&1&& &&&&& &&&&& &\\
& &&&&& &&&&& 1&&&&& &&1&&& &&2&&& &&&&& &\\
& &&&&& &&&&& &1&&1&& &&&&& 1&&&1&& &&&&& &\\
& &&&&& &&&&& &&&&& &&&&& &1&&&1& &&&&& &\\
& &&&&& &&&&& &&&&& &1&&&& &&&1&& &&&&& &\\
& &&&&& &&&&& &&&&& 1&&&&& &&1&&& &&2&&& &\\
& &&&&& &&&&& &&&&& &1&&1&& &&&&& 1&&&1&& &\\
& &&&&& &&&&& &&&&& &&&&& &&&&& &1&&&1& &\\
& &&&&& &&&&& &&&&& &&&&& &&&&& &&&&& &\\
& &&&&& &&&&& &&&&& &&&&& &&&&& &&&&& &\\
};
\draw[solid,lightgray](M-1-1.south west)--(M-1-32.south east);
\draw[dotted,lightgray](M-2-1.south west)--(M-2-32.south east);
\draw[dotted,lightgray](M-3-1.south west)--(M-3-32.south east);
\draw[dotted,lightgray](M-4-1.south west)--(M-4-32.south east);
\draw[-latex,thick,gray](M-5-1.south west)--(M-5-32.south east);
\node (j) at (M-6-33) {$j$};
\draw[dotted,lightgray](M-6-1.south west)--(M-6-32.south east);
\draw[dotted,lightgray](M-7-1.south west)--(M-7-32.south east);
\draw[dotted,lightgray](M-8-1.south west)--(M-8-32.south east);
\draw[solid,lightgray](M-9-1.south west)--(M-9-32.south east);
\draw[dotted,lightgray](M-10-1.south west)--(M-10-32.south east);
\draw[dotted,lightgray](M-11-1.south west)--(M-11-32.south east);
\draw[dotted,lightgray](M-12-1.south west)--(M-12-32.south east);
\draw[solid,lightgray](M-13-1.south west)--(M-13-32.south east);
\draw[dotted,lightgray](M-14-1.south west)--(M-14-32.south east);
\draw[dotted,lightgray](M-15-1.south west)--(M-15-32.south east);
\draw[dotted,lightgray](M-16-1.south west)--(M-16-32.south east);
\draw[solid,lightgray](M-17-1.south west)--(M-17-32.south east);
\draw[solid,lightgray](M-1-1.north east)--(M-18-1.south east);
\draw[dotted,lightgray](M-1-2.north east)--(M-18-2.south east);
\draw[dotted,lightgray](M-1-3.north east)--(M-18-3.south east);
\draw[dotted,lightgray](M-1-4.north east)--(M-18-4.south east);
\draw[dotted,lightgray](M-1-5.north east)--(M-18-5.south east);
\draw[solid,lightgray](M-1-6.north east)--(M-18-6.south east);
\draw[dotted,lightgray](M-1-7.north east)--(M-18-7.south east);
\draw[dotted,lightgray](M-1-8.north east)--(M-18-8.south east);
\draw[dotted,lightgray](M-1-9.north east)--(M-18-9.south east);
\draw[dotted,lightgray](M-1-10.north east)--(M-18-10.south east);
\draw[-latex,thick,gray](M-1-11.north east)--(M-18-11.south east);
\node (i) at (M-19-12) {$i$};
\draw[dotted,lightgray](M-1-12.north east)--(M-18-12.south east);
\draw[dotted,lightgray](M-1-13.north east)--(M-18-13.south east);
\draw[dotted,lightgray](M-1-14.north east)--(M-18-14.south east);
\draw[dotted,lightgray](M-1-15.north east)--(M-18-15.south east);
\draw[solid,lightgray](M-1-16.north east)--(M-18-16.south east);
\draw[dotted,lightgray](M-1-17.north east)--(M-18-17.south east);
\draw[dotted,lightgray](M-1-18.north east)--(M-18-18.south east);
\draw[dotted,lightgray](M-1-19.north east)--(M-18-19.south east);
\draw[dotted,lightgray](M-1-20.north east)--(M-18-20.south east);
\draw[solid,lightgray](M-1-21.north east)--(M-18-21.south east);
\draw[dotted,lightgray](M-1-22.north east)--(M-18-22.south east);
\draw[dotted,lightgray](M-1-23.north east)--(M-18-23.south east);
\draw[dotted,lightgray](M-1-24.north east)--(M-18-24.south east);
\draw[dotted,lightgray](M-1-25.north east)--(M-18-25.south east);
\draw[solid,lightgray](M-1-26.north east)--(M-18-26.south east);
\draw[dotted,lightgray](M-1-27.north east)--(M-18-27.south east);
\draw[dotted,lightgray](M-1-28.north east)--(M-18-28.south east);
\draw[dotted,lightgray](M-1-29.north east)--(M-18-29.south east);
\draw[dotted,lightgray](M-1-30.north east)--(M-18-30.south east);
\draw[solid,lightgray](M-1-31.north east)--(M-18-31.south east);
\begin{scope}[on background layer]
\draw[thin,red](M-18-15.south) -- (M-12-15.center) -| (M-8-17.center) -| (M-7-19.center)  -- (M-7-32.east);
\draw[thin,red](M-18-17.south) -- (M-15-17.center) -| (M-9-18.center) -| (M-8-20.center)  -- (M-8-32.east);
\draw[thin,red](M-18-18.south) -- (M-16-18.center) -| (M-11-19.center) -| (M-10-20.center) -| (M-9-21.center)  -- (M-9-32.east);
\end{scope}
\node[anchor=east,scale=1.25] at (M.west) {$A=$};
\end{tikzpicture}
\end{center}
where the red lines denote the zig-zags in Example \ref{ex:zigzag of A}.
The stream consisting of the back-post corners of ${\bf z}_k$ ($k\in\Z$) is
$${\bf s}^{(1)} = \{ \cdots >_{\tt NW} (2, 4) >_{\tt NW} (3, 6) >_{\tt NW} (4, 7) >_{\tt NW} \cdots \}.$$
By subtracting one at each inner corner and adding one at each outer corner of ${\bf z}_k$, we obtain $A^{\flat}$ as follows:
\vskip 2mm

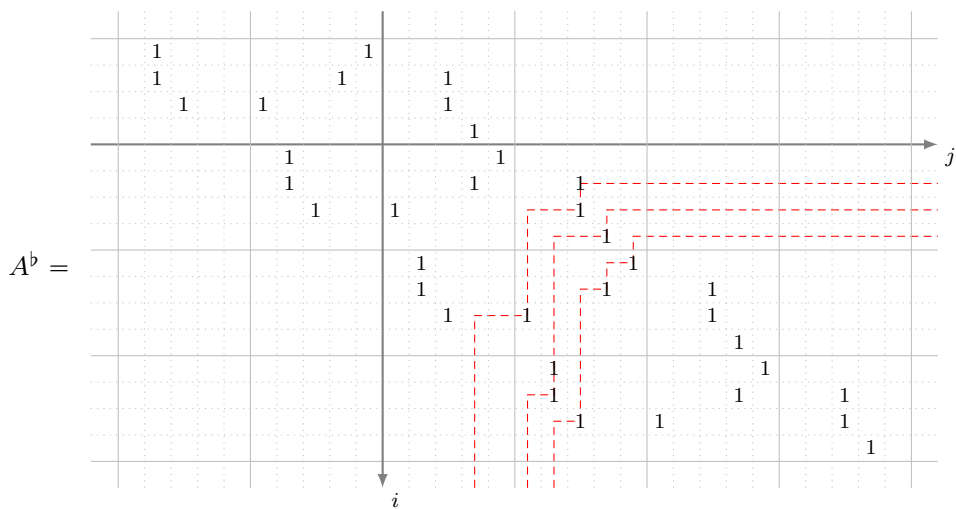
\begin{center}
\begin{tikzpicture}[every node/.style={font=\footnotesize,scale=1}]
\matrix (M)[matrix of math nodes,nodes in empty cells,nodes={rectangle,minimum height=1.0em,minimum width=1.0em,inner sep=0pt,anchor=center,align=center}]
{
& &&&&& &&&&& &&&&& &&&&& &&&&& &&&&& &\\
& &1&&&& &&&&1& &&&&& &&&&& &&&&& &&&&& &\\
& &1&&&& &&&1&& &&1&&& &&&&& &&&&& &&&&& &\\
& &&1&&& 1&&&&& &&1&&& &&&&& &&&&& &&&&& &\\
& &&&&& &&&&& &&&1&& &&&&& &&&&& &&&&& &\\
& &&&&& &1&&&& &&&&1& &&&&& &&&&& &&&&& &\\
& &&&&& &1&&&& &&&1&& &&1&&& &&&&& &&&&& &\\
& &&&&& &&1&&& 1&&&&& &&1&&& &&&&& &&&&& &\\
& &&&&& &&&&& &&&&& &&&1&& &&&&& &&&&& &\\
& &&&&& &&&&& &1&&&& &&&&1& &&&&& &&&&& &\\
& &&&&& &&&&& &1&&&& &&&1&& &&1&&& &&&&& &\\
& &&&&& &&&&& &&1&&& 1&&&&& &&1&&& &&&&& &\\
& &&&&& &&&&& &&&&& &&&&& &&&1&& &&&&& &\\
& &&&&& &&&&& &&&&& &1&&&& &&&&1& &&&&& &\\
& &&&&& &&&&& &&&&& &1&&&& &&&1&& &&1&&& &\\
& &&&&& &&&&& &&&&& &&1&&& 1&&&&& &&1&&& &\\
& &&&&& &&&&& &&&&& &&&&& &&&&& &&&1&& &\\
& &&&&& &&&&& &&&&& &&&&& &&&&& &&&&& &\\
& &&&&& &&&&& &&&&& &&&&& &&&&& &&&&& &\\
};
\draw[solid,lightgray](M-1-1.south west)--(M-1-32.south east);
\draw[dotted,lightgray](M-2-1.south west)--(M-2-32.south east);
\draw[dotted,lightgray](M-3-1.south west)--(M-3-32.south east);
\draw[dotted,lightgray](M-4-1.south west)--(M-4-32.south east);
\draw[-latex,thick,gray](M-5-1.south west)--(M-5-32.south east);
\node (j) at (M-6-33) {$j$};
\draw[dotted,lightgray](M-6-1.south west)--(M-6-32.south east);
\draw[dotted,lightgray](M-7-1.south west)--(M-7-32.south east);
\draw[dotted,lightgray](M-8-1.south west)--(M-8-32.south east);
\draw[solid,lightgray](M-9-1.south west)--(M-9-32.south east);
\draw[dotted,lightgray](M-10-1.south west)--(M-10-32.south east);
\draw[dotted,lightgray](M-11-1.south west)--(M-11-32.south east);
\draw[dotted,lightgray](M-12-1.south west)--(M-12-32.south east);
\draw[solid,lightgray](M-13-1.south west)--(M-13-32.south east);
\draw[dotted,lightgray](M-14-1.south west)--(M-14-32.south east);
\draw[dotted,lightgray](M-15-1.south west)--(M-15-32.south east);
\draw[dotted,lightgray](M-16-1.south west)--(M-16-32.south east);
\draw[solid,lightgray](M-17-1.south west)--(M-17-32.south east);
\draw[solid,lightgray](M-1-1.north east)--(M-18-1.south east);
\draw[dotted,lightgray](M-1-2.north east)--(M-18-2.south east);
\draw[dotted,lightgray](M-1-3.north east)--(M-18-3.south east);
\draw[dotted,lightgray](M-1-4.north east)--(M-18-4.south east);
\draw[dotted,lightgray](M-1-5.north east)--(M-18-5.south east);
\draw[solid,lightgray](M-1-6.north east)--(M-18-6.south east);
\draw[dotted,lightgray](M-1-7.north east)--(M-18-7.south east);
\draw[dotted,lightgray](M-1-8.north east)--(M-18-8.south east);
\draw[dotted,lightgray](M-1-9.north east)--(M-18-9.south east);
\draw[dotted,lightgray](M-1-10.north east)--(M-18-10.south east);
\draw[-latex,thick,gray](M-1-11.north east)--(M-18-11.south east);
\node (i) at (M-19-12) {$i$};
\draw[dotted,lightgray](M-1-12.north east)--(M-18-12.south east);
\draw[dotted,lightgray](M-1-13.north east)--(M-18-13.south east);
\draw[dotted,lightgray](M-1-14.north east)--(M-18-14.south east);
\draw[dotted,lightgray](M-1-15.north east)--(M-18-15.south east);
\draw[solid,lightgray](M-1-16.north east)--(M-18-16.south east);
\draw[dotted,lightgray](M-1-17.north east)--(M-18-17.south east);
\draw[dotted,lightgray](M-1-18.north east)--(M-18-18.south east);
\draw[dotted,lightgray](M-1-19.north east)--(M-18-19.south east);
\draw[dotted,lightgray](M-1-20.north east)--(M-18-20.south east);
\draw[solid,lightgray](M-1-21.north east)--(M-18-21.south east);
\draw[dotted,lightgray](M-1-22.north east)--(M-18-22.south east);
\draw[dotted,lightgray](M-1-23.north east)--(M-18-23.south east);
\draw[dotted,lightgray](M-1-24.north east)--(M-18-24.south east);
\draw[dotted,lightgray](M-1-25.north east)--(M-18-25.south east);
\draw[solid,lightgray](M-1-26.north east)--(M-18-26.south east);
\draw[dotted,lightgray](M-1-27.north east)--(M-18-27.south east);
\draw[dotted,lightgray](M-1-28.north east)--(M-18-28.south east);
\draw[dotted,lightgray](M-1-29.north east)--(M-18-29.south east);
\draw[dotted,lightgray](M-1-30.north east)--(M-18-30.south east);
\draw[solid,lightgray](M-1-31.north east)--(M-18-31.south east);
\begin{scope}[on background layer]
\draw[densely dashed,red](M-18-15.south) -- (M-12-15.center) -| (M-8-17.center) -| (M-7-19.center)  -- (M-7-32.east);
\draw[densely dashed,red](M-18-17.south) -- (M-15-17.center) -| (M-9-18.center) -| (M-8-20.center)  -- (M-8-32.east);
\draw[densely dashed,red](M-18-18.south) -- (M-16-18.center) -| (M-11-19.center) -| (M-10-20.center) -| (M-9-21.center)  -- (M-9-32.east);
\end{scope}
\node[anchor=east,scale=1.25] at (M.west) {$A^{\flat}=$};
\end{tikzpicture}
\end{center}
Repeating this process, we see that \vskip 2mm

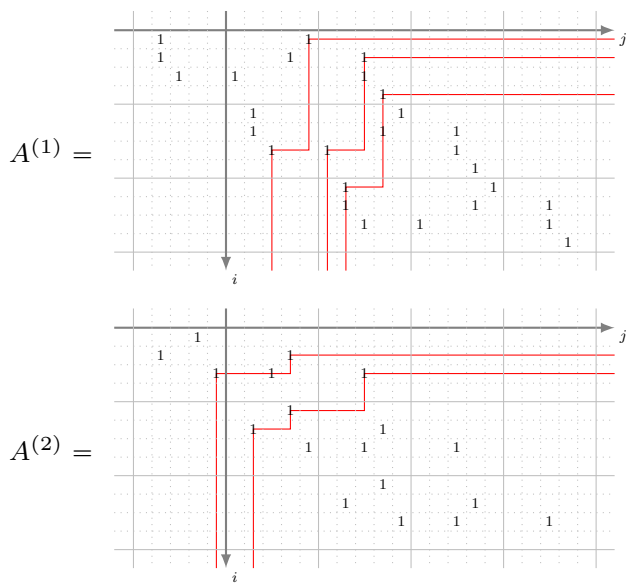
\begin{center}
\begin{tikzpicture}[every node/.style={font=\footnotesize,scale=0.70}]
\matrix (M)[matrix of math nodes,nodes in empty cells,nodes={rectangle,minimum height=1.0em,minimum width=1.0em,inner sep=0pt,anchor=center,align=center}]
{
& &&&&& &&&&& &&&&& &&&&& &&&&& &\\
& &1&&&& &&&&1& &&&&& &&&&& &&&&& &\\
& &1&&&& &&&1&& &&1&&& &&&&& &&&&& &\\
& &&1&&& 1&&&&& &&1&&& &&&&& &&&&& &\\
& &&&&& &&&&& &&&1&& &&&&& &&&&& &\\
& &&&&& &1&&&& &&&&1& &&&&& &&&&& &\\
& &&&&& &1&&&& &&&1&& &&1&&& &&&&& &\\
& &&&&& &&1&&& 1&&&&& &&1&&& &&&&& &\\
& &&&&& &&&&& &&&&& &&&1&& &&&&& &\\
& &&&&& &&&&& &1&&&& &&&&1& &&&&& &\\
& &&&&& &&&&& &1&&&& &&&1&& &&1&&& &\\
& &&&&& &&&&& &&1&&& 1&&&&& &&1&&& &\\
& &&&&& &&&&& &&&&& &&&&& &&&1&& &\\
& &&&&& &&&&& &&&&& &&&&& &&&&& &\\
& &&&&& &&&&& &&&&& &&&&& &&&&& &\\
};
\draw[-latex,thick,gray](M-1-1.south west)--(M-1-27.south east);
\node (j) at (M-2-28) {$j$};
\draw[dotted,lightgray](M-2-1.south west)--(M-2-27.south east);
\draw[dotted,lightgray](M-3-1.south west)--(M-3-27.south east);
\draw[dotted,lightgray](M-4-1.south west)--(M-4-27.south east);
\draw[solid,lightgray](M-5-1.south west)--(M-5-27.south east);
\draw[dotted,lightgray](M-6-1.south west)--(M-6-27.south east);
\draw[dotted,lightgray](M-7-1.south west)--(M-7-27.south east);
\draw[dotted,lightgray](M-8-1.south west)--(M-8-27.south east);
\draw[solid,lightgray](M-9-1.south west)--(M-9-27.south east);
\draw[dotted,lightgray](M-10-1.south west)--(M-10-27.south east);
\draw[dotted,lightgray](M-11-1.south west)--(M-11-27.south east);
\draw[dotted,lightgray](M-12-1.south west)--(M-12-27.south east);
\draw[solid,lightgray](M-13-1.south west)--(M-13-27.south east);
\draw[solid,lightgray](M-1-1.north east)--(M-14-1.south east);
\draw[dotted,lightgray](M-1-2.north east)--(M-14-2.south east);
\draw[dotted,lightgray](M-1-3.north east)--(M-14-3.south east);
\draw[dotted,lightgray](M-1-4.north east)--(M-14-4.south east);
\draw[dotted,lightgray](M-1-5.north east)--(M-14-5.south east);
\draw[-latex,thick,gray](M-1-6.north east)--(M-14-6.south east);
\node (i) at (M-15-7) {$i$};
\draw[dotted,lightgray](M-1-7.north east)--(M-14-7.south east);
\draw[dotted,lightgray](M-1-8.north east)--(M-14-8.south east);
\draw[dotted,lightgray](M-1-9.north east)--(M-14-9.south east);
\draw[dotted,lightgray](M-1-10.north east)--(M-14-10.south east);
\draw[solid,lightgray](M-1-11.north east)--(M-14-11.south east);
\draw[dotted,lightgray](M-1-12.north east)--(M-14-12.south east);
\draw[dotted,lightgray](M-1-13.north east)--(M-14-13.south east);
\draw[dotted,lightgray](M-1-14.north east)--(M-14-14.south east);
\draw[dotted,lightgray](M-1-15.north east)--(M-14-15.south east);
\draw[solid,lightgray](M-1-16.north east)--(M-14-16.south east);
\draw[dotted,lightgray](M-1-17.north east)--(M-14-17.south east);
\draw[dotted,lightgray](M-1-18.north east)--(M-14-18.south east);
\draw[dotted,lightgray](M-1-19.north east)--(M-14-19.south east);
\draw[dotted,lightgray](M-1-20.north east)--(M-14-20.south east);
\draw[solid,lightgray](M-1-21.north east)--(M-14-21.south east);
\draw[dotted,lightgray](M-1-22.north east)--(M-14-22.south east);
\draw[dotted,lightgray](M-1-23.north east)--(M-14-23.south east);
\draw[dotted,lightgray](M-1-24.north east)--(M-14-24.south east);
\draw[dotted,lightgray](M-1-25.north east)--(M-14-25.south east);
\draw[solid,lightgray](M-1-26.north east)--(M-14-26.south east);
\begin{scope}[on background layer]
\draw[thin,red](M-14-9.south) -- (M-8-9.center) -| (M-2-11.center)  -- (M-2-27.east);
\draw[thin,red](M-14-12.south) -- (M-8-12.center) -| (M-3-14.center)  -- (M-3-27.east);
\draw[thin,red](M-14-13.south) -- (M-10-13.center) -| (M-5-15.center)  -- (M-5-27.east);
\end{scope}
\node[anchor=east,scale=1.8] at (M.west) {$A^{(1)}=$};
\end{tikzpicture}

\begin{tikzpicture}[every node/.style={font=\footnotesize,scale=0.70}]
\matrix (M)[matrix of math nodes,nodes in empty cells,nodes={rectangle,minimum height=1.0em,minimum width=1.0em,inner sep=0pt,anchor=center,align=center}]
{
& &&&&& &&&&& &&&&& &&&&& &&&&& &\\
& &&&1&& &&&&& &&&&& &&&&& &&&&& &\\
& &1&&&& &&&1&& &&&&& &&&&& &&&&& &\\
& &&&&1& &&1&&& &&1&&& &&&&& &&&&& &\\
& &&&&& &&&&& &&&&& &&&&& &&&&& &\\
& &&&&& &&&1&& &&&&& &&&&& &&&&& &\\
& &&&&& &1&&&& &&&1&& &&&&& &&&&& &\\
& &&&&& &&&&1& &&1&&& &&1&&& &&&&& &\\
& &&&&& &&&&& &&&&& &&&&& &&&&& &\\
& &&&&& &&&&& &&&1&& &&&&& &&&&& &\\
& &&&&& &&&&& &1&&&& &&&1&& &&&&& &\\
& &&&&& &&&&& &&&&1& &&1&&& &&1&&& &\\
& &&&&& &&&&& &&&&& &&&&& &&&&& &\\
& &&&&& &&&&& &&&&& &&&&& &&&&& &\\
& &&&&& &&&&& &&&&& &&&&& &&&&& &\\
};
\draw[-latex,thick,gray](M-1-1.south west)--(M-1-27.south east);
\node (j) at (M-2-28) {$j$};
\draw[dotted,lightgray](M-2-1.south west)--(M-2-27.south east);
\draw[dotted,lightgray](M-3-1.south west)--(M-3-27.south east);
\draw[dotted,lightgray](M-4-1.south west)--(M-4-27.south east);
\draw[solid,lightgray](M-5-1.south west)--(M-5-27.south east);
\draw[dotted,lightgray](M-6-1.south west)--(M-6-27.south east);
\draw[dotted,lightgray](M-7-1.south west)--(M-7-27.south east);
\draw[dotted,lightgray](M-8-1.south west)--(M-8-27.south east);
\draw[solid,lightgray](M-9-1.south west)--(M-9-27.south east);
\draw[dotted,lightgray](M-10-1.south west)--(M-10-27.south east);
\draw[dotted,lightgray](M-11-1.south west)--(M-11-27.south east);
\draw[dotted,lightgray](M-12-1.south west)--(M-12-27.south east);
\draw[solid,lightgray](M-13-1.south west)--(M-13-27.south east);
\draw[solid,lightgray](M-1-1.north east)--(M-14-1.south east);
\draw[dotted,lightgray](M-1-2.north east)--(M-14-2.south east);
\draw[dotted,lightgray](M-1-3.north east)--(M-14-3.south east);
\draw[dotted,lightgray](M-1-4.north east)--(M-14-4.south east);
\draw[dotted,lightgray](M-1-5.north east)--(M-14-5.south east);
\draw[-latex,thick,gray](M-1-6.north east)--(M-14-6.south east);
\node (i) at (M-15-7) {$i$};
\draw[dotted,lightgray](M-1-7.north east)--(M-14-7.south east);
\draw[dotted,lightgray](M-1-8.north east)--(M-14-8.south east);
\draw[dotted,lightgray](M-1-9.north east)--(M-14-9.south east);
\draw[dotted,lightgray](M-1-10.north east)--(M-14-10.south east);
\draw[solid,lightgray](M-1-11.north east)--(M-14-11.south east);
\draw[dotted,lightgray](M-1-12.north east)--(M-14-12.south east);
\draw[dotted,lightgray](M-1-13.north east)--(M-14-13.south east);
\draw[dotted,lightgray](M-1-14.north east)--(M-14-14.south east);
\draw[dotted,lightgray](M-1-15.north east)--(M-14-15.south east);
\draw[solid,lightgray](M-1-16.north east)--(M-14-16.south east);
\draw[dotted,lightgray](M-1-17.north east)--(M-14-17.south east);
\draw[dotted,lightgray](M-1-18.north east)--(M-14-18.south east);
\draw[dotted,lightgray](M-1-19.north east)--(M-14-19.south east);
\draw[dotted,lightgray](M-1-20.north east)--(M-14-20.south east);
\draw[solid,lightgray](M-1-21.north east)--(M-14-21.south east);
\draw[dotted,lightgray](M-1-22.north east)--(M-14-22.south east);
\draw[dotted,lightgray](M-1-23.north east)--(M-14-23.south east);
\draw[dotted,lightgray](M-1-24.north east)--(M-14-24.south east);
\draw[dotted,lightgray](M-1-25.north east)--(M-14-25.south east);
\draw[solid,lightgray](M-1-26.north east)--(M-14-26.south east);
\begin{scope}[on background layer]
\draw[thin,red](M-14-6.south) -- (M-4-6.center) -| (M-3-10.center)  -- (M-3-27.east);
\draw[thin,red](M-14-8.south) -- (M-7-8.center) -| (M-6-10.center) -| (M-4-14.center)  -- (M-4-27.east);
\end{scope}
\node[anchor=east,scale=1.8] at (M.west) {$A^{(2)}=$};
\end{tikzpicture} 
\end{center}

\begin{center}
\begin{tikzpicture}[every node/.style={font=\footnotesize,scale=0.70}]
\matrix (M)[matrix of math nodes,nodes in empty cells,nodes={rectangle,minimum height=1.0em,minimum width=1.0em,inner sep=0pt,anchor=center,align=center}]
{
& &&&&& &&&&& &&&&& &&&&& &&&&& &\\
& &&&&& &&1&&& &&&&& &&&&& &&&&& &\\
& &&&1&& &&&&& &&&&& &&&&& &&&&& &\\
& &&&&& &&1&1&& &&&&& &&&&& &&&&& &\\
& &&&&& &&&&& &&&&& &&&&& &&&&& &\\
& &&&&& &&&&& &&1&&& &&&&& &&&&& &\\
& &&&&& &&&1&& &&&&& &&&&& &&&&& &\\
& &&&&& &&&&& &&1&1&& &&&&& &&&&& &\\
& &&&&& &&&&& &&&&& &&&&& &&&&& &\\
& &&&&& &&&&& &&&&& &&1&&& &&&&& &\\
& &&&&& &&&&& &&&1&& &&&&& &&&&& &\\
& &&&&& &&&&& &&&&& &&1&1&& &&&&& &\\
& &&&&& &&&&& &&&&& &&&&& &&&&& &\\
& &&&&& &&&&& &&&&& &&&&& &&&&& &\\
& &&&&& &&&&& &&&&& &&&&& &&&&& &\\
};
\draw[-latex,thick,gray](M-1-1.south west)--(M-1-27.south east);
\node (j) at (M-2-28) {$j$};
\draw[dotted,lightgray](M-2-1.south west)--(M-2-27.south east);
\draw[dotted,lightgray](M-3-1.south west)--(M-3-27.south east);
\draw[dotted,lightgray](M-4-1.south west)--(M-4-27.south east);
\draw[solid,lightgray](M-5-1.south west)--(M-5-27.south east);
\draw[dotted,lightgray](M-6-1.south west)--(M-6-27.south east);
\draw[dotted,lightgray](M-7-1.south west)--(M-7-27.south east);
\draw[dotted,lightgray](M-8-1.south west)--(M-8-27.south east);
\draw[solid,lightgray](M-9-1.south west)--(M-9-27.south east);
\draw[dotted,lightgray](M-10-1.south west)--(M-10-27.south east);
\draw[dotted,lightgray](M-11-1.south west)--(M-11-27.south east);
\draw[dotted,lightgray](M-12-1.south west)--(M-12-27.south east);
\draw[solid,lightgray](M-13-1.south west)--(M-13-27.south east);
\draw[solid,lightgray](M-1-1.north east)--(M-14-1.south east);
\draw[dotted,lightgray](M-1-2.north east)--(M-14-2.south east);
\draw[dotted,lightgray](M-1-3.north east)--(M-14-3.south east);
\draw[dotted,lightgray](M-1-4.north east)--(M-14-4.south east);
\draw[dotted,lightgray](M-1-5.north east)--(M-14-5.south east);
\draw[-latex,thick,gray](M-1-6.north east)--(M-14-6.south east);
\node (i) at (M-15-7) {$i$};
\draw[dotted,lightgray](M-1-7.north east)--(M-14-7.south east);
\draw[dotted,lightgray](M-1-8.north east)--(M-14-8.south east);
\draw[dotted,lightgray](M-1-9.north east)--(M-14-9.south east);
\draw[dotted,lightgray](M-1-10.north east)--(M-14-10.south east);
\draw[solid,lightgray](M-1-11.north east)--(M-14-11.south east);
\draw[dotted,lightgray](M-1-12.north east)--(M-14-12.south east);
\draw[dotted,lightgray](M-1-13.north east)--(M-14-13.south east);
\draw[dotted,lightgray](M-1-14.north east)--(M-14-14.south east);
\draw[dotted,lightgray](M-1-15.north east)--(M-14-15.south east);
\draw[solid,lightgray](M-1-16.north east)--(M-14-16.south east);
\draw[dotted,lightgray](M-1-17.north east)--(M-14-17.south east);
\draw[dotted,lightgray](M-1-18.north east)--(M-14-18.south east);
\draw[dotted,lightgray](M-1-19.north east)--(M-14-19.south east);
\draw[dotted,lightgray](M-1-20.north east)--(M-14-20.south east);
\draw[solid,lightgray](M-1-21.north east)--(M-14-21.south east);
\draw[dotted,lightgray](M-1-22.north east)--(M-14-22.south east);
\draw[dotted,lightgray](M-1-23.north east)--(M-14-23.south east);
\draw[dotted,lightgray](M-1-24.north east)--(M-14-24.south east);
\draw[dotted,lightgray](M-1-25.north east)--(M-14-25.south east);
\draw[solid,lightgray](M-1-26.north east)--(M-14-26.south east);
\begin{scope}[on background layer]
\draw[thin,red](M-14-5.south) -- (M-3-5.center) -| (M-2-9.center)  -- (M-2-27.east);
\draw[thin,red](M-14-9.south) -- (M-4-9.center)  -- (M-4-27.east);
\end{scope}
\node[anchor=east,scale=1.8] at (M.west) {$A^{(3)}=$};
\end{tikzpicture}

\begin{tikzpicture}[every node/.style={font=\footnotesize,scale=0.70}]
\matrix (M)[matrix of math nodes,nodes in empty cells,nodes={rectangle,minimum height=1.0em,minimum width=1.0em,inner sep=0pt,anchor=center,align=center}]
{
& &&&&& &&&&& &&&&& &&&&& &&&&& &\\
& &&&&& &&&&& &&&&& &&&&& &&&&& &\\
& &&&&& &&1&&& &&&&& &&&&& &&&&& &\\
& &&&&& &&&1&& &&&&& &&&&& &&&&& &\\
& &&&&& &&&&& &&&&& &&&&& &&&&& &\\
& &&&&& &&&&& &&&&& &&&&& &&&&& &\\
& &&&&& &&&&& &&1&&& &&&&& &&&&& &\\
& &&&&& &&&&& &&&1&& &&&&& &&&&& &\\
& &&&&& &&&&& &&&&& &&&&& &&&&& &\\
& &&&&& &&&&& &&&&& &&&&& &&&&& &\\
& &&&&& &&&&& &&&&& &&1&&& &&&&& &\\
& &&&&& &&&&& &&&&& &&&1&& &&&&& &\\
& &&&&& &&&&& &&&&& &&&&& &&&&& &\\
& &&&&& &&&&& &&&&& &&&&& &&&&& &\\
& &&&&& &&&&& &&&&& &&&&& &&&&& &\\
};
\draw[-latex,thick,gray](M-1-1.south west)--(M-1-27.south east);
\node (j) at (M-2-28) {$j$};
\draw[dotted,lightgray](M-2-1.south west)--(M-2-27.south east);
\draw[dotted,lightgray](M-3-1.south west)--(M-3-27.south east);
\draw[dotted,lightgray](M-4-1.south west)--(M-4-27.south east);
\draw[solid,lightgray](M-5-1.south west)--(M-5-27.south east);
\draw[dotted,lightgray](M-6-1.south west)--(M-6-27.south east);
\draw[dotted,lightgray](M-7-1.south west)--(M-7-27.south east);
\draw[dotted,lightgray](M-8-1.south west)--(M-8-27.south east);
\draw[solid,lightgray](M-9-1.south west)--(M-9-27.south east);
\draw[dotted,lightgray](M-10-1.south west)--(M-10-27.south east);
\draw[dotted,lightgray](M-11-1.south west)--(M-11-27.south east);
\draw[dotted,lightgray](M-12-1.south west)--(M-12-27.south east);
\draw[solid,lightgray](M-13-1.south west)--(M-13-27.south east);
\draw[solid,lightgray](M-1-1.north east)--(M-14-1.south east);
\draw[dotted,lightgray](M-1-2.north east)--(M-14-2.south east);
\draw[dotted,lightgray](M-1-3.north east)--(M-14-3.south east);
\draw[dotted,lightgray](M-1-4.north east)--(M-14-4.south east);
\draw[dotted,lightgray](M-1-5.north east)--(M-14-5.south east);
\draw[-latex,thick,gray](M-1-6.north east)--(M-14-6.south east);
\node (i) at (M-15-7) {$i$};
\draw[dotted,lightgray](M-1-7.north east)--(M-14-7.south east);
\draw[dotted,lightgray](M-1-8.north east)--(M-14-8.south east);
\draw[dotted,lightgray](M-1-9.north east)--(M-14-9.south east);
\draw[dotted,lightgray](M-1-10.north east)--(M-14-10.south east);
\draw[solid,lightgray](M-1-11.north east)--(M-14-11.south east);
\draw[dotted,lightgray](M-1-12.north east)--(M-14-12.south east);
\draw[dotted,lightgray](M-1-13.north east)--(M-14-13.south east);
\draw[dotted,lightgray](M-1-14.north east)--(M-14-14.south east);
\draw[dotted,lightgray](M-1-15.north east)--(M-14-15.south east);
\draw[solid,lightgray](M-1-16.north east)--(M-14-16.south east);
\draw[dotted,lightgray](M-1-17.north east)--(M-14-17.south east);
\draw[dotted,lightgray](M-1-18.north east)--(M-14-18.south east);
\draw[dotted,lightgray](M-1-19.north east)--(M-14-19.south east);
\draw[dotted,lightgray](M-1-20.north east)--(M-14-20.south east);
\draw[solid,lightgray](M-1-21.north east)--(M-14-21.south east);
\draw[dotted,lightgray](M-1-22.north east)--(M-14-22.south east);
\draw[dotted,lightgray](M-1-23.north east)--(M-14-23.south east);
\draw[dotted,lightgray](M-1-24.north east)--(M-14-24.south east);
\draw[dotted,lightgray](M-1-25.north east)--(M-14-25.south east);
\draw[solid,lightgray](M-1-26.north east)--(M-14-26.south east);
\begin{scope}[on background layer]
\draw[thin,red](M-14-9.south) -- (M-3-9.center)  -- (M-3-27.east);
\draw[thin,red](M-14-10.south) -- (M-4-10.center)  -- (M-4-27.east);
\end{scope}
\node[anchor=east,scale=1.8] at (M.west) {$A^{(4)}=$};
\end{tikzpicture}
\end{center}
with $A^{(5)} = {\mathbb O}$, and
\begin{equation*}
\begin{split}
{\bf s}^{(2)} & = \{ \cdots >_{\tt NW} (1, 3) >_{\tt NW} (2, 6) >_{\tt NW} (4, 7) >_{\tt NW} \cdots \}, \\
{\bf s}^{(3)} & = \{ \cdots >_{\tt NW} (2, 0) >_{\tt NW} (3, 2) >_{\tt NW} \cdots \}, \\
{\bf s}^{(4)} & = \{ \cdots >_{\tt NW} (1, -1) >_{\tt NW} (3, 3) >_{\tt NW} \cdots \}, \\
{\bf s}^{(5)} & = \{ \cdots >_{\tt NW} (2, 3) >_{\tt NW} (3, 4) >_{\tt NW} \cdots \}.
\end{split}
\end{equation*}
Hence $\kappa_0(A)=(P_0, Q_0, \rho)$, where

\begin{equation*}
P_0 = 
\begin{ytableau}
2&1&2&1&2\\
3&2&3&3&3\\
4&4
\end{ytableau}\,, \quad
Q_0 = 
\begin{ytableau}
1&1&2&3&3\\
2&2&5&4&4\\
4&3
\end{ytableau}\,, \quad
\rho = (2, 2, -1, -1, 0).
\end{equation*}
}
\end{ex}
\vskip 2mm

\begin{rem}\label{rem:fw and st commute}
{\rm
We have shown that $(d^{\tt sw}_A)^{\tt st}=d^{\tt sw}_{A^{\tt st}}$ in \eqref{eq:sw proper A and Ast}. This implies that  
 \begin{equation*}\label{eq:flat and st}
 (A^{\tt st})^\flat = (A^\flat)^{\tt st},
 \end{equation*}
{in the sense that the right-hand side is obtained by removing all zero rows and columns} on the left-hand side (cf.~Section \ref{subsec:standardization} for the convention of row and column indices of $A^{\tt st}$)
}
\end{rem}

\subsection{Semistandard tableaux of rectangular shape}\label{subsec:aff.tableaux}
Let $\la\in\cP$ be given.
For $T\in SST_\Z(\la)$, let us write $T=(T^s,\dots,T^2,T^1)$, where $T^j$ denotes the $j$-th column of $T$ from the right. For any $1\le j<k\le s$, we may regard $(T^{k},\dots,T^{j+1},T^{j})$ as a semistandard subtableau of $T$ with the columns $T^{k},\dots,T^{j+1},T^{j}$.

For $a\ge 1$ and $1\le b\le n$, let 
\begin{equation*}
 R = (a^b)=(\underbrace{a,\dots,a}_{b})\in \cP_n
\end{equation*}
be a Young diagram of rectangular shape.
Let
\begin{equation*}
\mc{B}(R)=
\left\{\,T\in SST_\Z(R)\ \Big\vert \ T^j(b)-T^j(1)<n\ (1\le j\le a)\,\right\},
\end{equation*}
where $T^j(i)$ denotes the entry of $T^j$ at the $i$-th row from the top.

For $T\in \mc{B}((1^b))$, let $\tau_n(T)=\tau(T)$ be the column semistandard tableau obtained from $T$ by replacing its entries $T(1)<T(2)<\dots<T(b)$ with 
$T(2)<T(3)<\dots <T(b)<T(1)+n.$ 
Hence $\tau$ defines a bijection on $\mc{B}((1^b))$ and $\tau^{-1}$ denotes its inverse.
In general, for $\alpha=(\alpha_1,\dots,\alpha_a)\in \Z^a$ and $(T_1,\dots,T_a)\in \mc{B}((1^b))\times \dots \times \mc{B}((1^b))$, we define
\begin{equation}\label{eq:tau}
\tau^{\alpha}(T_1,\dots,T_a) = \left(\tau^{\alpha_1}(T_1),\dots,\tau^{\alpha_a}(T_a)\right),
\end{equation}
which gives a bijection on $\mc{B}((1^b))\times \dots \times \mc{B}((1^b))$.
Regarding $\mc{B}(R)\subset \mc{B}((1^b))\times\dots\times \mc{B}((1^b))$, we may define $\tau^\alpha$ on $\mc{B}(R)$.
Let
\begin{equation*}
\begin{split}
&\mc{B}(R)_0 =
\left\{\,T \in \mc{B}(R)\ \Big\vert \ \text{$(T^{j+1},\tau^{-1}(T^{j}))$ is not semistandard for $1\le j\le a-1$}\,\right\}.
\end{split}
\end{equation*}

\begin{lem}\label{lem:S0 and S}
We have a bijection
\begin{equation*}
\xymatrixcolsep{3pc}\xymatrixrowsep{0.5pc}\xymatrix{
\ \mc{B}(R)_0\times \cP_{a-1}  \ \ar@{->}[r]  & \ \
 \mc{B}(R)  \\
\quad (T,\nu) \  \ar@{|->}[r]  &\ \ \tau^{\nu_{\tt rev}}(T)},
\end{equation*}
for $T=(T^a,\dots,T^1)$ and $\nu=(\nu_1,\dots,\nu_{a-1})\in\cP_{a-1}$,
where $\nu_{\tt rev}=(0,\nu_{a-1},\dots,\nu_2,\nu_1)$.
\end{lem}

Let us also regard $CSST_{[n]}(R)\subset \mc{B}((1^b))\times\dots\times \mc{B}((1^b))$, and define $\tau^\alpha : CSST_{[n]}(R) \longrightarrow \mc{B}((1^b))^a$ for $\alpha=(\alpha_1,\dots,\alpha_a)\in \Z^a$ as in \eqref{eq:tau}.

\begin{df}\label{df:offset constant}{\rm
Let $T \in CSST_{[n]}(R)$ be given. For $1\le j\le a-1$, let $r_j$ be the minimal non-negative integer such that 
$$(T^{j+1},\tau^{r_j}(T^{j}))$$ 
is $\Z$-semistandard,
and put $\eta_j = r_j + r_{j+1} +\dots + r_{a-1}$. We call $(r_1, \dots, r_{a-1})$ the {\em offset vector} and $\eta = (\eta_1, \dots, \eta_{a-1})$ the {\em symmetrized offset vector} of $T$.}
\end{df}

The symmetrized offset vector is the unique $\eta \in \cP_{a-1}$ such that $\tau^{\eta_{\tt rev}}(T) \in \cB(R)_0$. The following is immediate from the definition of $\eta$.

\begin{lem}\label{lem:CSST and B}
Let $T \in CSST_{[n]}(R)$ and $\alpha=(\alpha_1, \dots, \alpha_a) \in \Z^a$ be given. Then $\tau^{\alpha_{\tt rev}}(T) \in \cB(R)$ with $\alpha_{\tt rev} = (\alpha_a, \dots, \alpha_1)$ if and only if
$$ \alpha_a  \le \alpha_{a-1} - \eta_{a-1} \le \cdots \le \alpha_1 - \eta_1 .$$
\end{lem}

By Lemma \ref{lem:CSST and B}, we have a bijection
\begin{equation}\label{eq:CSST and S}
\xymatrixcolsep{3pc}\xymatrixrowsep{0.5pc}\xymatrix{
\ CSST_{[n]}(R) \times \mathcal{P}_a  \ \ar@{->}[r]  & \ \
 \mc{B}(R)  \\
\quad (T,\nu) \  \ar@{|->}[r]  &\ \ \tau^{\nu_{\tt rev} + \eta_{\tt rev}}(T)},
\end{equation}
where $\mc{P}_a = \{\, \nu=(\nu_1, \dots, \nu_a) \in \Z^a \, \vert \, \nu_1 \ge \cdots \ge \nu_a \, \}$ is the set of generalized partitions of length $a$ and $\eta$ is the symmetrized offset vector of $T$.
Let
\begin{equation}\label{eq:tau on rectangle}
\xymatrixcolsep{2.5pc}\xymatrixrowsep{0.5pc}\xymatrix{
\uptau : \mc{B}(R)  \ \ar@{->}[r]  & \ \
 \mc{B}(R)\\
 \ \ T =(T^a,\dots, T^1)\  \ar@{|->}[r]  &\ \ (\tau(T^a),\dots, \tau(T^1))}
\end{equation}
be the bijection given by applying $\tau$ to each column of the tableaux in $\mc{B}(R)$, which induces a $\Z$-action on $\mc{B}(R)$ and $\mc{B}(R)_0$. 
Let $\cB(R)_0 \big/ \Z$ denote the set of equivalence classes under this $\Z$-action. We may identify $\cB(R)_0 \big/ \Z$ with the set of $T \in \cB(R)_0$ such that the first column {has} entries in $[n]$. 
Hence, we have another bijection 
\begin{equation}\label{eq:CSST and S0}
\xymatrixcolsep{3pc}\xymatrixrowsep{0.5pc}\xymatrix{
\ CSST_{[n]}(R) \ar@{->}[r]  & \ \
 \mc{B}(R)_0\big/\Z  \\
\quad T \  \ar@{|->}[r]  &\ \ [\tau^{\eta_{\tt rev}}(T)]},
\end{equation} 
where $[T]$ denotes the equivalence class of $T$ and $\eta$ is the symmetrized offset vector of $T$.


\subsection{Rectangular decomposition}\label{subsec: rec decomp}

Let $\la\in\cP_n$ be given.
{We decompose $\la$ (its Young diagram) into diagrams of rectangular shapes $R_i$ defined by}
\begin{equation}\label{eq:rec decomp for lambda}
R_i=(\underbrace{m_i,\dots,m_i}_{i})\quad (1\le i\le l),
\end{equation}
where $m_i$ is the number of occurrences of $i$ in $\la'$ and $l=\ell(\la)$.
Here we assume that $R_i$ is empty when $m_i=0$. {For example, if $\la = (6,4,1,1)$, then see that $R_1 = (1)$, $R_2 = (3^2)$, $R_3 = \emptyset$ and $R_4 = (2^4)$ as illustrated in the following figure.}\vskip 2mm

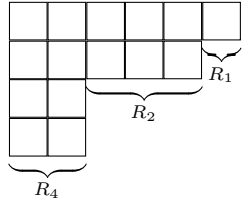
\begin{center}
\begin{tikzpicture}[
BC/.style = {decorate, decoration={calligraphic brace, amplitude=5pt, raise=1mm},
        thick, pen colour={black}
            },
                    ]
\matrix (m) [matrix of math nodes,
             nodes={draw, minimum size=5mm, anchor=center},
             ]
{
~&~&~&~&~&~\\
~&~&~&~&~& \\
~&~& & & & \\
~&~& & & & \\
};
\draw[BC] (m-4-2.south east) -- node[below=2.2mm] {${}_{R_4}$} (m-4-1.south west);
\draw[BC] (m-2-5.south east) -- node[below=2.2mm] {${}_{R_2}$} (m-2-3.south west);
\draw[BC] (m-1-6.south east) -- node[below=2.2mm] {${}_{R_1}$} (m-1-6.south west);
\end{tikzpicture}
\end{center}

Let
\begin{equation}\label{eq:B(la)}
\begin{split}
 & \mc{B}_n(\la)_0=\mc{B}(\la)_0= \mc{B}(R_l)_0 \times \dots \times \mc{B}(R_1)_0,\\
 & \mc{B}_n(\la)=\mc{B}(\la)= \mc{B}(R_l) \times \dots \times \mc{B}(R_1).
\end{split} 
\end{equation}
If we put $\cP(\la)=\cP_{m_l-1}\times \dots \times \cP_{m_2-1}\times \cP_{m_1-1}$, where we take the product over $m_i \ge 1$, 
then we have a bijection
\begin{equation*}
\xymatrixcolsep{3pc}\xymatrixrowsep{0.5pc}\xymatrix{
\ \mc{B}(\la)_0\times \cP(\la)  \ \ar@{->}[r]  & \ \
 \mc{B}(\la) \\
\left(\left(T^{(i)}\right)_{1\le i\le l},(\nu^{(i)})_{1\le i\le l}\right) \  \ar@{|->}[r]  &\ \ \left(\tau^{\nu^{(i)}_{\tt rev}}(T^{(i)})\right)_{1\le i\le l}},
\end{equation*}
by applying Lemma \ref{lem:S0 and S} to each component, where $T^{(i)}\in \mc{B}(R_i)_0$ and $\nu^{(i)}\in \cP_{m_i-1}$. 

Similarly, if we let $\mc{P}(\la) = \mc{P}_{m_l}\times \dots \times \mc{P}_{m_1}$, where we take the product over $m_i \ge 1$, 
and regard
\begin{equation*}
CSST_{[n]}(\la) = CSST_{[n]}(R_l) \times \dots \times CSST_{[n]}(R_1),
\end{equation*} 
then by \eqref{eq:tau on rectangle} we have a bijection
\begin{equation}\label{eq:CSST and B general}
\xymatrixcolsep{3pc}\xymatrixrowsep{0.5pc}\xymatrix{
\ CSST_{[n]}(\la) \times \mc{P}(\la)  \ \ar@{->}[r]  & \ \
 \mc{B}(\la) \\
 \left(\left(T^{(i)}\right)_{1\le i\le l},(\nu^{(i)})_{1\le i\le l}\right) \  \ar@{|->}[r]  &\ \ \left(\tau^{\nu^{(i)}_{\tt rev} + \eta^{(i)}_{\tt rev}}(T^{(i)})\right)_{1\le i\le l}},
\end{equation} 
where $\eta^{(i)} \in \cP_{m_i-1}$ is the symmetrized offset vector of $T^{(i)}$

\begin{ex}{\rm
Let $P_0 \in CSST_{[4]}(\lambda)$ be the tableau in Example \ref{ex:kappa_0 of A}, i.e.,
\begin{equation*}
P_0 = 
\begin{ytableau}
2&1&2&1&2\\
3&2&3&3&3\\
4&4
\end{ytableau}\    
\end{equation*}
where $\lambda=(5,5,2)$. Then $\lambda$ is decomposed into $R_2 = (3,3)$ and $R_3 = (2,2,2)$, and the corresponding decompositions of $P_0$ are 

\begin{equation*}
P^{(3)}_0 = \begin{ytableau}
2&1\\
3&2\\
4&4
\end{ytableau}\ , \quad\quad
P^{(2)}_0 = \begin{ytableau}
2&1&2\\
3&3&3
\end{ytableau}.
\end{equation*}
Since

\begin{equation*}
\tau^{(0,1)}(P^{(3)}_0) = 
\begin{ytableau}
2&2\\
3&4\\
4&5\end{ytableau} \ \in \mc{B}(R_3)_0,
\end{equation*} 

\noindent the symmetrized offset vector of $P^{(3)}_0$ is $\eta^{(3)}=(1) \in \cP_{1}$. Similarly, we have $\eta^{(2)}=(1,1) \in \cP_{2}$. For $\nu = ( (1, -1), ((2,1,0)) \in \mathcal{P}(\lambda)$, the image of $(P_0, \nu)$ under the bijection \eqref{eq:CSST and B general} is

\begin{equation*}
    \begin{ytableau}
0&4&2&5&7\\
2&5&3&7&10\\
3&6
\end{ytableau} \ \in \mc{B}(\lambda).
\end{equation*}

} 
\end{ex}\vskip 2mm
 
Let 
\begin{equation*}
 \mc{B}(\la)_0 \big/{\Z^l} = \mc{B}(R_l)_0 \big/ {\Z} \times \dots \times \mc{B}(R_1)_0 \big/{\Z},
\end{equation*}
where each $\mc{B}(R_i)_0 \big/{\Z}$ is the set of equivalence classes under the $\Z$-action \eqref{eq:tau on rectangle}. 
Then we also have a bijection
\begin{equation*}
\xymatrixcolsep{3pc}\xymatrixrowsep{0.5pc}\xymatrix{
\ CSST_{[n]}(\la)  \ \ar@{->}[r]  & \ \
 \mc{B}(\la)_0 \big/ {\Z^l} \\
 \left(T^{(i)}\right)_{1\le i\le l} \  \ar@{|->}[r]  &\ \ \left([\tau^{\eta^{(i)}_{\tt rev}}(T^{(i)})]\right)_{1\le i\le l} \\
},
\end{equation*}
where $\eta^{(i)} \in \cP_{m_i-1}$ is the symmetrized offset vector of $T^{(i)}$.

\subsection{Affine RSK correspondence}\label{subsec:affine RSK}

Suppose that $A\in \wh{\M}_{m\times n}$ is given. We keep the notations in Sections \ref{subsec:mbc} and \ref{subsec: rec decomp}.

Let $(P_0,Q_0,\rho)$ be given as in \eqref{eq:Aff RSK 0}. For $1\le i\le l$ with $m_i \ge 1$, let
\begin{itemize}
\item[$\bullet$] $P^{(i)}_0$, $Q^{(i)}_0$ : the subtableaux of $P_0$ and $Q_0$ corresponding to $R_i$, respectively,

\item[$\bullet$] $\rho^{(i)}\in \Z^{m_i}$ : the subsequence of $\rho$ corresponding to the columns of $R_i$, 

\item[$\bullet$] $\eta^{(i)}\in \cP_{m_i-1}$ : the symmetrized offset vector of $P^{(i)}_0$.

\end{itemize}
Then we define 
\begin{equation}\label{eq:(P,Q)}
\begin{split}
Q&=\left(\tau^{\rho^{(l)}+\eta^{(l)}_{\tt rev}}(Q^{(l)}_0),\dots,\tau^{\rho^{(1)}+\eta^{(1)}_{\tt rev}}(Q^{(1)}_0)\right)
=\left(Q^{(l)},\dots,Q^{(1)} \right).
\end{split}
\end{equation}
Note that the action of $\tau$ on $Q_0$ should be understood as $\tau_n$.
The following is one of the main results in this paper. The proof will be given in Section \ref{sec:Proof of affine RSK}.
\begin{thm}\label{thm:main-1}
We have a bijection
\begin{equation*}
\xymatrixcolsep{3pc}\xymatrixrowsep{0.5pc}\xymatrix{
\kappa:\ \ \wh{\M}_{m\times n} \ \ar@{->}[r]  & \ \
\displaystyle{\bigsqcup_{\la\in \cP_m\cap\cP_n} CSST_{[m]}(\la) \times \mc{B}_n(\la)} \\
\ \ A \  \ar@{|->}[r]  &\ \ (P_0,Q)}.
\end{equation*}
\end{thm}

\begin{rem}{\rm
Applying the bijection \eqref{eq:CSST and B general} to $\mc{B}_n(\la)$, we have a bijection
\begin{equation*}
\xymatrixcolsep{3pc}\xymatrixrowsep{0.5pc}\xymatrix{
\kappa^\diamond:\ \ \wh{\M}_{m\times n} \ \ar@{->}[r]  & \ \
\displaystyle{\bigsqcup_{\la\in \cP_m\cap\cP_n} CSST_{[m]}(\la) \times CSST_{[n]}(\la) \times \mc{P}(\la)} }.
\end{equation*}
Again applying the inverse of \eqref{eq:CSST and B general} to $CSST_{[m]}(\la) \times \mc{P}(\la)$, we have a bijection
\begin{equation*}
\xymatrixcolsep{3pc}\xymatrixrowsep{0.5pc}\xymatrix{
\kappa' :\ \ \wh{\M}_{m\times n} \ \ar@{->}[r]  & \ \
\displaystyle{\bigsqcup_{\la\in \cP_m\cap\cP_n} \mc{B}_m(\la) \times  CSST_{[n]}(\la)}}.
\end{equation*}
 
} 
\end{rem}

\begin{ex}\label{ex:kappa of A}{\rm
Let $\kappa_0(A) = (P_0, Q_0, \rho)$ be as in Example \ref{ex:kappa_0 of A} {($m=4$, $n=5$)}.
The rectangular decompositions of $P_0$, $Q_0$ and $\rho$ are

\begin{align*}
P^{(2)}_0 & = \begin{ytableau}
2&1&2\\
3&3&3
\end{ytableau}\,, & Q^{(2)}_0 & = \begin{ytableau}
2&3&3\\
5&4&4
\end{ytableau}\,, & \rho^{(2)} & = (-1, -1, 0),\\ & & & \\
P^{(3)}_0 & = \begin{ytableau}
2&1\\
3&2\\
4&4
\end{ytableau}\,, & Q^{(3)}_0 & = \begin{ytableau}
1&1\\
2&2\\
4&3
\end{ytableau}\,, & \rho^{(3)} &= (2,2),
\end{align*}
\vskip 2mm

\noindent where $R_2$ and $R_3$ are the only non-trivial rectangles in this decomposition.
The symmetrized offset vectors of $P^{(2)}_0$ and $P^{(3)}_0$ are
$$\eta^{(2)}_{\tt rev} = (0, 1, 1),\quad \eta^{(3)}_{\tt rev} = (0, 1),$$ 
and hence
$$Q = 
\begin{ytableau}
4&6&0&3&4\\
6&7&2&4&8\\
7&8
\end{ytableau}\,.
$$
}

\end{ex}

\begin{rem}{\rm
Let us give some comments on $\kappa_0$ and the bijection in \cite{IMS}.
Let $\ov{\mathbb{M}}_{m\times n}$ be the set of $\ov{M} = \big(\ov{M}_{j,i}(k) \big)$ ($i\in [m]$, $j\in [n]$, $k \in \Z$) with $\ov{M}_{j,i}(k)\in\Z_{\ge 0}$ and $\ov{M}_{j,i}(k)=0$ for $|k|\gg 0$ \cite[(2.5)]{IMS} and let $\ov{\mathbb{M}}^+_{m\times n}$ be the subset of $\ov{\mathbb{M}}_{m\times n}$ consisting of $\ov{M}$ such that $\ov{M}_{j,i}(k)=0$ for $k<0$.
For $A = (a_{ij})_{i,j \in \Z} \in \wh{\M}_{m\times n}$, we define
$\ov{M}_A = \big(\ov{M}_{j,i}(k) \big)$ by 
\begin{equation*}
    \ov{M}_{j,i}(k) = a_{i-km,\, n + 1 - j}.
\end{equation*}
Then the map sending $A$ to $\ov{M}_A$ gives a bijection from $\wh{\M}_{m\times n}$ to $\ov{\mathbb{M}}_{m\times n}$. 

Let $\td{\Upsilon}$ denote the bijection 
\begin{equation*}
\xymatrixcolsep{3pc}\xymatrixrowsep{0.5pc}\xymatrix{
\ov{\mathbb{M}}^+_{m\times n} \ \ar@{->}[r]  & \ \
\displaystyle{\bigsqcup_{\la\in \cP_m\cap\cP_n}  CSST_{[m]}(\la)}\times CSST_{[n]}(\la)\times \mc{K}(\la)},
\end{equation*}
given in \cite[Corollary 8.2]{IMS}.
 
Let $A = (a_{ij})_{i,j \in \Z} \in \wh{\M}_{m\times n}$ be given such that $\ov{M}_A\in \ov{\mathbb{M}}^+_{m\times n}$. Applying $\kappa_0$ and $\td{\Upsilon}$  directly to $A$ and $\ov{M}_A$, respectively, do not seem to give the same result in general. This may happen due to conventions for generalized affine permutations.  
For example, let $A\in \wh{\M}_{5\times 6}$ be as follows:
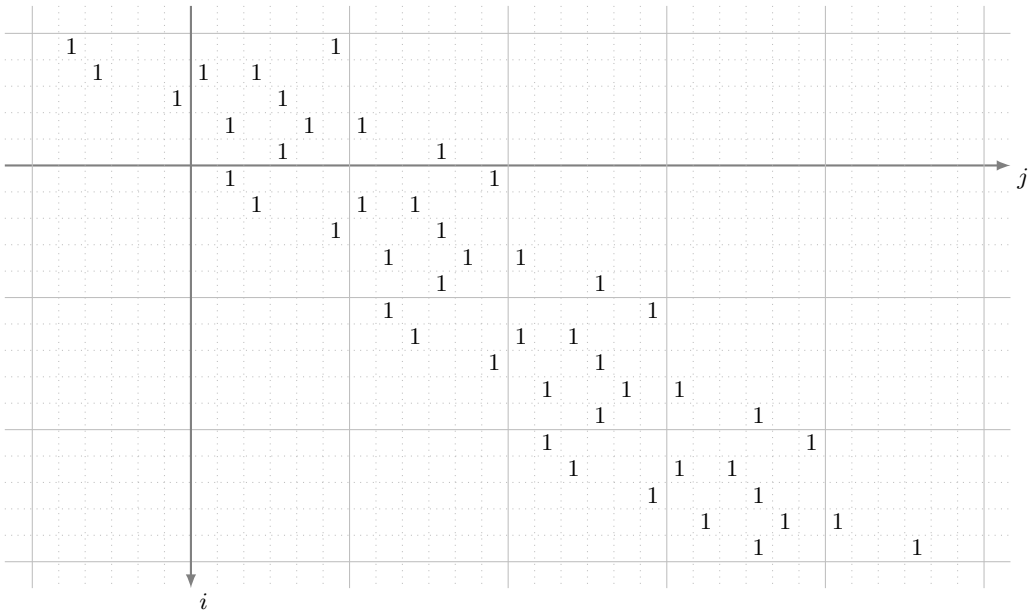
\begin{center}
\begin{tikzpicture}[every node/.style={font=\small}]
\matrix (M)[matrix of math nodes,nodes in empty cells,nodes={rectangle,minimum height=1.0em,minimum width=1.0em,inner sep=0pt,anchor=center,align=center}]
{
& &&&&&& &&&&&& &&&&&& &&&&&& &&&&&& &&&&&& &\\
& &1&&&&& &&&&&1& &&&&&& &&&&&& &&&&&& &&&&&& &\\
& &&1&&&& 1&&1&&&& &&&&&& &&&&&& &&&&&& &&&&&& &\\
& &&&&&1& &&&1&&& &&&&&& &&&&&& &&&&&& &&&&&& &\\
& &&&&&& &1&&&1&& 1&&&&&& &&&&&& &&&&&& &&&&&& &\\
& &&&&&& &&&1&&& &&&1&&& &&&&&& &&&&&& &&&&&& &\\
& &&&&&& &1&&&&& &&&&&1& &&&&&& &&&&&& &&&&&& &\\
& &&&&&& &&1&&&& 1&&1&&&& &&&&&& &&&&&& &&&&&& &\\
& &&&&&& &&&&&1& &&&1&&& &&&&&& &&&&&& &&&&&& &\\
& &&&&&& &&&&&& &1&&&1&& 1&&&&&& &&&&&& &&&&&& &\\
& &&&&&& &&&&&& &&&1&&& &&&1&&& &&&&&& &&&&&& &\\
& &&&&&& &&&&&& &1&&&&& &&&&&1& &&&&&& &&&&&& &\\
& &&&&&& &&&&&& &&1&&&& 1&&1&&&& &&&&&& &&&&&& &\\
& &&&&&& &&&&&& &&&&&1& &&&1&&& &&&&&& &&&&&& &\\
& &&&&&& &&&&&& &&&&&& &1&&&1&& 1&&&&&& &&&&&& &\\
& &&&&&& &&&&&& &&&&&& &&&1&&& &&&1&&& &&&&&& &\\
& &&&&&& &&&&&& &&&&&& &1&&&&& &&&&&1& &&&&&& &\\
& &&&&&& &&&&&& &&&&&& &&1&&&& 1&&1&&&& &&&&&& &\\
& &&&&&& &&&&&& &&&&&& &&&&&1& &&&1&&& &&&&&& &\\
& &&&&&& &&&&&& &&&&&& &&&&&& &1&&&1&& 1&&&&&& &\\
& &&&&&& &&&&&& &&&&&& &&&&&& &&&1&&& &&&1&&& &\\
& &&&&&& &&&&&& &&&&&& &&&&&& &&&&&& &&&&&& &\\
& &&&&&& &&&&&& &&&&&& &&&&&& &&&&&& &&&&&& &\\
};
\draw[solid,lightgray](M-1-1.south west)--(M-1-38.south east);
\draw[dotted,lightgray](M-2-1.south west)--(M-2-38.south east);
\draw[dotted,lightgray](M-3-1.south west)--(M-3-38.south east);
\draw[dotted,lightgray](M-4-1.south west)--(M-4-38.south east);
\draw[dotted,lightgray](M-5-1.south west)--(M-5-38.south east);
\draw[-latex,thick,gray](M-6-1.south west)--(M-6-38.south east);
\node (j) at (M-7-39) {$j$};
\draw[dotted,lightgray](M-7-1.south west)--(M-7-38.south east);
\draw[dotted,lightgray](M-8-1.south west)--(M-8-38.south east);
\draw[dotted,lightgray](M-9-1.south west)--(M-9-38.south east);
\draw[dotted,lightgray](M-10-1.south west)--(M-10-38.south east);
\draw[solid,lightgray](M-11-1.south west)--(M-11-38.south east);
\draw[dotted,lightgray](M-12-1.south west)--(M-12-38.south east);
\draw[dotted,lightgray](M-13-1.south west)--(M-13-38.south east);
\draw[dotted,lightgray](M-14-1.south west)--(M-14-38.south east);
\draw[dotted,lightgray](M-15-1.south west)--(M-15-38.south east);
\draw[solid,lightgray](M-16-1.south west)--(M-16-38.south east);
\draw[dotted,lightgray](M-17-1.south west)--(M-17-38.south east);
\draw[dotted,lightgray](M-18-1.south west)--(M-18-38.south east);
\draw[dotted,lightgray](M-19-1.south west)--(M-19-38.south east);
\draw[dotted,lightgray](M-20-1.south west)--(M-20-38.south east);
\draw[solid,lightgray](M-21-1.south west)--(M-21-38.south east);
\draw[solid,lightgray](M-1-1.north east)--(M-22-1.south east);
\draw[dotted,lightgray](M-1-2.north east)--(M-22-2.south east);
\draw[dotted,lightgray](M-1-3.north east)--(M-22-3.south east);
\draw[dotted,lightgray](M-1-4.north east)--(M-22-4.south east);
\draw[dotted,lightgray](M-1-5.north east)--(M-22-5.south east);
\draw[dotted,lightgray](M-1-6.north east)--(M-22-6.south east);
\draw[-latex,thick,gray](M-1-7.north east)--(M-22-7.south east);
\node (i) at (M-23-8) {$i$};
\draw[dotted,lightgray](M-1-8.north east)--(M-22-8.south east);
\draw[dotted,lightgray](M-1-9.north east)--(M-22-9.south east);
\draw[dotted,lightgray](M-1-10.north east)--(M-22-10.south east);
\draw[dotted,lightgray](M-1-11.north east)--(M-22-11.south east);
\draw[dotted,lightgray](M-1-12.north east)--(M-22-12.south east);
\draw[solid,lightgray](M-1-13.north east)--(M-22-13.south east);
\draw[dotted,lightgray](M-1-14.north east)--(M-22-14.south east);
\draw[dotted,lightgray](M-1-15.north east)--(M-22-15.south east);
\draw[dotted,lightgray](M-1-16.north east)--(M-22-16.south east);
\draw[dotted,lightgray](M-1-17.north east)--(M-22-17.south east);
\draw[dotted,lightgray](M-1-18.north east)--(M-22-18.south east);
\draw[solid,lightgray](M-1-19.north east)--(M-22-19.south east);
\draw[dotted,lightgray](M-1-20.north east)--(M-22-20.south east);
\draw[dotted,lightgray](M-1-21.north east)--(M-22-21.south east);
\draw[dotted,lightgray](M-1-22.north east)--(M-22-22.south east);
\draw[dotted,lightgray](M-1-23.north east)--(M-22-23.south east);
\draw[dotted,lightgray](M-1-24.north east)--(M-22-24.south east);
\draw[solid,lightgray](M-1-25.north east)--(M-22-25.south east);
\draw[dotted,lightgray](M-1-26.north east)--(M-22-26.south east);
\draw[dotted,lightgray](M-1-27.north east)--(M-22-27.south east);
\draw[dotted,lightgray](M-1-28.north east)--(M-22-28.south east);
\draw[dotted,lightgray](M-1-29.north east)--(M-22-29.south east);
\draw[dotted,lightgray](M-1-30.north east)--(M-22-30.south east);
\draw[solid,lightgray](M-1-31.north east)--(M-22-31.south east);
\draw[dotted,lightgray](M-1-32.north east)--(M-22-32.south east);
\draw[dotted,lightgray](M-1-33.north east)--(M-22-33.south east);
\draw[dotted,lightgray](M-1-34.north east)--(M-22-34.south east);
\draw[dotted,lightgray](M-1-35.north east)--(M-22-35.south east);
\draw[dotted,lightgray](M-1-36.north east)--(M-22-36.south east);
\draw[solid,lightgray](M-1-37.north east)--(M-22-37.south east);
\end{tikzpicture}
\end{center}
Then we have $\kappa_0(A)=(P_0,Q_0,\rho)$, where

\begin{equation*}
    P_0 = \ytableaushort{1221,4342,5453}\ , \quad Q_0 = \ytableaushort{2131,3244,6456}\ .
\end{equation*}

\noindent On the other hand, we have $\td{\Upsilon}(\ov{M}_A)=(V,W,\xi)$, where 

\begin{equation*}
    V = \ytableaushort{1221,4342,5453}\ , \quad W = \ytableaushort{1231,3354,6465}\ .
\end{equation*} 

\noindent Hence $P_0=V$ but $Q_0\neq W$, while we observe that $W^j(i) = 7 - Q_0^{5-j}(4-i)$ for $1\le i \le 3$ and $1\le j\le 4$. Recall that $W^j$ and $Q_0^j$ denote the $j$-th columns from the right. 

In general, one may expect that $P_0=V$ and $Q_0 = e_nW$, where $e_n$ is an operator on $CSST_{[n]}(\la)$ given by
\begin{equation*}
(e_nW)^j(i) = n + 1 - W^{a + 1 - j}(b + 1 - i) \quad (1\le i\le b, 1\le j\le a),
\end{equation*}
in case of a rectangular shape $\la=(a^b)$. The operator $e_n$ can be viewed as a generalization of the affine evacuation in \cite{CFKLY}. We do not know yet a precise relation between $\rho$ and $\xi$.
} 
\end{rem}


\section{Affine crystals}\label{sec:aff bicrystal} 

\subsection{Crystals} 
Let us give a brief review on crystals (see \cite{HK} for more details). 
Let $\g$ be the Kac-Moody algebra associated to a symmetrizable generalized Cartan matrix $A =(a_{ij})_{i,j\in I}$ indexed by a set $I$. Let $P^\vee$ be the dual weight lattice, $P = {\rm Hom}_\Z( P^\vee,\Z)$ the weight lattice, $\Pi^\vee=\{\,h_i\,|\,i\in I\,\}\subset P^\vee$ the set of simple coroots, and $\Pi=\{\,\alpha_i\,|\,i\in I\,\}\subset P$ the set of simple roots of $\g$ such that $\langle \alpha_j,h_i\rangle=a_{ij}$ for $i,j\in I$.

A {\it $\g$-crystal} (or crystal if there is no confusion on $\mf g$) is a set $B$ together with the maps ${\rm wt} : B \longrightarrow P$, $\varepsilon_i, \varphi_i: B \longrightarrow \mathbb{Z}\cup\{-\infty\}$ and $\te_i, \tf_i: B \longrightarrow B\cup\{{\bf 0}\}$ for $i\in I$ such that
for $b\in B$ and $i\in I$
\begin{itemize}
\item[(1)] $\varphi_i(b) =\langle {\rm wt}(b),h_i \rangle +
\varepsilon_i(b),$

\item[(2)]  $\varepsilon_i(\te_i b) = \varepsilon_i(b) - 1$, $\varphi_i(\te_i b) =
\varphi_i(b) + 1$, ${\rm wt}(\te_ib)={\rm wt}(b)+\alpha_i$ if $\te_i b \neq {\bf 0}$,

\item[(3)] $\varepsilon_i(\tf_i b) = \varepsilon_i(b) + 1$, $\varphi_i(\tf_i b) =
\varphi_i(b) - 1$, ${\rm wt}({\tf_i}b)={\rm wt}(b)-\alpha_i$ if $\tf_i b \neq {\bf 0}$,

\item[(4)] $\tf_i b = b'$ if and only if $b = \te_i b'$ for $b, b' \in B$, 

\item[(5)] $\te_ib=\tf_ib={\bf 0}$ if $\varphi_i(b)=-\infty$,
\end{itemize}
where ${\bf 0}$ is a formal symbol. Here we assume that
$-\infty+n=-\infty$ for all $n\in\Z$. A crystal $B$ can be viewed as an $I$-colored graph where $b\stackrel{i}{\rightarrow}b'$ if and only if $\tf_ib=b'$ for  $b, b'\in B$ and $i\in I$. A crystal $B$ is called {\em normal} if $\varepsilon_i(b)=\max\{\,k\,|\,\te_i^k b \neq {\bf 0}\,\}$ and  $\varphi_i(b)=\max\{\,k\,|\,\tf_i^k b \neq {\bf 0}\,\}$ for $b\in B$ and $i\in I$.

Let $B_1$ and $B_2$ be crystals. A {morphism}
$\psi : B_1 \longrightarrow B_2$ is a map from $B_1\cup\{{\bf 0}\}$ to
$B_2\cup\{{\bf 0}\}$ such that
\begin{itemize}
\item[(1)] $\psi(\bf{0})=\bf{0}$,

\item[(2)] ${\rm wt}(\psi(b))={\rm wt}(b)$,
$\varepsilon_i(\psi(b))=\varepsilon_i(b)$, and
$\varphi_i(\psi(b))=\varphi_i(b)$ if $\psi(b)\neq \bf{0}$,

\item[(3)] $\psi(\te_i b)=\te_i\psi(b)$ if $\psi(b)\neq \bf{0}$ and
$\psi(\te_i b)\neq \bf{0}$,

\item[(4)] $\psi(\tf_i
b)=\tf_i\psi(b)$ if $\psi(b)\neq
\bf{0}$ and $\psi(\tf_i b)\neq \bf{0}$,
\end{itemize}
for $b\in B_1$ and $i\in I$.
We call $\psi$ an {\em embedding} when $\psi$ is injective, and call $\psi$ {\em strict} if
$\psi : B_1\cup\{{\bf 0}\} \longrightarrow B_2\cup\{{\bf 0}\}$ commutes
with $\te_i$ and $\tf_i$ for all $i\in I$, where we assume that $\te_i{\bf
0}=\tf_i{\bf 0}={\bf 0}$. When $\psi$ is a bijection, it is called an isomorphism.

Given $b_1 \in B_1$ and $b_2 \in B_2$,
we say that $b_1$ is {\em $\mf g$-crystal equivalent to $b_2$} if there is an isomorphism of crystals $C(b_1) \longrightarrow C(b_2)$ mapping $b_1$ to $b_2$, where $C(b_i)$ denotes the connected component of $b_i$ in $B_i$ for $i=1,2$.

The {\em tensor product $B_1 \ot B_2$} is defined to be $B_1 \times B_2$ as a set with elements denoted by $b_1 \ot b_2$, where
{\allowdisplaybreaks
\begin{equation} \label{eq:tensor_product_rule}
\begin{split}
	& {\rm wt}(b_1 \ot b_2) = {\rm wt}(b_1) + {\rm wt}(b_2), \\
	& \varepsilon_{i}(b_{1} \ot b_{2}) = \max\{ \varepsilon_{i}(b_{1}), \varepsilon_{i}(b_{2})-\langle {\rm wt}(b_{1}), h_i \rangle \}, \\
	& \varphi_{i}(b_{1} \ot b_{2}) = \max\{ \varphi_{i}(b_{1})+\langle {\rm wt}(b_2), h_i \rangle, \varphi_{i}(b_{2}) \}, \\
	& \tilde{e}_{i}(b_{1} \ot b_{2}) = \left\{ \begin{array}{cc} \tilde{e}_{i}b_{1} \ot b_{2} & \textrm{if} \ \varphi_{i}(b_{1}) \ge \varepsilon_{i}(b_{2}), \\ b_{1} \ot \tilde{e}_{i}b_{2} & \textrm{if} \ \varphi_{i}(b_{1}) < \varepsilon_{i}(b_{2}), \end{array} \right. \\
	& \tilde{f}_{i}(b_{1} \ot b_{2}) = \left\{ \begin{array}{cc} \tilde{f}_{i}b_{1} \ot b_{2} & \textrm{if} \ \varphi_{i}(b_{1}) > \varepsilon_{i}(b_{2}), \\ b_{1} \ot \tilde{f}_{i}b_{2} & \textrm{if} \ \varphi_{i}(b_{1}) \le \varepsilon_{i}(b_{2}),\end{array} \right.
\end{split}	
\end{equation}} 
for $i \in I$. Here, we assume that $\textbf{0} \ot b_{2} = b_{1} \ot \textbf{0} = \textbf{0}$. Then $B_1 \ot B_2$ is a crystal. 

{Let us recall a well-known combinatorial rule, which is often used to describe a crystal structure on the tensor product of crystals, and which will be used later in this paper.} Suppose that $\sigma=(\ldots,\sigma_{-2},\sigma_{-1},\sigma_0,\sigma_{1},\sigma_2,\ldots)$ is a sequence (not necessarily finite) where $\sigma_{k}\in \{\,+\,,\,-\, , \ \cdot\ \}$ and $\sigma_k=\,\cdot\,$ except for finitely many $k$.  
In the sequence $\sigma$, we replace a pair $(\sigma_{s},\sigma_{s'})=(+,-)$, where $s<s'$ and $\sigma_t=\,\cdot\,$ for $s<t<s'$, with $(\,\cdot\,,\,\cdot\,)$, and repeat this process as far as possible until we get a sequence with no $-$ placed to the right of $+$. 
We denote the resulting reduced sequence by $\td{\sigma}$, where no $+$ is placed to the left of any $-$ in the sequence. 

\subsection{$A_{n-1}^{(1)}$-crystals}
Suppose that $\mf g$ is an affine Kac-Moody algebra of type $A_{n-1}^{(1)}$ for $n\ge 2$. 
We assume that 
$I=\{\,0,1,\dots,n-1\,\}$, 
$P^\vee=\bigoplus_{i\in I}\Z h_i\oplus\Z d$, and 
$P=\bigoplus_{i\in I} \Z\La_i\oplus\Z\delta$, 
where $\delta=\sum_{i\in I}\alpha_i$ is the null root such that 
$\langle \de,h_j \rangle =0$ and $\langle \de,d \rangle =1$,
and $\La_i$ is the $i$-th fundamental weight such that $\langle \La_i,h_j \rangle =\delta_{ij}$ and $\langle \La_i,d \rangle =0$. 
Let us call a $\mf{g}$-crystal $B$ an $A_{n-1}^{(1)}$-crystal for simplicity.
For $\La\in P$, let $B(\La)$ denote the crystal of the {\em extremal weight module} generated by an extremal vector of weight $\La$ \cite{Kas94'}.

We let
\begin{equation*}
\begin{split}
 P^0&=\Z\varpi_1 \oplus \dots \oplus\Z\varpi_{n-1}\oplus\Z\de,\ \
 P_{\rm cl}=P/\Z\delta,\ \
 P^0_{\rm cl}=P^0/\Z\delta, 
\end{split}
\end{equation*}
where $\varpi_i=\La_i-\La_0$ is the $i$-th level zero fundamental weight for $i\in \{\,1,\dots,n-1\,\}$. 
Let ${\rm cl} : P \longrightarrow P_{\rm cl}$ be the canonical projection.
Define  $\epsilon_1=\varpi_1$ and $\epsilon_{i}=\varpi_i-\varpi_{i-1}$ for $i\in \{\,2,\dots,n-1\,\}$. Put $\epsilon_n=-\epsilon_1-\dots-\epsilon_{n-1}$. We have $\alpha_i=\epsilon_i-\epsilon_{i+1}$ for $i\in \{\,1,\dots,n-1\,\}$ and $\alpha_0=\epsilon_n-\epsilon_1+\de$.

In this paper, we also need to consider an $A_{n-1}^{(1)}$-crystal $B$ with ${\rm wt}$ defined on $P^0$ and $P_{\rm cl}^0$ instead of $P$. In this case, we call $B$  
an {\em $A_{n-1}^{(1)}$-crystal with $P^0$-weights and $P_{\rm cl}^0$-weights}, respectively.



\begin{ex}\label{ex:crystal of level 1 fundamental weight}
{\rm

Let $\mc{F}$ be the set of sequences ${\bf v}=(v_k)_{k\in \Z}$ such that $v_k\in \{\,0,1\,\}$, $v_k=1$ for $k\ll 0$ and $v_k=0$ for $k\gg 0$. 
For $i\in I$, 
consider a sequence
\begin{equation*}
 \sigma_{{\bf v},i} = (\dots, 
 \sigma_{i-n},\sigma_{i+1-n},
 \sigma_i,\sigma_{i+1},
 \sigma_{i+n},\sigma_{i+1+n},
 \dots ),
\end{equation*}
where we assign for $k\equiv i \pmod{n}$
\begin{equation*}
 (\sigma_{k},\sigma_{k+1})=
\begin{cases}
 (\,+\, ,\ \cdot\ ) & \text{if $(v_{k}, v_{k+1})=(1,0)$},\\
 (\ \cdot\ ,\,-\,) & \text{if $(v_{k}, v_{k+1})=(0,1)$},\\
 (\ \cdot\ ,\ \cdot\ ) & \text{if $(v_{k}, v_{k+1})=(0,0)$ or $(1,1)$}.
\end{cases}
\end{equation*}

We define $\te_i{\bf v}$ to be the sequence obtained from ${\bf v}$ by replacing $(v_k,v_{k+1})=(0,1)$ with $(1,0)$, where $k+1\equiv i+1 \pmod{n}$ is the index corresponding to the rightmost $-$ in the reduced sequence $\td{\sigma_{{\bf v},i}}$. We define $\te_i{\bf v}={\bf 0}$ if there is no such $k$.
Similarly, we define $\tf_i{\bf v}$ to be the sequence obtained from ${\bf v}$ by replacing $(v_k,v_{k+1})=(1,0)$ with $(0,1)$, where $k\equiv i \pmod{n}$ is the index corresponding to the leftmost $+$ in $\td{\sigma_{{\bf v},i}}$. We define $\tf_i{\bf v}={\bf 0}$ if there is no such $k$.

For $j\in I$, let ${\bf v}_{\La_j}=(a_k)_{k\in \Z}$ be such that $v_k=1$ if and only if $k\le j$. It is well-known that the connected component of ${\bf v}_{\La_j}$ in $\mc F$ with respect to $\te_i$ and $\tf_i$ for $i\in I$ is isomorphic to $B(\La_j)$, where we put ${\rm wt}({\bf v}_{\La_j})=\La_j$ \cite{MM}. Moreover, {if we put $\Lambda_j = \Lambda_{j+n}$ for $j \in \Z$,} we have $\mc{F} \cong \bigoplus_{j \in \Z} B(\Lambda_j)$.
  
Let $\mc{F}^\vee$ be the set of sequences ${\bf v}=(v_k)_{k\in \Z}$ such that $v_k\in \{\,0,1\,\}$, $v_k=0$ for $k\ll 0$ and $v_k=1$ for $k\gg 0$ (which is the dual of $\mc{F}$). 
Let ${\bf v}_{-\La_j}=(v_k)_{k\in \Z}$ be such that $a_k=0$ if and only if $k\le j$. Then the connected component of ${\bf v}_{-\La_j}$ in $\mc{F}^\vee$ with respect to $\te_i$ and $\tf_i$ for $i \in I$ is isomorphic to $B(-\La_j)$, where ${\rm wt}({\bf v}_{-\La_j})=-\La_j$.
 
}
\end{ex}

\subsection{Bicrystal structure on $\wh{\M}_{m \times n}$}\label{subsec:bicrystal}

In this subsection, we define two crystal structures on $\wh{\M}_{m \times n}$ of type $A_{m-1}^{(1)}$ and $A_{n-1}^{(1)}$.

For $i,j \in\Z$, let $E_{ij}$ denote the elementary matrix in $\wh{\M}_{\Z \times \Z}$ with $1$ at the $(i,j)$-position and $0$ elsewhere and put
\begin{equation*}
 \wh{E}_{ij}= \sum_{k\in \Z}E_{i+km\, j+kn}\in \wh{\M}_{m \times n}.
\end{equation*}

Let us first describe an $A_{m-1}^{(1)}$-crystal structure on $\wh{\M}_{m \times n}$ for $m \ge 2$. 
Suppose that $A=(a_{ij})_{i,j\in\Z}\in \wh{\M}_{m \times n}$ is given. 
For $i\in \{\,0,1,\dots,m-1\,\}$, we define $\te_i A$ and $\tf_i A$ as follows:

\begin{itemize}
 \item[(1)] Let $\sigma$ be a sequence of $\{\,+,-\,\}$ given by 
\begin{equation*}
\sigma = (\ \cdots,\ \underbrace{-\dots -}_{a_{i+1\, j}},\ \underbrace{+\dots +}_{a_{i\, j}},\ \underbrace{-\dots -}_{a_{i+1\, j+1}},\ \underbrace{+\dots +}_{a_{i\, j+1}},\ \cdots \ ),
\end{equation*}
and let $\td{\sigma}$ be the reduced one, which is well-defined since $\sigma$ has only finitely many $+$'s and $-$'s. 
 
 \item[(2)] If $\td{\sigma}$ has at least one $-$, then we define
\begin{equation*}
 \te_i A =A + \wh{E}_{i j_0} - \wh{E}_{i+1\, j_0},
\end{equation*}
where $j_0$ is the column index of $A$ corresponding the rightmost $-$ in $\td{\sigma}$. If $\td{\sigma}$ has no $-$, then we define $\te_iA={\bf 0}$. Similarly, if $\td{\sigma}$ has at least one $+$, then we define
\begin{equation*}
 \tf_i A =A - \wh{E}_{i j_1} + \wh{E}_{i+1\, j_1},
\end{equation*}
where $j_1$ is the column index of $A$ corresponding the leftmost $+$ in $\td{\sigma}$. If $\td{\sigma}$ has no $+$, then we define $\tf_iA={\bf 0}$.
 
\end{itemize}

Put
\begin{equation}\label{eq:wt A}
\begin{split}
 {\rm wt}^0(A) & = 
 \sum_{k\in\Z}\left(\sum_{j=1}^na_{1+km\, j}\right)(\epsilon_1-k\de)+\dots+
 \sum_{k\in\Z}\left(\sum_{j=1}^na_{m+km\, j}\right)(\epsilon_m-k\de)
 \in P^0,\\
\end{split}
\end{equation}
\begin{equation}\label{eq:wt_0 A}
\begin{split}
 {\rm wt}^0_{\rm cl}(A) & = {\rm cl}({\rm wt}^0(A))=
 \left(\sum_{j\in\Z}a_{1j}\right){\rm cl}(\epsilon_1)+\dots+\left(\sum_{j\in\Z}a_{mj}\right){\rm cl}(\epsilon_m)\in P_{\rm cl}^0,
\end{split}
\end{equation}
\begin{equation*}
\varepsilon_i(A)=\max\{\,k\,|\,\te_i^k A \neq {\bf 0}\,\}, \quad \varphi_i(A)=\max\{\,k\,|\,\tf_i^k A \neq {\bf 0}\,\},\quad (i\in I).
\end{equation*} 
Both $\varepsilon_i(A)$ and $\varphi_i(A)$ are finite since $\sigma$ has only finitely many $+$'s and $-$'s. Moreover, we have
\begin{equation*}
\varphi_i(A) - \varepsilon_i(A) = \sum_{j\in\Z}a_{i+1 \, j} - \sum_{j\in\Z}a_{ij} = \langle {\rm wt}^0(A),h_i \rangle = \langle {\rm wt}^0_{\rm cl}(A),h_i \rangle.
\end{equation*}
Hence we have the following lemma.

\begin{lem}
The set $\wh{\M}_{m \times n}$ is a normal $A_{m-1}^{(1)}$-crystal with $P^0$ and $P^0_{\rm cl}$-weights with respect to ${\rm wt}, \te_i, \tf_i$ for $i\in \{\,0,1,\dots,m-1\,\}$, where ${\rm wt}={\rm wt}^0$ and ${\rm wt}^0_{\rm cl}$, respectively.
\end{lem}

\begin{ex}{\rm
Let $A$ be the generalized affine permutation in Example \ref{ex:A}. For $i = 2$, the associated sequence $\sigma$ and its reduced one $\td{\sigma}$ are
\begin{align*}
\sigma & = (+,-\,,-\,,+\,,-\,,+\,,+\,,-\,), \\
\td{\sigma} & = (\, \cdot\ ,\ \cdot\ ,-\,,\ \cdot\ ,\ \cdot\ ,+\,,\ \cdot\ ,\ \cdot\, ).
\end{align*}
Hence the cell $(2, 8)$ is the position corresponding to the leftmost $+$ in $\td{\sigma}$, and $\td{f}_i A = A - \wh{E}_{2\,8} + \wh{E}_{3\,8}$. {The followings are the submatrices of $A$ and $\tf_2A$ with the row indices in $[4]$.}
\vskip 2mm

\begin{center}
\begin{tikzpicture}[every node/.style={font=\small, scale=1.0}]
\matrix (M)[matrix of math nodes,nodes in empty cells,nodes={rectangle,minimum height=1.0em,minimum width=1.0em,inner sep=0pt,anchor=center,align=center}]
{
& &&&&& &&&&& &&&&& &&&&& &&&&& &&&&& &\\
& &&&&& &1&&&& &&&1&& &&&&& &&&&& &&&&& &\\
& &&&&& 1&&&&& &&1&&& &&\red{2}&&& &&&&& &&&&& &\\
& &&&&& &1&&1&& &&&&& 1&&&1&& &&&&& &&&&& &\\
& &&&&& &&&&& &&&&& &1&&&1& &&&&& &&&&& &\\
& &&&&& &&&&& &&&&& &&&&& &&&&& &&&&& &\\
& &&&&& &&&&& &&&&& &&&&& &&&&& &&&&& &\\
};
\draw[-latex,thick,gray](M-1-1.south west)--(M-1-32.south east);
\node (j) at (M-2-33) {$j$};
\draw[dotted,lightgray](M-2-1.south west)--(M-2-32.south east);
\draw[dotted,lightgray](M-3-1.south west)--(M-3-32.south east);
\draw[dotted,lightgray](M-4-1.south west)--(M-4-32.south east);
\draw[solid,lightgray](M-5-1.south west)--(M-5-32.south east);
\draw[solid,lightgray](M-1-1.north east)--(M-6-1.south east);
\draw[dotted,lightgray](M-1-2.north east)--(M-6-2.south east);
\draw[dotted,lightgray](M-1-3.north east)--(M-6-3.south east);
\draw[dotted,lightgray](M-1-4.north east)--(M-6-4.south east);
\draw[dotted,lightgray](M-1-5.north east)--(M-6-5.south east);
\draw[solid,lightgray](M-1-6.north east)--(M-6-6.south east);
\draw[dotted,lightgray](M-1-7.north east)--(M-6-7.south east);
\draw[dotted,lightgray](M-1-8.north east)--(M-6-8.south east);
\draw[dotted,lightgray](M-1-9.north east)--(M-6-9.south east);
\draw[dotted,lightgray](M-1-10.north east)--(M-6-10.south east);
\draw[-latex,thick,gray](M-1-11.north east)--(M-6-11.south east);
\node (i) at (M-7-12) {$i$};
\draw[dotted,lightgray](M-1-12.north east)--(M-6-12.south east);
\draw[dotted,lightgray](M-1-13.north east)--(M-6-13.south east);
\draw[dotted,lightgray](M-1-14.north east)--(M-6-14.south east);
\draw[dotted,lightgray](M-1-15.north east)--(M-6-15.south east);
\draw[solid,lightgray](M-1-16.north east)--(M-6-16.south east);
\draw[dotted,lightgray](M-1-17.north east)--(M-6-17.south east);
\draw[dotted,lightgray](M-1-18.north east)--(M-6-18.south east);
\draw[dotted,lightgray](M-1-19.north east)--(M-6-19.south east);
\draw[dotted,lightgray](M-1-20.north east)--(M-6-20.south east);
\draw[solid,lightgray](M-1-21.north east)--(M-6-21.south east);
\draw[dotted,lightgray](M-1-22.north east)--(M-6-22.south east);
\draw[dotted,lightgray](M-1-23.north east)--(M-6-23.south east);
\draw[dotted,lightgray](M-1-24.north east)--(M-6-24.south east);
\draw[dotted,lightgray](M-1-25.north east)--(M-6-25.south east);
\draw[solid,lightgray](M-1-26.north east)--(M-6-26.south east);
\draw[dotted,lightgray](M-1-27.north east)--(M-6-27.south east);
\draw[dotted,lightgray](M-1-28.north east)--(M-6-28.south east);
\draw[dotted,lightgray](M-1-29.north east)--(M-6-29.south east);
\draw[dotted,lightgray](M-1-30.north east)--(M-6-30.south east);
\draw[solid,lightgray](M-1-31.north east)--(M-6-31.south east);
\node[minimum width=3em,anchor=east] at (M.west) {$A=$};
\end{tikzpicture}

\begin{tikzpicture}[every node/.style={font=\small, scale=1.0}]
\matrix (M)[matrix of math nodes,nodes in empty cells,nodes={rectangle,minimum height=1.0em,minimum width=1.0em,inner sep=0pt,anchor=center,align=center}]
{
& &&&&& &&&&& &&&&& &&&&& &&&&& &&&&& &\\
& &&&&& &1&&&& &&&1&& &&&&& &&&&& &&&&& &\\
& &&&&& 1&&&&& &&1&&& &&1&&& &&&&& &&&&& &\\
& &&&&& &1&&1&& &&&&& 1&&\red{1}&1&& &&&&& &&&&& &\\
& &&&&& &&&&& &&&&& &1&&&1& &&&&& &&&&& &\\
& &&&&& &&&&& &&&&& &&&&& &&&&& &&&&& &\\
& &&&&& &&&&& &&&&& &&&&& &&&&& &&&&& &\\
};
\draw[-latex,thick,gray](M-1-1.south west)--(M-1-32.south east);
\node (j) at (M-2-33) {$j$};
\draw[dotted,lightgray](M-2-1.south west)--(M-2-32.south east);
\draw[dotted,lightgray](M-3-1.south west)--(M-3-32.south east);
\draw[dotted,lightgray](M-4-1.south west)--(M-4-32.south east);
\draw[solid,lightgray](M-5-1.south west)--(M-5-32.south east);
\draw[solid,lightgray](M-1-1.north east)--(M-6-1.south east);
\draw[dotted,lightgray](M-1-2.north east)--(M-6-2.south east);
\draw[dotted,lightgray](M-1-3.north east)--(M-6-3.south east);
\draw[dotted,lightgray](M-1-4.north east)--(M-6-4.south east);
\draw[dotted,lightgray](M-1-5.north east)--(M-6-5.south east);
\draw[solid,lightgray](M-1-6.north east)--(M-6-6.south east);
\draw[dotted,lightgray](M-1-7.north east)--(M-6-7.south east);
\draw[dotted,lightgray](M-1-8.north east)--(M-6-8.south east);
\draw[dotted,lightgray](M-1-9.north east)--(M-6-9.south east);
\draw[dotted,lightgray](M-1-10.north east)--(M-6-10.south east);
\draw[-latex,thick,gray](M-1-11.north east)--(M-6-11.south east);
\node (i) at (M-7-12) {$i$};
\draw[dotted,lightgray](M-1-12.north east)--(M-6-12.south east);
\draw[dotted,lightgray](M-1-13.north east)--(M-6-13.south east);
\draw[dotted,lightgray](M-1-14.north east)--(M-6-14.south east);
\draw[dotted,lightgray](M-1-15.north east)--(M-6-15.south east);
\draw[solid,lightgray](M-1-16.north east)--(M-6-16.south east);
\draw[dotted,lightgray](M-1-17.north east)--(M-6-17.south east);
\draw[dotted,lightgray](M-1-18.north east)--(M-6-18.south east);
\draw[dotted,lightgray](M-1-19.north east)--(M-6-19.south east);
\draw[dotted,lightgray](M-1-20.north east)--(M-6-20.south east);
\draw[solid,lightgray](M-1-21.north east)--(M-6-21.south east);
\draw[dotted,lightgray](M-1-22.north east)--(M-6-22.south east);
\draw[dotted,lightgray](M-1-23.north east)--(M-6-23.south east);
\draw[dotted,lightgray](M-1-24.north east)--(M-6-24.south east);
\draw[dotted,lightgray](M-1-25.north east)--(M-6-25.south east);
\draw[solid,lightgray](M-1-26.north east)--(M-6-26.south east);
\draw[dotted,lightgray](M-1-27.north east)--(M-6-27.south east);
\draw[dotted,lightgray](M-1-28.north east)--(M-6-28.south east);
\draw[dotted,lightgray](M-1-29.north east)--(M-6-29.south east);
\draw[dotted,lightgray](M-1-30.north east)--(M-6-30.south east);
\draw[solid,lightgray](M-1-31.north east)--(M-6-31.south east);
\node[minimum width=3em,anchor=east] at (M.west) {$\td{f}_2A=$};
\end{tikzpicture}
\end{center}

} 
\end{ex}

Next, we define an $A_{n-1}^{(1)}$-crystal structure on $\wh{\M}_{m \times n}$ for $n \ge 2$, say ${\rm wt}^t, \varepsilon^t_i, \varphi^t_i, \te^t_i, \tf^t_i$ for $i \in \{\,0,1,\dots,n-1\,\}$, by applying the $A_{n-1}^{(1)}$-crystal structure on $\wh{\M}_{n \times m}$ to the transpose $A^t$ of $A\in \wh{\M}_{m \times n}$
\begin{gather*}
{\rm wt}^t(A)={\rm wt}(A^t),\\ \te_i^t A = (\te_iA^t)^t,\quad  \tf_i^t A = (\tf_iA^t)^t,\\
 \varepsilon^t_i(A)= \varepsilon_i(A^t),\quad \varphi^t_i(A)= \varphi_i(A^t).
\end{gather*}

\begin{prop}\label{prop:bocrystal}
The operators $\te_i$ and $\tf_i$ for  $i\in \{\,0,1,\dots,m-1\,\}$ commute with $\te_j^t$ and $\tf_j^t$ for $j\in \{\,0,1,\dots,n-1\,\}$ on $\wh{\mc{M}}_{m\times n}\cup\{{\bf 0}\}$.
\end{prop}
\pf Let $A \in \wh{\mc{M}}_{m\times n}$ be given. Suppose that for each $1 \le i \le m$, we have $a_{ij}=0$ unless $1 \le j \le n$. Then $A$ can be viewed as an $m \times n$ matrix and it is well known that the proposition holds for $A$ when $i, j \neq 0$ (see for examples \cite[Lemma 3.4]{K07} or \cite[Lemma 1.4.7]{vL}). For an arbitrary $A \in \wh{\mc{M}}_{m\times n}$ and $i, j$, we may apply the same argument.
\qed

\begin{rem}\label{rem:bicrystal}
{\rm
If we assume that 
\begin{equation*}
{\rm wt}(A)={\rm wt}^0_{\rm cl}(A){,} \quad {\rm wt}^t(A)={\rm wt}^0_{\rm cl}(A^t),
\end{equation*}
for $A\in \wh{\mc{M}}_{m\times n}$, then $\wh{\M}_{m \times n}$ is an {\em $(A_{m-1}^{(1)},A_{n-1}^{(1)})$-bicrystal} or $(A_{m-1}^{(1)}\times A_{n-1}^{(1)})$-crystal. In other words, 
$\td{x}_i$ and $\td{y}_j^t$ on $\wh{\mc{M}}_{m\times n}\cup\{{\bf 0}\}$ are strict morphisms of $A_{n-1}^{(1)}$-crystals and $A_{m-1}^{(1)}$-crystals, respectively for $x,y\in \{\,e,f\,\}$, $i\in \{\,0,1,\dots,m-1\,\}$, and $j\in \{\,0,1,\dots,n-1\,\}$.
On the other hand, if we assume that ${\rm wt}={\rm wt}^0$ on $\wh{\mc M}_{m\times n}$ as an $A_{m-1}^{(1)}$-crystal, then $\te_j^t, \tf_j^t$ preserves the weights only when $j\neq 0$. 
In case of $\te_0^t, \tf_0^t$, we have 
\begin{equation*}
 {\rm wt}(\te_0^tA)={\rm wt}(A)-\de,\quad {\rm wt}(\tf_0^tA)={\rm wt}(A)+\de,
\end{equation*}
for $A\in \wh{\mc{M}}_{m\times n}$. The same holds for $\te_i$ and $\tf_i$.
}
\end{rem}

\subsection{Crystals of level zero extremal weight modules}
Let $n \ge 2$. In this subsection, we describe a structure of $A_{n-1}^{(1)}$-crystal with $P^0$-weights on $\mc{B}(\la)=\mc{B}_n(\la)$ for $\la\in \cP_n$.

First, consider $\mc{B}((1))=SST_\Z((1))$, which can be identified with $\Z$.
For $\boxed{k}\in \mc{B}((1))$ and $i\in I=\{\,0,1,\dots,n-1\,\}$, we define
\begin{gather*}
\text{${\rm wt}^0(\boxed{k})= \epsilon_{r}-q\delta$,\quad  where $k=nq+r$ for $q\in \Z$ and $1\le r\le n$},\\
 \te_i \,\boxed{k} =
 \begin{cases}
 \boxed{k-1} & \text{if $k\equiv i+1\pmod n$},\\
 {\bf 0} & \text{otherwise},
 \end{cases}\quad\
  \tf_i\, \boxed{k} =
 \begin{cases}
 \boxed{k+1} & \text{if $k\equiv i\pmod n$},\\
 {\bf 0} & \text{otherwise}.
 \end{cases}
\end{gather*}
Then $\mc{B}((1))$ is a normal $A_{n-1}^{(1)}$-crystal with $P^0$-weights with respect to ${\rm wt}^0$, $\te_i$, $\tf_i$ for $i\in I$. 
 
Next, consider $\mc{B}((1^b))$ for $2\le b\le n$. Let $T\in \mc{B}((1^b))$ be given with entries $T(1)<T(2)<\dots<T(b)$. For $i\in I$, we define $\te_iT$ as follows:
\begin{itemize}
 \item[(1)] If $T(k)\equiv i+1 \pmod{n}$ and $T(k-1)<T(k)-1$ for some $2\le k\le b$, then we define $\te_i T$ to be the tableau obtained by replacing $T(k)$ with $T(k)-1$.
 
 \item[(2)] If $T(1)\equiv i+1 \pmod{n}$ and $T(b)< T(1)-1+n$, then we define $\te_i T$ to be the tableau obtained by replacing $T(1)$ with $T(1)-1$.  
 
 \item[(3)] Otherwise, we define $\te_i T={\bf 0}$.
 
\end{itemize}
We define $\tf_i T$ in a similar way. We put ${\rm wt}^0(T)=\sum_{k=1}^b {\rm wt}^0(T(k))$, where we regard $T(k)$ as an element in $\mc{B}((1))$.
 
\begin{prop}\label{eq:crystal of level 0 fundamental}
For $1\le b\le n$, $\mc{B}((1^b))$ is a connected normal $A_{n-1}^{(1)}$-crystal with $P^0$-weights. For $1\le b\le n-1$, $\mc{B}((1^b))$ is isomorphic to $B(\varpi_b)$.
\end{prop}
\pf It is easy to see that $\mc{B}((1^b))$ is a connected normal $A_{n-1}^{(1)}$-crystal. In particular, $\mc{B}((1^n))$ is the crystal {consisting} of single element of weight $0$.

Suppose that $1\le b\le n-1$. Let us keep the notations in Example \ref{ex:crystal of level 1 fundamental weight}. For $i\in \Z$, let ${\bf e}_i=(\de_{ik})_{k\in \Z}$. 
Let $T\in \mc{B}((1^b))$ be given with entries $T(1)<\dots<T(b)$. 
Let $k$ be the positive integer such that $1\le k\le b-1$ and
$$T(1)<\dots<T(k)\le 0< T(k+1)<\dots <T(b).$$ 
In particular, $k=0$ if $0<T(1)$ and $k=b$ if $T(b)\le 0$.
Let 
\begin{equation*}
\begin{split}
 {\bf v}_1 &= {\bf v}_{\La_k} + {\bf e}_{T(k+1)}+\dots+ {\bf e}_{T(b)},\\
 {\bf v}_2 &= {\bf v}_{-\La_k} + {\bf e}_{T(1)}+\dots+ {\bf e}_{T(k)}.
\end{split}
\end{equation*}
Recall that $u_{(1^b)} \in  \mc{B}((1^b))$ is the column tableau having $i$ in its $i$-th row. Then the map sending $T$ to ${\bf v}_1\ot {\bf v}_2$ gives a strict embedding 
\begin{equation*}
 \xymatrixcolsep{3pc}\xymatrixrowsep{0.5pc}\xymatrix{
\mc{B}((1^b)) \ \ar@{->}[r]  & \ \ B(\La_b)\ot B(-\La_0)},
\end{equation*}
and it maps $u_{(1^b)}$ to ${\bf v}_{\La_b}\ot {\bf v}_{-\La_0}$.
This implies that $\mc{B}((1^b))$ is isomorphic to the connected component of ${\bf v}_{\La_b}\ot {\bf v}_{-\La_0}$ in $B(\La_b)\ot B(-\La_0)$. Hence $\mc{B}((1^b))$ is isomorphic to $B(\varpi_b)$ by \cite{Kas94'} and \cite[Proposition 5.4]{Kas02}.
\qed\newline

Recall the bijection $\tau$ on $\mc{B}((1^b))$, defined in Section\ref{subsec:aff.tableaux}. Then we have the following.
\begin{cor}\label{cor:tau on fundamental}
For $T\in \mc{B}((1^b))$, we have $\tau(\te_i T)=\te_i\tau(T)$ and $\tau(\tf_i T)=\tf_i\tau(T)$ for $i\in I$, where ${\rm wt}^0(\tau(T))={\rm wt}^0(T)-\de$. 
\end{cor}
\pf It follows directly from the crystal structure on $\mc{B}((1^b))$.
\qed\newline

Suppose that $R = (a^b)\in \cP_n$ is given for $a\ge 1$ and $1\le b\le n$. Define a map 
\begin{equation*}
\xymatrixcolsep{3pc}\xymatrixrowsep{0.5pc}\xymatrix{
\psi:\ \ \mc{B}(R) \ \ar@{->}[r]  & \ \
 \mc{B}((1^b))^{\ot a} \\
\ \ T =(T^a,\dots, T^1)\  \ar@{|->}[r]  &\ \ T^1\ot \cdots \ot T^a }.
\end{equation*}
Let $u_{R}$ denote the unique element of $\mc{B}(R)$ such that the $i$-th row is filled with $i$ for $1\le i\le b$.
 
\begin{lem}\label{lem:crystal B(R)}
Under the above hypothesis, we have the following.
\begin{itemize}
\item[(1)] The image of $\mc{B}(R)\cup\{{\bf 0}\}$ under $\psi$ is {closed} under $\te_i$ and $\tf_i$ for $i\in I$. Hence $\mc{B}(R)$ is an $A_{n-1}^{(1)}$-crystal with $P^0$-weights such that $\psi$ is a strict embedding.
 
\item[(2)] The connected component of $u_{R}$ in $\mc{B}(R)$ is $\mc{B}(R)_0$ when $b<n$.
 
\item[(3)] The map 
\begin{equation*}
\xymatrixcolsep{3pc}\xymatrixrowsep{0.5pc}\xymatrix{
\ \mc{B}(R)_0\times \cP_{a-1}  \ \ar@{->}[r]  & \ \ \mc{B}(R)}
\end{equation*}
in Lemma \ref{lem:S0 and S} is an isomorphism of $A_{n-1}^{(1)}$-crystals with $P^0$-weights. 
Here we understand $\mc{B}(R)_0\times \cP_{a-1}$ as $\mc{B}(R)_0\ot \cP_{a-1}$ and $\cP_{a-1}$ as an $A_{n-1}^{(1)}$-crystal with $P^0$-weights, where ${\rm wt}(\nu)=-|\nu|\de$, $\te_i\nu=\tf_i\nu={\bf 0}$, and $\varepsilon_i(\nu)=\varphi_i(\nu)=0$ for $\nu\in\cP_{a-1}$ and $i\in I$.
 
\end{itemize}
\end{lem}
\pf (1) Since {$\mc{B}((1^b)) \cup \{{\bf 0}\}$} is {closed} under $\te_i$ and $\tf_i$ for $i\in I$, it is enough to show that $\te_i T$ and $\tf_i T$ are semistandard {for $T \in \mc{B}(R)$}. The proof is similar to the case of type $A_{n-1}$ (cf.~\cite{HK}).

(2) Let $T=(T^a,\dots,T^1)\in \mc{B}(R)_0$ be given. 

Suppose that $\te_i T\not\in \mc{B}(R)_0\cup\{{\bf 0}\}$ for some $i$. 
Then $\te_i T = (\dots, \te_i T^j,\dots)$ for some $1\le j\le a$, and $\te_i T^j$ is obtained from $T^j$ by replacing $T^j(k)$ with $T^j(k)-1$ for some $1\le k\le b$.
Since $T\in \mc{B}(R)_0$ but $\te_i T\not\in \mc{B}(R)_0\cup\{{\bf 0}\}$, at least one of $(\te_i T^j,\tau^{-1}(T^{j-1}))$ or $(T^{j+1},\tau^{-1}(\te_i T^j))$ is semistandard. Suppose that $(\te_i T^j,\tau^{-1}(T^{j-1}))$ is semistandard. Considering $\te_i T \in \mc{B}(R)$,  
it is straightforward to see that $(T^j,\tau^{-1}(T^{j-1}))$ is also semistandard, which is a contradiction. For the other cases, we have similar contradiction. By the same arguments, we also have $\tf_i T\in \mc{B}(R)_0\cup\{{\bf 0}\}$. Hence {$\mc{B}(R)_0\cup \{{\bf 0}\}$ is closed} under $\te_i$ and $\tf_i$ for $i \in I$.

Conversely, we claim that any $T\in \mc{B}(R)_0$ is connected to $u_R$.
Let $t=T^a(1)$. We first show that $T\in \mc{B}(R)_0$ is connected to $T_0=u_R^{(t)}$, where $u_{R}^{(t)}$ is the element of $\mc{B}(R)$ such that the $i$-th row is filled with $t+i-1$ for $1\le i\le b$. 
Suppose that $T\neq T_0$. Let $d(T)=\sum_{j=1}^a\sum_{k=1}^b(T^j(k)-T_0^j(k))\ge 0$. 
Since $T\in \mc{B}(R)_0$ and {$\mc{B}(R)_0\cup \{{\bf 0}\}$ is closed} under $\te_i$ and $\tf_i$ for $i\in I$, we have $0 \le d(\te_iT)< d(T)$ for $i$ such that $\te_iT\neq {\bf 0}$.
Note that there exists at least one $i$ such that $\te_iT\neq {\bf 0}$. 
For example, choose the smallest $k$ such that $T^1(k)\neq T^1_0(k)$ and then $i=T^1(k)-1$. By induction on $d(T)$, we conclude that $T$ is connected to $T_0$. 

Next, we have
\begin{equation*}
\begin{cases}
 \te_{t+b-2}^a\dots \te_{t}^a\te_{t-1}^a T_0=u_R^{(t-1)} & (t>1),\\
 \tf_{t}^a\dots \tf_{t+b-2}^a\tf_{t+b-1}^a T_0=u_R^{(t+1)} & (t<1).
\end{cases}
\end{equation*}
Repeating this step, we conclude that $T_0$ is connected to $u_R=u_R^{(0)}$.
Therefore, $\mc{B}(R)_0$ is the connected component of $u_{R}$ in $\mc{B}(R)$.

(3) It follows from Corollary \ref{cor:tau on fundamental}.
\qed

\begin{rem}\label{rem:trivial rectangle}
{\rm
Suppose that $b=n$. In this case, $\mc{B}(R)_0=\{\,\uptau^t (u_R)\,|\,t\in\Z\,\}$, where each $\{\,\uptau^t (u_R)\,\}$ forms a trivial crystal of weight $-at\de$.
} 
\end{rem}

\begin{prop}
For $R=(a^b)$ with $a\ge1$ and $1\le b<n$, $\mc{B}(R)$ is isomorphic to $B(a\varpi_b)$ as an $A_{n-1}^{(1)}$-crystal with $P^0$-weights.
\end{prop}
\pf It follows from Proposition \ref{eq:crystal of level 0 fundamental}, Lemma \ref{lem:crystal B(R)}, and \cite[Theorem 4.16(a)]{BN}.
\qed 

\begin{rem} 
{\rm
The isomorphism in Lemma \ref{lem:crystal B(R)}(3) is also proved in \cite[Theorem 2]{NS} for an affine Lie algebra $\mf{g}$ in terms of Lakshmibai-Seshadri paths.} 
\end{rem}

Let $\la\in \cP_n$ be given and let $\mc{B}(\la)$ as given in \eqref{eq:B(la)}. We regard $\mc{B}(\la)$ as an $A_{n-1}^{(1)}$-crystal with $P^0$-weights by identifying 
\begin{equation*}
\mc{B}(\la)= \mc{B}(R_1) \otimes \dots \otimes \mc{B}(R_l). 
\end{equation*}
Let 
\begin{equation*}
\varpi_\la=m_1\varpi_1+\dots+m_{n-1}\varpi_{n-1}, 
\end{equation*}
where $m_i$ is the multiplicity of $i$ in $\la$.
By \cite[Theorem 4.16]{BN} and \cite[Conjectures 13.1, 13.2]{Kas02} (see also \cite[Remark 4.17]{BN}), we have the following.
\begin{prop}\label{prop:extremal crystal}ove
 For $\la\in\cP_{n-1}$, $\mc{B}(\la)$ is isomorphic to $B(\varpi_\la)$.
\end{prop} 

\begin{rem}{\rm
(1) In \cite{I}, another proof of Proposition \ref{prop:extremal crystal} is given using the standard monomial theory for semi-infinite Lakshmibai–Seshadri paths \cite{INS}.

(2) Let $u_{\la}=u_{R_1}\ot{\cdots}\ot u_{R_l}$.
Then $\mc{B}(\la)_0=\mc{B}(R_1)_0 \otimes \dots \otimes \mc{B}(R_l)_0$ is the connected component of $u_\la$ in $\mc{B}(\la)$ (see \cite[Remark 4.17]{BN}). In case of $A_{n-1}^{(1)}$, we can prove it directly using the crystal structure described here.

(3) Suppose that $\la\in \cP_n$ with $n=\ell(\la)$. We have $\mc{B}(\la)=\mc{B}(\mu)\ot \mc{B}(R_n)$, where $\mu=(\la_1,\dots,\la_{n-1})$ (see Remark \ref{rem:trivial rectangle}).

(4) We may also regard $\mc{B}(\la)$ as an $A_{n-1}^{(1)}$-crystal with $P^0_{\rm cl}$-weights by using ${\rm wt}^0_{\rm cl}={\rm cl}({\rm wt}^0_{\rm cl})$.
}
\end{rem}
 
\begin{rem}\label{rem:finite crystal}
{\rm 
Let $R = (a^b)$. Then the set {$\mc{B}(R)_0/\Z \cup \{{\bf 0}\}$ is closed} under $\te_i$ and $\tf_i$ for $i\in I$ since $\uptau$ commutes with $\te_i$ and $\tf_i$, and ${\rm wt}^0_{\rm cl}$ is constant on any equivalence class. 
Hence by \eqref{eq:CSST and S0}, we may regard  $CSST_{[n]}(R)$ as a normal $A_{n-1}^{(1)}$-crystal with $P^0_{\rm cl}$-weights. Then $CSST_{[n]}(R)$ is isomorphic to $SST_{[n]}((1^b))^{\ot a}$. Note that $SST_{[n]}((1^b))=CSST_{[n]}((1^b))$ and it is isomorphic to the crystal of Kirillov-Reshetikhin module or level one perfect crystal associated to ${\rm cl}(\varpi_b)$ \cite{KMN2}.

Furthermore, if $(r_1,\dots,r_{a-1})$ is the offset vector of $T\in CSST_{[n]}(R)$, then $-r_j$ is equal to the value of the local energy function on $SST_{[n]}((1^b))^{\ot 2}$ at $(T^{j+1},T^j)$ or $T^j\ot T^{j+1}$ for $1\le j\le a-1$ (see \cite[Section 3]{NaYa}).

In general, for $\la\in \cP_{n-1}$ we regard as an $A_{n-1}^{(1)}$-crystal with $P^0_{\rm cl}$-weights
\begin{equation*}
 CSST_{[n]}(\la) = CSST_{[n]}(R_1) \otimes \dots \otimes CSST_{[n]}(R_l).
\end{equation*}
Hence if we regard $\mc{B}(\la)$ as an $A_{n-1}^{(1)}$-crystal with $P^0_{\rm cl}$-weights, then the bijection \eqref{eq:CSST and B general}
\begin{equation*}
\xymatrixcolsep{3pc}\xymatrixrowsep{0.5pc}\xymatrix{
\ CSST_{[n]}(\la) \times \mc{P}(\la)  \ \ar@{->}[r]  & \ \
 \mc{B}(\la)}
\end{equation*} 
is an isomorphism, where the left-hand side is understood as a disjoint union of $CSST_{[n]}(\la)$ parametrized by $\mc{P}(\la)$.

} 
\end{rem}

\subsection{Isomorphism of crystals}

Let 
\begin{equation}\label{eq:T mn}
\mc{T}_{m\times n}= \bigsqcup_{\la\in \cP_m\cap\cP_n} 
CSST_{[m]}(\la)\times \mc{B}_n(\la).
\end{equation}
We first assume that $\mc{T}_{m\times n}$ is a normal $A_{m-1}^{(1)}$-crystal with $P^0_{\rm cl}$-weights with respect to $\te_i$, $\tf_i$ for $i\in \{\,0,1,\dots,m-1\,\}$ for $m \ge 2$, where $\te_i$, $\tf_i$ are the Kashiwara operators on $CSST_{[m]}(\la)$ (see Remark \ref{rem:finite crystal}).
Also we assume that $\mc{T}_{m\times n}$ is a normal $A_{n-1}^{(1)}$-crystal with $P^0$-weights with respect to $\te^t_j$, $\tf^t_j$ for $j\in \{\,0,1,\dots,n-1\,\}$ for $n \ge 2$, where $\te^t_j$, $\tf^t_j$ denote the Kashiwara operators on $\mc{B}_n(\la)$.

The following is the second main result in this paper. 
The proof is given in Section \ref{sec:proof of main-2}.

\begin{thm}\label{thm:main-2}
The bijection 
\begin{equation*}
\xymatrixcolsep{3pc}\xymatrixrowsep{0.5pc}\xymatrix{
\kappa:\ \ \wh{\M}_{m\times n} \ \ar@{->}[r]  & \ \ \mc{T}_{m\times n}}
\end{equation*}
commutes with $\te_i$, $\tf_i$ for $i\in \{\,1,\dots,m-1\,\}$ and $\te^t_j$, $\tf^t_j$ for $j\in \{\,0,1,\dots,n-1\,\}$.
\end{thm}

We remark that the map $\kappa$ does not commute with $\te_0$ and $\tf_0$, but $\kappa_1:=\pi_1\circ\kappa$ does, where $\pi_1$ is the projection of $\mc{T}_{m\times n}$ along the first component (see Remark \ref{rem:e_0 and f_0}).
Since ${\rm wt}^0_{\rm cl}(A)={\rm wt}^0_{\rm cl}(P_0)$, $\kappa_1$ induces the following. 

\begin{cor}\label{cor:row equiv}
A generalized affine permutation $A\in \wh{\M}_{m\times n}$ is $A_{m-1}^{(1)}$-crystal equivalent to $P_0$, where $\kappa(A)=(P_0,Q)$. 
\end{cor}

Note that $\kappa$ is an isomorphism of $A_{n-1}^{(1)}$-crystals with $P^0_{\rm cl}$-weights, but it does not preserve $P^0$-weights. More precisely, for $A\in \wh{\M}_{m\times n}$ with $\kappa(A)=(P_0,Q)$, we see from the definition of $\kappa$ that
\begin{equation*}
{\rm wt}^t(A)= {\rm wt}^0(A^t) = {\rm wt}^0(Q) - \left(\sum_{i=1}^l|\eta^{(i)}| \right) \de,
\end{equation*}
where $\eta^{(i)}$ is the symmetrized offset vectors of $P_0^{(i)}$ in Section \ref{subsec:affine RSK}. 
So in order to have a morphism of $A_{n-1}^{(1)}$ crystals with $P^0$-weights, we may modify the weight function on $\wh{\M}_{m\times n}$ in $P^0$ by 
\begin{equation}\label{eq:new aff wt}
{\rm\bf wt}^t(A)= {\rm wt}^0(A^t) + H_m(A)\de,
\end{equation}
where $H_m(A)=\sum_{i=1}^l|\eta^{(i)}|$. 
Then we have the following.
\begin{cor}\label{cor:decomp M into An crystal}
If we regard $\wh{\M}_{m\times n}$ as an $A_{m-1}$-crystal with $P^0_{\rm cl}$-weights and as an $A_{n-1}^{(1)}$-crystal with $P^0$-weights with respect to ${\rm\bf wt}^t$, then it is an $(A_{m-1},A_{n-1}^{(1)})$-bicrystal and $\kappa$ is an isomorphism of $(A_{m-1},A_{n-1}^{(1)})$-bicrystals. In particular, a generalized affine permutation $A\in \wh{\M}_{m\times n}$ is $A_{n-1}^{(1)}$-crystal equivalent to $Q$, where $\kappa(A)=(P_0,Q)$. 
\end{cor}

We remark that {both $m$ and $n$ do not need to be greater} than $1$ for Theorem \ref{thm:main-2} and its corollaries. In particular, Corollary \ref{cor:row equiv} holds for $n=1$ and Corollary \ref{cor:decomp M into An crystal} holds for $m=1$.

\begin{rem} {\rm
The function $H_m(\cdot)$ in \eqref{eq:new aff wt} is related to the intrinsic energy function on $A^{(1)}_{m-1}$-crystals with $P^0_{\rm cl}$-weights as follows.
Let $T \in CSST_{[m]}(R)$ be given where $R=(a^b)$. Let $r=(r_1, \dots, r_{a-1})$ be the offset vector and $\eta = (\eta_1, \dots, \eta_{a-1})$ the symmetrized offset vector of $T$.  
Then we have
\begin{equation*}
\mc{H}_m(T)= |\eta| - a|r|,
\end{equation*}
where $|r|=r_1+\dots+r_{a-1}$ and $\mc{H}_m(\cdot)$ is the intrinsic energy function on $CSST_{[n]}(R)$ with  $\mc{H}_m(u_R)=0$ (cf.~\cite{Sh} for its definition). 
Hence for $A\in \wh{\M}_{m\times n}$ with $\kappa(A)=(P_0,Q)$, we have
\begin{equation*}
 H_m(A) = \sum_{i=1}^l\left( \mc{H}_m(P_0^{(i)}) + m_i |r^{(i)}| \right),
\end{equation*}
where $R_i=(a_i^{m_i})$ is the shape and $r^{(i)}$ is the offset vector of $P_0^{(i)}$ in Section \ref{subsec:affine RSK}, respectively.
} 
\end{rem}

\section{Proof of Theorem \ref{thm:main-1}}\label{sec:Proof of affine RSK}

\subsection{Standardization and $\kappa_0$}\label{subsec:standardization 2}In this subsection, we show the map $\kappa_0$ in \eqref{eq:Aff RSK 0} is compatible with the standardization \eqref{eq:standardization of A}.

Let $\alpha=(\alpha_1,\dots,\alpha_m)\in \Z_{\geq 0}^m$ and $\beta=(\beta_1,\dots,\beta_n)\in \Z_{\geq 0}^n$ be given such that 
$|\alpha| = |\beta|$, where $|\alpha|=\alpha_1 + \dots + \alpha_m$.
Let
\begin{equation*}
 \wh{\M}_{m \times n}(\alpha,\beta)=\{\,A\,|\,A\in \wh{\M}_{m\times n},\ {\rm row}(A)=\alpha, \ {\rm col}(A)=\beta\,\},
\end{equation*} where ${\rm row}(A)$ and ${\rm col}(A)$ are the row and column contents of $A$ given in \eqref{eq:row and col of A}.
We have $\wh{\M}_{m \times n}=\bigsqcup_{\alpha,\beta}\wh{\M}_{m \times n}(\alpha,\beta)$.

For $K \ge 1$, recall that $A=(a_{ij})_{i, j \in \Z} \in \wh{\M}_{K\times K}$ is an extended affine permutation of $K$ if $\sum_{j \in \Z} a_{ij} = \sum_{i \in \Z}a_{ij} = 1$ for each $i$ and $j$. Let $\td{W}_K\subset \wh{\M}_{K\times K}$ denote the set of extended affine permutations of $K$.
Let $B \in \td{W}_K$ be given. 
We say that $i \in [K]$ is a {\em descent} of $B$ if $(i, j), (i+1, j') \in {\rm supp}(B)$ and $j > j'$. Let $\alpha=(\alpha_1,\dots,\alpha_m)\in \Z_{\geq 0}^m$ and $\beta=(\beta_1,\dots,\beta_n)\in \Z_{\geq 0}^n$ be given with $|\alpha| = |\beta| = K$. 
We say that $B \in \td{W}_K$ is {\em $\alpha$-descending} if for any {$k \in [m]$} and $i$ with 
\begin{equation*}
\alpha_1 + \dots + \alpha_{k-1} < i < \alpha_1 + \dots + \alpha_{k-1} + \alpha_k,
\end{equation*}
$i$ is a descent of $B$, where we understand the empty sum is $0$. We say that $B$ is {\em $(\alpha, \beta)$-descending} if $B$ is $\alpha$-descending and $B^t$ is $\beta$-descending. 
Let $\td{W}_{K,(\alpha,\beta)}$ denote the set of extended affine permutations which are $(\alpha, \beta)$-descending. 

\begin{lem}\label{lem:inverse std}
Under the above hypothesis, we have a bijection
\begin{equation*}
\xymatrixcolsep{3pc}\xymatrixrowsep{0.5pc}\xymatrix{
\wh{\M}_{m \times n }(\alpha,\beta) \ \ar@{->}[r]  & \ \
 \td{W}_{K,(\alpha,\beta)} \\
\ \ A \  \ar@{|->}[r]  &\ \ A^{\tt st} }. 
\end{equation*}
\end{lem}
\pf For $A\in \wh{\M}_{m \times n }(\alpha,\beta)$, we have $A^{\tt st}\in \td{W}_{K,(\alpha,\beta)}$ by Lemma \ref{lem: descending Ast}, hence the map is well-defined.

Let us prove the existence of its inverse.
Let $B=(b_{kl})_{k,l\in\Z} \in \td{W}_{K,(\alpha,\beta)}$ be given. 
Consider the partition $\{\,I_i \times J_j\,\}_{i,j \in \Z}$ of $\Z\times \Z$ associated to $\alpha$ and $\beta$ in \eqref{eq:row and column for A perm}. 
Then we define $A=(a_{ij})_{i, j \in \Z}$ to be a $\Z \times \Z$ matrix where
\begin{equation*}
a_{ij} = \sum_{(k,l) \in I_i \times J_j}b_{kl}
\end{equation*}
Since $B\in \td{W}_{K,(\alpha,\beta)}$, we see from the definition that $I_i\times J_j$ submatrix of $B$ contains an $a_{ij}\times a_{ij}$ matrix with 1 on the antidiagonal and 0 elsewhere. 
Then it is easy to see that $A$ is the unique matrix in $\wh{\M}_{m \times n}(\alpha,\beta)$ such that $A^{\tt st}=B$. This proves the bijectivity.
\qed
\newline

Now, let us consider standardization on the side of tableaux in \eqref{eq:Aff RSK 0}.

Let $\la\in \cP_n$ and $\alpha=(\alpha_1, \dots, \alpha_n) \in \Z^n_+$ be given. Let $CSST_{[n]}(\la)_\alpha$ be the set of $T\in CSST_{[n]}(\la)$ with content $\alpha$, that is, $\alpha_k$ is the number of occurrences of $k \in [n]$ in $T$. 
For $T\in CSST_{[n]}(\la)_\alpha$, we define $T^{\tt st}$ to be a tableau obtained from $T$ by replacing each $k \in [n]$ in $T$ with $\alpha_k\neq 0$ by the consecutive numbers
\begin{equation*}
\alpha_1 + \dots + \alpha_{k-1} + 1 < \dots < \alpha_1 + \dots + \alpha_{k-1} + \alpha_{k} 
\end{equation*}
from left to right.
Let $K=|\la|$. By definition, we have 
\begin{equation*}
T^{\tt st}\in CST_{[K]}(\la)
\end{equation*}
(cf.~ Section \ref{subsec:tableaux}).
We call $T^{\tt st}$ the {\em standardization of $T$}.
The following lemma implies that the symmetrized offset vector of $T$ is invariant under standardization.

\begin{lem}\label{lem:std and msstd}
Let $T \in CSST_{[n]}(R)$ be given where $R=(a^b) \in \cP_n$ for some $a, b \ge 1$. For $\nu \in \cP_{a-1}$, we have
\begin{equation*}
\tau_n^{\nu_{\tt rev}}(T) \in \cB_n(R)_0 \quad\text{if and only if}\quad \tau_K^{\nu_{\tt rev}}(T^{\tt st}) \in \cB_K(R)_0,
\end{equation*} where $K=ab$.
\end{lem}
\pf
Let $T=(T^a, \dots, T^1)$ and $T^{\tt st}=S=(S^a, \dots, S^1)$. 
For $1\le j\le a-1$, let $r_j$ be the smallest integer satisfying $(T^{j+1},\tau_n^{r_j}(T^{j}))$ is semistandard. It suffices to show that $(S^{j+1}, \tau_K^{r_j}(S^j))$ is semistandard but  $(S^{j+1}, \tau_K^{r_j-1}(S^j))$ is not. It is straightforward to see from the definition of $T^{\tt st}$.
\qed\newline

Let $T \in CST_{[K]}(\la)$ be given with $|\la|=K$. We say that $i \in [K]$ is a {\em (column) descent} of $T$ if the letter $i+1$ appears to the right of $i$ in $T$ ($i+1$ lies in the column with a smaller index than that of $i$). For $\alpha=(\alpha_1, \dots, \alpha_n) \in \Z^n_+$ with $K=|\alpha|$, we say that $T \in CST_{[K]}(\la)$ is {\em $\alpha$-descending} if for any $k \in [n]$ and $i$ with
$$\alpha_1 + \dots + \alpha_{k-1} < i < \alpha_1 + \dots + \alpha_{k-1} + \alpha_k,$$
$i$ is a descent of $T$.
Let $CST_{[K],\alpha}(\la)$ be the set of $T\in CST_{[K]}(\la)$ which are $\alpha$-descending. Then we can check the following.

\begin{lem}\label{lem:std of tableaux}
Under the above hypothesis, we have a bijection
\begin{equation*}
\xymatrixcolsep{3pc}\xymatrixrowsep{0.5pc}\xymatrix{
CSST_{[n]}(\la)_\alpha \ \ar@{->}[r]  & \ \
 CST_{[K],\alpha}(\la) \\
\ \ T \  \ar@{|->}[r]  &\ \ T^{\tt st} }. 
\end{equation*}
\end{lem}

The following generalizes a well-known property of the usual RS correspondence (cf. \cite[Section 1.1]{Ful}).
\begin{lem}\label{lem:RSK and descent}
Let $B \in \td{W}_K$ be given with $\kappa_0(B) = (P_0, Q_0, \rho)$. Then $i\in [K]$ is a descent of $B$ if and only if $i$ is a descent of $P_0$. Similarly, $j \in [K]$ is a descent of $B^t$ if and only if $j$ is a descent of $Q_0$.
\end{lem}

\pf  It follows from \cite[Proposition 3.6]{CPY}.
\qed\newline

\begin{cor}\label{cor:affine RS preserves descents}
Let $B \in \td{W}_K$ be given with $\kappa_0(B) = (P_0, Q_0, \rho)$, where $P_0, Q_0 \in CST_{[K]}(\la)$ and $\rho\in \Z_{\geq 0}^{\la_1}$ with $\la\in \cP_m\cap\cP_n$. 
For $\alpha\in \Z_{\geq 0}^m$ and $\beta\in \Z_{\geq 0}^n$ with $|\alpha| = |\beta| = K$,
the following are equivalent:
\begin{itemize}
    \item[(1)] $B\in \td{W}_{K,(\alpha,\beta)}$,
    \item[(2)] $P_0\in CST_{[K],\alpha}(\la)$ and $Q_0 \in CST_{[K],\beta}(\la)$.
\end{itemize}
\end{cor}

It is clear that $\kappa_0$ in \eqref{eq:Aff RSK 0} preserves contents, that is, if $A \in \wh{\M}_{m \times n}(\alpha,\beta)$ with $\kappa_0(A)=(P_0, Q_0, \rho)$, then $P_0 \in CSST_{[m]}(\la)_\alpha$ and $Q_0 \in CSST_{[n]}(\la)_\beta$.
The following proposition shows that $\kappa_0$ is compatible with standardizations.

\begin{prop}\label{prop:affine RSK and standardization}
For $A \in \wh{\M}_{m \times n}$ with $\kappa_0(A) = (P_0, Q_0, \rho)$, we have 
$$\kappa_0(A^{\tt st}) = (P^{\tt st}_0, Q^{\tt st}_0, \rho).$$
\end{prop}

\pf 
Assume that $A\in \wh{\mc M}_{m\times n}(\alpha,\beta)$ for some $\alpha\in\Z_{\geq 0}^m$ and $\beta\in\Z_{\geq 0}^n$ with $|\alpha| = |\beta|=K$.

Let $\kappa_0(A^{\tt st}) = (P_1,Q_1,\varrho)$. 
By Lemma \ref{lem:inverse std} and Corollary \ref{cor:affine RS preserves descents}, we have 
\begin{equation*}
 P_1\in CST_{[K],\alpha}(\la),\quad Q_1\in CST_{[K],\beta}(\la)
\end{equation*}
for some $\la\in \cP$..

We claim that $P_0$ (resp. $Q_0$) is the image of $P_1$ (resp. $Q_1$) under the inverse of the bijection in Lemma \ref{lem:std of tableaux} and $\rho = \varrho$. 
Let $\{\,I_i \times J_j\,\}_{i,j \in \Z}$ be the partition of $\Z\times \Z$ associated to $\alpha$ and $\beta$ in \eqref{eq:row and column for A perm}. 
Let ${\bf s}$ and ${\bf s}^{\tt st}$ be the stream given from the the back-post corners of the {sets of} zig-zags $\{{\bf z}_k\}_{k\in\Z}$ and $\{{\bf z}^{\tt st}_k\}_{k\in\Z}$ associated to $d^{\tt sw}_A$ and $d^{\tt sw}_{A^{\tt st}}$, respectively. 
Since we have $(d^{\tt sw}_A)^{\tt st}=d^{\tt sw}_{A^{\tt st}}$ by \eqref{eq:sw proper A and Ast}, the flow of ${\bf s}$ is equal to that of ${\bf s}^{\tt st}$, and ${\bf s}$ is obtained from ${\bf s}^{\tt st}$ by replacing each $(k, l)$ in ${\bf s}^{\tt st}$ with $(i, j)$ when $(k, l) \in I_i\times J_j$. Since $(A^\flat)^{\tt st}=(A^{\tt st})^\flat$ (see Remark \ref{rem:fw and st commute}), we use induction on $K$ to prove the claim. Hence we have $P^{\tt st}_0 = P_1$ and $Q^{\tt st}_0 = Q_1$.
%
\qed
\newline

\subsection{Dominance condition and $\kappa$}\label{subsec:dominance}
In this subsection, we characterize the image of $\kappa_0$, and then prove Theorem \ref{thm:main-1}.

\subsubsection{Offset constants}\label{subsec:two columns}
Let $\la = (2^l)$ be of rectangular shape with two columns for $1\le l\le \min\{m,n\}$.

Suppose that $P_0 = (P_0^2,P_0^1) \in CSST_{[m]}(\la)$ and $Q_0 = (Q_0^2,Q_0^1) \in CSST_{[n]}(\la)$ are given.
Let $(0,\eta)$ and $(0,\theta)$ be the offset vectors for $P_0$ and $Q_0$, that is,
\begin{equation}\label{eq:eta and theta}
\tau^{(0,\eta)}(P_0)=\left(P_0^2, \tau^{\eta}(P_0^1)\right) \in \mc{B}_m(\la)_0,\quad 
\tau^{(0,\theta)}(Q_0)= \left(Q_0^2, \tau^{\theta}(Q_0^1)\right) \in \mc{B}_n(\la)_0,
\end{equation}
 and put
\begin{equation}\label{eq:offset const}
r = \theta - \eta. 
\end{equation}
We claim that $r$ is equal to the {\em offset constant} introduced in \cite[Section 5.2]{CPY}, which plays an important role in characterizing the image of $\kappa_0$ in {the} case of extended affine permutations.

Write $P_0 = (\ba_1, \ba_2)$ and $Q_0 = (\bb_1, \bb_2)$.
For $(\rho_1, \rho_2) \in \Z^2$, let 
\begin{equation*}
 {\bf s}=\{\,c_i=(a_i, b_i)\,\}_{i\in\Z},\quad 
 {\bf s'}=\{\,c'_i=(a'_i, b'_i)\,\}_{i\in\Z}
\end{equation*}
be the streams of flow $l$, whose defining data are $(\ba_1, \bb_1, \rho_1)$ and $(\ba_2, \bb_2, \rho_2)$, respectively. 
We may assume that 
\begin{equation}\label{eq:a_1,a_2,b_1,b_2}
\begin{split}
 \ba_1=(a_1, \dots, a_l),\quad \bb_1 = (b_{1-\rho_1}, \dots, b_{l-\rho_1}),\\ 
 \ba_2=(a'_1, \dots, a'_l),\quad \bb_2=(b'_{1-\rho_2}, \dots, b'_{l-\rho_2}).  
\end{split}
\end{equation}
Let $k$ be the smallest integer such that $c_{i} \ge_{\tt nw} c'_{i+k}$ for all $i \in \Z$, and then consider the {two streams} 
\begin{equation*}
 {\bf t}=\{\,d_i=(a'_{i+k}, b_i)\,\}_{i\in\Z},\quad 
 {\bf t'}=\{\,d'_i=(a_i, b'_{i+k})\,\}_{i\in\Z}.
\end{equation*}

Let $A=(a_{ij})_{i,j \in \Z} \in \wh{\M}_{m\times n}$, where $a_{ij}$ is the number of occurrence of $(i, j)$ in ${\bf t}\cup{\bf t'}$. 
It is easy to observe that $l$ is the width of $A$ and ${\bf t}$ is the southwest channel of $A$. Let $d$ be a proper numbering on $A$ defined by $d(d_i)=d(d'_i)=i$ for $i \in \Z$.

\begin{lem}\label{lem:offset constant}
Under the above hypothesis, we have
$\rho_2 - \rho_1 \ge r$ if and only if $d=d^{\tt sw}_{A}$.
\end{lem}
\pf We observe that in order to show $d=d^{\tt sw}_{A}$, it suffices to find $i$ such that $d_i >_{\tt NW} d'_{i+1}$ by the condition (2) in Lemma \ref{lem:characterization of ch numbering}. 
Furthermore, $d_{i-1} >_{\tt NW} d'_i$ for some $i$ is equivalent to $a'_{i-1+k} < a_i$ since $b_{i-1} < b'_{i+k}$ is redundant by our choice of $k$.

By \eqref{eq:eta and theta}, the constants $\eta$ and $\theta$ are the smallest integers such that 
\begin{equation*}
 a_i \le a'_{i+\eta},\quad b_i \le b'_{i+\rho_1-\rho_2+\theta}
\end{equation*}
for {all $i \in \Z$} respectively. By the minimality of $\eta$, $\theta$, we have $\eta \le k$ and $\theta \le \rho_2 - \rho_1 + k$.

Suppose that $d=d^{\tt sw}_{A}$. 
Then by the above observation, there exists $i$ with $a_i > a'_{i+k-1}$. 
In this case, we have $k-1 < \eta \le k$ by the minimality of $\eta$. 
Hence $\eta = k$ and $$\rho_2 - \rho_1 \ge \theta - k = \theta - \eta =r.$$
Note that we have $\eta\ge \rho_1-\rho_2 +\theta$, which implies that $b_i\le b'_{i+\eta}$.

Conversely, suppose that $\rho_2 - \rho_1 \ge \theta - \eta$. 
Since $\rho_2 - \rho_1 + \eta \ge \theta$, we have 
\begin{equation}\label{eq:equiv cond for offset const}
 b_i \le b'_{i+\rho_1-\rho_2+(\rho_2-\rho_1+\eta)} = b'_{i+\eta}\quad (i\in \Z).
\end{equation}
Since $a_i\le a'_{i+\eta}$, we have $k\le \eta$ by the minimality of $k$, which implies $a'_{i+k-1} < a_i$ for some $i$. 
Hence $d=d^{\tt sw}_{A}$.
\qed\newline

The following is another characterization of $r$.
\begin{cor}\label{cor:offset constant}
Under the above hypothesis, we have $\rho_2 - \rho_1 \ge r$ if and only if
\begin{center}
$\tau^{(\rho_1,\rho_2 + \eta)}(Q_0) \in \mc{B}_n(\la)$.
\end{center}
\end{cor}
\pf We have shown in the proof of Lemma \ref{lem:offset constant} that $\rho_2-\rho_1\ge r$ is equivalent to $b_i\le b'_{i+\eta}$ \eqref{eq:equiv cond for offset const}. On the other hand, by \eqref{eq:a_1,a_2,b_1,b_2} it is equivalent to
\[
\pushQED{\qed} 
\tau^{(\rho_1,\rho_2 + \eta)}(Q_0)=(\tau^{\rho_1}\bb_1, \tau^{\rho_2 + \eta}\bb_2) \in \mc{B}_n(\la).\qedhere
\popQED
\]

\begin{rem}\label{rem:offser constant CPY}
{\rm
It is not difficult to see that the constants $\eta$ and $\theta$ are equal to the {\em local charges} of $P_0$ and $Q_0$, respectively \cite[Definition 5.3]{CLP}. 
Hence by \cite[Theorem 5.10]{CLP}, the constant $r=\theta-\eta$ in \eqref{eq:offset const} is equal to the offset constant of $(P_0,Q_0)$ in \cite{CPY} when $P_0, Q_0\in CST_{[K]}((2^l))$ for some $K$.
} 
\end{rem}

\subsubsection{Proof of Theorem \ref{thm:main-1}}
We keep the notations in Section \ref{subsec:affine RSK}.
Suppose that $(P_0, Q_0, \rho) \in CSST_{[m]}(\la)\times CSST_{[n]}(\la)\times \Z^{\la_1}$ is given for $\la\in \cP_m\cap\cP_n$.

For $1 \le i \le l$ with $m_i \ge 1$, we let 
\begin{itemize}
\item[$\bullet$] $\rho^{(i)}\in \Z^{m_i}$ : the subsequence of $\rho$ corresponding to the columns of $R_i$,

\item[$\bullet$] $\eta^{(i)}\in \cP_{m_i-1}$ : the symmetrized offset vector of $P^{(i)}_0$,

\item[$\bullet$] $\theta^{(i)}\in \cP_{m_i-1}$ : the symmetrized offset vector of $Q^{(i)}_0$, 

\item[$\bullet$] $\zeta^{(i)} = \theta^{(i)} - \eta^{(i)} \in \Z^{m_i}$ with the last component being zero.
\end{itemize}

We write $\zeta=(\zeta^{(1)},\dots,\zeta^{(l)})$.

\begin{df}\label{eq:dominant condition}
{\rm
Under the above hypothesis, we say that $(P_0, Q_0, \rho)$ is {\em dominant} if \begin{equation*}
\rho_{\tt rev} - \zeta \in \mc{P}(\la).
\end{equation*}
}
\end{df}

\begin{lem}\label{lem:std preserves dominance}
Under the above hypothesis, we have
$(P_0, Q_0, \rho)$ is dominant if and only if $(P^{\tt st}_0, Q^{\tt st}_0, \rho)$ is dominant.
\end{lem}
\pf For $1 \le i \le l$ with $m_i \ge 1$, let $\eta^{(i)}_{\tt st}$ and $\theta^{(i)}_{\tt st}$ be the symmetrized offset vectors of $(P^{\tt st}_0)^{(i)}$ and $(Q^{\tt st}_0)^{(i)}$ respectively. It is immediate from Lemma \ref{lem:std and msstd} that $\eta^{(i)} = \eta^{(i)}_{\tt st}$ and $\theta^{(i)} = \theta^{(i)}_{\tt st}$. Hence the assertion follows.
\qed

\begin{rem}\label{rem:dominance CPY}
{\rm
By Remark \ref{rem:offser constant CPY}, we see that $(P^{\tt st}_0, Q^{\tt st}_0, \rho)$ is dominant if and only if $\rho$ is dominant in the sense of \cite[Definition 5.8]{CPY}.
} 
\end{rem}

Now, we can describe the image of $\kappa_0$ as follows.
\begin{prop}\label{prop: kappa_0}
Let
\begin{equation*}\label{eq:the image}
\Omega_{\rm dom} = \left\{\,(P_0, Q_0, \rho) \ \Bigg\vert \ 
\begin{array}{l}
(1)\ P_0 \in CSST_{[m]}(\la), \ Q_0 \in CSST_{[n]}(\la) \ \ (\la \in \cP_m \cap \cP_n), \\ 
(2)\ \rho \in \Z^{\la_1},\\
(3)\ \text{$(P_0, Q_0, \rho)$ is dominant} \ \\ 
\end{array}
\,\right\}.
\end{equation*}
Then the map $\kappa_0$ gives a bijection 
\begin{equation*}
\xymatrixcolsep{2.5pc}\xymatrixrowsep{0.5pc}
\xymatrix{
\wh{\M}_{m\times n}  \ \ar@{->}[r]  & \ \ \Omega_{\rm dom}
}.
\end{equation*}
\end{prop}
\pf
Suppose that $A \in \wh{\M}_{m\times n}$ is given. 
We have $\kappa_0(A)=(P_0, Q_0, \rho)$ for some $P_0 \in CSST_{[m]}(\la)$, $Q_0 \in CSST_{[n]}(\la)$ with $\la \in \cP_m \cap \cP_n$, and $\rho \in \Z^{\la_1}$. 

Then we have $\kappa_0(A^{\tt st})=(P^{\tt st}_0, Q^{\tt st}, \rho)$ by Proposition \ref{prop:affine RSK and standardization}. 
Since $A^{\tt st}$ is an extended affine permutation, $(P^{\tt st}_0, Q^{\tt st},\rho)$ is dominant by \cite[Theorem 5.11]{CPY} (see Remark \ref{rem:dominance CPY}), and hence $\kappa_0(A)=(P_0, Q_0, \rho)$ is dominant by Lemma \ref{lem:std preserves dominance}.

Conversely, suppose that $(P_0, Q_0, \rho) \in \Omega_{\rm dom}$ is given. We claim that there exists unique $A \in \wh{\M}_{m\times n}$ with $\kappa_0(A)=(P_0, Q_0, \rho)$. 

Let $\alpha$ and $\beta$ be the contents of $P_0$ and $Q_0$ respectively, and let $K= |\alpha| = |\beta| = |\la|$. 
By Lemma \ref{lem:std preserves dominance}, $(P^{\tt st}_0, Q^{\tt st}_0, \rho)$ is dominant. 
Then it follows from \cite[Theorem 5.1, Theorem 5.12]{CPY} that there exists a unique extended affine permutation $B \in \wh{\M}_{K \times K}$ such that $\kappa_0(B)=(P^{\tt st}_0, Q^{\tt st}_0, \rho)$. 
By Corollary \ref{cor:affine RS preserves descents}, $B$ is $(\alpha, \beta)$-descending, and hence by Lemma \ref{lem:inverse std} there exists a unique $A \in \wh{\M}_{m \times n}(\alpha,\beta)$ such that $A^{\tt st}=B$. 
Finally, we have $\kappa_0(A)=(P_0, Q_0, \rho)$ by Lemma \ref{lem:std of tableaux}.
\qed

\begin{cor}
The map $\kappa$ is well-defined.
\end{cor}
\pf It follows from Lemma \ref{cor:offset constant}.
\qed

\begin{lem}\label{lem:kapp_0 to kappa}
For $\la \in \cP_m\cap\cP_n$, let $\Omega_{\rm dom}(\la)$ denote the set of $(P_0, Q_0, \rho) \in \Omega_{\rm dom}$ such that the shape of $P_0$ and $Q_0$ is $\la$. 
Then we have a bijection 
\begin{equation*}
\xymatrixcolsep{2.5pc}\xymatrixrowsep{0.5pc}\xymatrix{
\Omega_{\rm dom}(\la)  \ \ar@{->}[r]  & \ \
 CSST_{[m]}(\la) \times \mc{B}_n(\la) \\
 (P_0, Q_0, \rho) \ \ar@{|->}[r]  & \ \ (P_0,Q)
},
\end{equation*}
where $Q$ is given as in \eqref{eq:(P,Q)}.
\end{lem}

\pf Consider the bijection
\begin{equation}\label{eq:Ohm dom}
\xymatrixcolsep{2.5pc}\xymatrixrowsep{0.5pc}\xymatrix{
\Omega_{\rm dom}(\la)  \ \ar@{->}[r]  & \ \
 CSST_{[m]}(\la) \times  CSST_{[n]}(\la) \times \mc{P}(\la) \\
 (P_0, Q_0, \rho) \ \ar@{|->}[r]  & \ \ (P_0, Q_0, \rho_{\tt rev} - \zeta)
}.
\end{equation} 
Then the map $(P_0,Q_0,\rho)\longmapsto (P_0,Q)$ is a composition of \eqref{eq:Ohm dom} followed by the bijection \eqref{eq:CSST and B general} to the last two components. Hence the bijectivity follows.
\qed
\newline

\noindent Since $\Omega_{\rm dom}=\bigsqcup_{\la\in\cP_m\cap\cP_n}\Omega_{\rm dom}(\la)$, Theorem \ref{thm:main-1} follows from Proposition \ref{prop: kappa_0} and Lemma \ref{lem:kapp_0 to kappa}. This completes the proof.

\begin{rem}{\rm
The well-definedness of $\kappa$ also follows directly from Proposition \ref{prop: kappa_0} and Lemma \ref{lem:kapp_0 to kappa}.}
\end{rem}

\section{Proof of Theorem \ref{thm:main-2}}\label{sec:proof of main-2}

In this section, we prove 
\begin{equation}\label{eq:kappa commutes with e and f}
 \kappa(\td{x}_iA)=\td{x}_i\kappa(A),\quad  \kappa(\td{y}^t_jA)=\td{y}^t_j\kappa(A),
\end{equation}
for $A\in \wh{\mc M}_{m\times n}$, $i\in \{\,1,\dots,m-1\,\}, j\in \{\,0,1,\dots,n-1\,\}$ and $x,y\in \{\,e,f\,\}$. 

\subsection{Notations}
Let ${\bf s}$ be a stream of flow $l$ with defining data $({\bf a}, {\bf b}, r)$. We regard ${\bf s}$ as an element of $\mc{T}_{m\times n}$ in \eqref{eq:T mn}, which is both an $A^{(1)}_{m-1}$ and $A^{(1)}_{n-1}$ crystal with $P^0_{\rm cl}$-weights, {given by}
\begin{equation*}
 {\bf s}=({\bf a}, \tau^{r}{\bf b}) \in CSST_{[m]}((1^l)) \times \mc{B}_n((1^l)) \subset \mc{T}_{m\times n}.
\end{equation*}

For $A\in \wh{\mc M}_{m\times n}$ with the corresponding streams ${\bf s}^{(1)},\dots,{\bf s}^{(s)}$ in Section \ref{subsec:mbc}, we identify $\kappa(A)$ with ${\bf s}^{(s)} \otimes \cdots \otimes {\bf s}^{(1)}$.

Define a map  
\begin{equation*}
 \xymatrixcolsep{2.5pc}\xymatrixrowsep{0.5pc}\xymatrix{
\Psi : \  \wh{\mc M}_{m\times n}  \ \ar@{->}[r]  & \ \
\displaystyle{\bigsqcup_{l\ge 0}\ \wh{\mc M}_{m\times n}\otimes \left( CSST_{[m]}((1^l)) \times \mc{B}_n((1^l))\right)} \\
 A \ \ar@{|->}[r]  & \ \ A^\flat \otimes {\bf s}^{(1)}
}.
\end{equation*}
Since
\begin{equation}\label{eq:Psi composition}
 \left((\Psi \ot  {\rm id}^{\ot s-1}) \circ \cdots  \circ (\Psi \ot {\rm id}) \circ \Psi\right)(A) = {\mathbb O} \ot {\bf s}^{(s)} \otimes \cdots \otimes {\bf s}^{(1)}={\mathbb O} \ot \kappa(A),
\end{equation}
where ${\rm id}$ is the identity morphism, it suffices to show that $\Psi$ commutes with $\td{x}_i$ and $\td{y}^t_j$ for the proof of \eqref{eq:kappa commutes with e and f}.

In order to simplify the description of $\td{x}_i$ and $\td{y}^t_j$ on $A^{\flat} \ot {\bf s}$ (see \eqref{eq:f on A*}), let us introduce some additional notations and conventions.
Let $\Z^\ast = \Z\cup\{\infty\}$, where we understand that $a<\infty$ and $a + \infty = \infty$ for $a \in \Z$.
{Let $A\in \wh{\mc M}_{m\times n}$ be given. We define $A^\ast=({a}^*_{ij})_{i, j \in {\Z}^*}$ by}
\begin{equation*}
{a}^*_{ij} = 
\begin{cases}
a^{\flat}_{ij} & \text{if $(i, j) \in \Z \times \Z$},\\
1 & \text{if $i=\infty$ and $(k,j) \in {\bf s}^{(1)}$ for some $k\in \Z$},\\
1 & \text{if $j=\infty$ and $(i,k) \in {\bf s}^{(1)}$ for some $k\in \Z$},\\
0 & \text{otherwise},
\end{cases}
\end{equation*} where $A^{\flat}=(a^{\flat}_{ij})_{i,j \in \Z}$.
In other words, ${A}^*$ is an augmented matrix obtained from $A^{\flat}$ by
$${A}^* = A^{\flat} + \sum_{(i,j) \in {\bf s}^{(1)}}(E_{i\infty} + E_{\infty j}).$$ 
Note that ${A}^*$ satisfies ${a}^*_{i+m\, j+n} = {a}^*_{ij}$ for $(i,j)\in \Z^*\times\Z^*$.
%

Let ${\bf z}$ be a zig-zag of $A$ with the back-post corner $(i, j)$. 
Let ${\bf z}^* = {\bf z} \cup \{(\infty, j), (i, \infty) \}$ and regard $(\infty, j)$ and $(i, \infty)$ as outer corners of ${\bf z}^*$. 
Then we may understand $\Psi(A) = {A}^*$ as a $\Z^*\times\Z^*$-matrix obtained by 
\begin{itemize}
 \item[$\bullet$] identifying $A=(a_{ij})_{i,j\in\Z}$ with $(a_{ij})_{i, j \in {\Z}^*}$ where $a_{\infty j} = a_{i \infty} = 0$ for $i, j \in {\Z}^*$,
 
 \item[$\bullet$] applying the same rules for $\flat$ in Section \ref{subsec:mbc} to $A$ along ${\bf z}^*_k$ instead of ${\bf z}_k$,
\end{itemize} 
where $\{{\bf z}_k\}_{k \in \Z}$ is the set of zig-zags associated to $d^{\tt sw}_A$. Note that we can recover $A^\flat$ and ${\bf s}^{(1)}$ from $A^*$ and $\{{\bf z}^*_k\}_{k \in \Z}$.
From now on, we assume that a matrix is a $\Z^*\times\Z^*$-matrix and a zig-zag is of the form ${\bf z}^*$.

\begin{ex}\label{ex:A*}{\rm
Let $A$ be the generalized affine permutation in Example \ref{ex:A}. We regard $A$ as a $\Z^* \times \Z^*$ matrix as follows:\vskip 2mm

\begin{center}
\begin{tikzpicture}[every node/.style={font=\footnotesize, scale=1}]
\matrix (M)[matrix of math nodes,nodes in empty cells,nodes={rectangle,minimum height=1.0em,minimum width=1.0em,inner sep=0pt,anchor=center,align=center}]
{
& &&&&& &&&&& &&&&& &&&&& &&&&& &&&&& &&&\\
& &1&&&& &&&1&& &&&&& &&&&& &&&&& &&&&& &&&\\
& 1&&&&& &&1&&& &&2&&& &&&&& &&&&& &&&&& &&&\\
& &1&&1&& &&&&& 1&&&1&& &&&&& &&&&& &&&&& &&&\\
& &&&&& &&&&& &1&&&1& &&&&& &&&&& &&&&& &&&\\
& &&&&& &1&&&& &&&1&& &&&&& &&&&& &&&&& &&&\\
& &&&&& 1&&&&& &&1&&& &&2&&& &&&&& &&&&& &&&\\
& &&&&& &1&&1&& &&&&& 1&&&1&& &&&&& &&&&& &&&\\
& &&&&& &&&&& &&&&& &1&&&1& &&&&& &&&&& &&&\\
& &&&&& &&&&& &1&&&& &&&1&& &&&&& &&&&& &&&\\
& &&&&& &&&&& 1&&&&& &&1&&& &&2&&& &&&&& &&&\\
& &&&&& &&&&& &1&&1&& &&&&& 1&&&1&& &&&&& &&&\\
& &&&&& &&&&& &&&&& &&&&& &1&&&1& &&&&& &&&\\
& &&&&& &&&&& &&&&& &1&&&& &&&1&& &&&&& &&&\\
& &&&&& &&&&& &&&&& 1&&&&& &&1&&& &&2&&& &&&\\
& &&&&& &&&&& &&&&& &1&&1&& &&&&& 1&&&1&& &&&\\
& &&&&& &&&&& &&&&& &&&&& &&&&& &1&&&1& &&&\\
& &&&&& &&&&& &&&&& &&&&& &&&&& &&&&& &&&\\
& &&&&& &&&&& &&&&& &&&&& &&&&& &&&&& &&&\\
& &&&&& &&&&& &&&&& &&&&& &&&&& &&&&& &&&\\
& &&&&& &&&&& &&&&& &&&&& &&&&& &&&&& &&&\\
};
\draw[solid,lightgray](M-1-1.south west)--(M-1-32.south east);
\draw[solid,lightgray](M-1-35.south west)--(M-1-35.south east);
\draw[dotted,lightgray](M-2-1.south west)--(M-2-32.south east);
\draw[dotted,lightgray](M-2-35.south west)--(M-2-35.south east);
\draw[dotted,lightgray](M-3-1.south west)--(M-3-32.south east);
\draw[dotted,lightgray](M-3-35.south west)--(M-3-35.south east);
\draw[dotted,lightgray](M-4-1.south west)--(M-4-32.south east);
\draw[dotted,lightgray](M-4-35.south west)--(M-4-35.south east);
\draw[-latex,thick,gray](M-5-1.south west)--(M-5-32.south east);
\node (jinf) at (M-1-35.north) {$j=\infty$};
\draw[thick,gray](M-5-35.south west)--(M-5-35.south east);
\draw[dotted,lightgray](M-6-1.south west)--(M-6-32.south east);
\draw[dotted,lightgray](M-6-35.south west)--(M-6-35.south east);
\draw[dotted,lightgray](M-7-1.south west)--(M-7-32.south east);
\draw[dotted,lightgray](M-7-35.south west)--(M-7-35.south east);
\draw[dotted,lightgray](M-8-1.south west)--(M-8-32.south east);
\draw[dotted,lightgray](M-8-35.south west)--(M-8-35.south east);
\draw[solid,lightgray](M-9-1.south west)--(M-9-32.south east);
\draw[solid,lightgray](M-9-35.south west)--(M-9-35.south east);
\draw[dotted,lightgray](M-10-1.south west)--(M-10-32.south east);
\draw[dotted,lightgray](M-10-35.south west)--(M-10-35.south east);
\draw[dotted,lightgray](M-11-1.south west)--(M-11-32.south east);
\draw[dotted,lightgray](M-11-35.south west)--(M-11-35.south east);
\draw[dotted,lightgray](M-12-1.south west)--(M-12-32.south east);
\draw[dotted,lightgray](M-12-35.south west)--(M-12-35.south east);
\draw[solid,lightgray](M-13-1.south west)--(M-13-32.south east);
\draw[solid,lightgray](M-13-35.south west)--(M-13-35.south east);
\draw[dotted,lightgray](M-14-1.south west)--(M-14-32.south east);
\draw[dotted,lightgray](M-14-35.south west)--(M-14-35.south east);
\draw[dotted,lightgray](M-15-1.south west)--(M-15-32.south east);
\draw[dotted,lightgray](M-15-35.south west)--(M-15-35.south east);
\draw[dotted,lightgray](M-16-1.south west)--(M-16-32.south east);
\draw[dotted,lightgray](M-16-35.south west)--(M-16-35.south east);
\draw[solid,lightgray](M-17-1.south west)--(M-17-32.south east);
\draw[solid,lightgray](M-17-35.south west)--(M-17-35.south east);
\draw[solid,lightgray](M-1-1.north east)--(M-18-1.south east);
\draw[solid,lightgray](M-21-1.north east)--(M-21-1.south east);
\draw[dotted,lightgray](M-1-2.north east)--(M-18-2.south east);
\draw[dotted,lightgray](M-21-2.north east)--(M-21-2.south east);
\draw[dotted,lightgray](M-1-3.north east)--(M-18-3.south east);
\draw[dotted,lightgray](M-21-3.north east)--(M-21-3.south east);
\draw[dotted,lightgray](M-1-4.north east)--(M-18-4.south east);
\draw[dotted,lightgray](M-21-4.north east)--(M-21-4.south east);
\draw[dotted,lightgray](M-1-5.north east)--(M-18-5.south east);
\draw[dotted,lightgray](M-21-5.north east)--(M-21-5.south east);
\draw[solid,lightgray](M-1-6.north east)--(M-18-6.south east);
\draw[solid,lightgray](M-21-6.north east)--(M-21-6.south east);
\draw[dotted,lightgray](M-1-7.north east)--(M-18-7.south east);
\draw[dotted,lightgray](M-21-7.north east)--(M-21-7.south east);
\draw[dotted,lightgray](M-1-8.north east)--(M-18-8.south east);
\draw[dotted,lightgray](M-21-8.north east)--(M-21-8.south east);
\draw[dotted,lightgray](M-1-9.north east)--(M-18-9.south east);
\draw[dotted,lightgray](M-21-9.north east)--(M-21-9.south east);
\draw[dotted,lightgray](M-1-10.north east)--(M-18-10.south east);
\draw[dotted,lightgray](M-21-10.north east)--(M-21-10.south east);
\draw[-latex,thick,gray](M-1-11.north east)--(M-18-11.south east);
\node [anchor=east] (iinf) at (M-21-1.west) {$i=\infty$};
\draw[thick,gray](M-21-11.north east)--(M-21-11.south east);
\draw[dotted,lightgray](M-1-12.north east)--(M-18-12.south east);
\draw[dotted,lightgray](M-21-12.north east)--(M-21-12.south east);
\draw[dotted,lightgray](M-1-13.north east)--(M-18-13.south east);
\draw[dotted,lightgray](M-21-13.north east)--(M-21-13.south east);
\draw[dotted,lightgray](M-1-14.north east)--(M-18-14.south east);
\draw[dotted,lightgray](M-21-14.north east)--(M-21-14.south east);
\draw[dotted,lightgray](M-1-15.north east)--(M-18-15.south east);
\draw[dotted,lightgray](M-21-15.north east)--(M-21-15.south east);
\draw[solid,lightgray](M-1-16.north east)--(M-18-16.south east);
\draw[solid,lightgray](M-21-16.north east)--(M-21-16.south east);
\draw[dotted,lightgray](M-1-17.north east)--(M-18-17.south east);
\draw[dotted,lightgray](M-21-17.north east)--(M-21-17.south east);
\draw[dotted,lightgray](M-1-18.north east)--(M-18-18.south east);
\draw[dotted,lightgray](M-21-18.north east)--(M-21-18.south east);
\draw[dotted,lightgray](M-1-19.north east)--(M-18-19.south east);
\draw[dotted,lightgray](M-21-19.north east)--(M-21-19.south east);
\draw[dotted,lightgray](M-1-20.north east)--(M-18-20.south east);
\draw[dotted,lightgray](M-21-20.north east)--(M-21-20.south east);
\draw[solid,lightgray](M-1-21.north east)--(M-18-21.south east);
\draw[solid,lightgray](M-21-21.north east)--(M-21-21.south east);
\draw[dotted,lightgray](M-1-22.north east)--(M-18-22.south east);
\draw[dotted,lightgray](M-21-22.north east)--(M-21-22.south east);
\draw[dotted,lightgray](M-1-23.north east)--(M-18-23.south east);
\draw[dotted,lightgray](M-21-23.north east)--(M-21-23.south east);
\draw[dotted,lightgray](M-1-24.north east)--(M-18-24.south east);
\draw[dotted,lightgray](M-21-24.north east)--(M-21-24.south east);
\draw[dotted,lightgray](M-1-25.north east)--(M-18-25.south east);
\draw[dotted,lightgray](M-21-25.north east)--(M-21-25.south east);
\draw[solid,lightgray](M-1-26.north east)--(M-18-26.south east);
\draw[solid,lightgray](M-21-26.north east)--(M-21-26.south east);
\draw[dotted,lightgray](M-1-27.north east)--(M-18-27.south east);
\draw[dotted,lightgray](M-21-27.north east)--(M-21-27.south east);
\draw[dotted,lightgray](M-1-28.north east)--(M-18-28.south east);
\draw[dotted,lightgray](M-21-28.north east)--(M-21-28.south east);
\draw[dotted,lightgray](M-1-29.north east)--(M-18-29.south east);
\draw[dotted,lightgray](M-21-29.north east)--(M-21-29.south east);
\draw[dotted,lightgray](M-1-30.north east)--(M-18-30.south east);
\draw[dotted,lightgray](M-21-30.north east)--(M-21-30.south east);
\draw[solid,lightgray](M-1-31.north east)--(M-18-31.south east);
\draw[solid,lightgray](M-21-31.north east)--(M-21-31.south east);
\begin{scope}[on background layer]
\draw[thin,red](M-21-15.center) -- (M-12-15.center) -| (M-8-17.center) -| (M-7-19.center)  -- (M-7-35.center);
\draw[thin,red](M-21-17.center) -- (M-15-17.center) -| (M-9-18.center) -| (M-8-20.center)  -- (M-8-35.center);
\draw[thin,red](M-21-18.center) -- (M-16-18.center) -| (M-11-19.center) -| (M-10-20.center) -| (M-9-21.center)  -- (M-9-35.center);
\end{scope}
\node [anchor=east] at (M.west) {$A=$};
\end{tikzpicture}
\end{center}
where the red lines denote the zig-zags ${\bf z}_1, {\bf z}_2, {\bf z}_3$ associated to $d^{\tt sw}_A$. Then $A^*$ is given as follows.\vskip 2mm

\begin{center}
\begin{tikzpicture}[every node/.style={font=\footnotesize,scale=1}]
\matrix (M)[matrix of math nodes,nodes in empty cells,nodes={rectangle,minimum height=1.0em,minimum width=1.0em,inner sep=0pt,anchor=center,align=center}]
{
& &&&&& &&&&& &&&&& &&&&& &&&&& &&&&& &&&\\
& &1&&&& &&&&1& &&&&& &&&&& &&&&& &&&&& &&&\\
& &1&&&& &&&1&& &&1&&& &&&&& &&&&& &&&&& &&&1\\
& &&1&&& 1&&&&& &&1&&& &&&&& &&&&& &&&&& &&&1\\
& &&&&& &&&&& &&&1&& &&&&& &&&&& &&&&& &&&1\\
& &&&&& &1&&&& &&&&1& &&&&& &&&&& &&&&& &&&\\
& &&&&& &1&&&& &&&1&& &&1&&& &&&&& &&&&& &&&1\\
& &&&&& &&1&&& 1&&&&& &&1&&& &&&&& &&&&& &&&1\\
& &&&&& &&&&& &&&&& &&&1&& &&&&& &&&&& &&&1\\
& &&&&& &&&&& &1&&&& &&&&1& &&&&& &&&&& &&&\\
& &&&&& &&&&& &1&&&& &&&1&& &&1&&& &&&&& &&&1\\
& &&&&& &&&&& &&1&&& 1&&&&& &&1&&& &&&&& &&&1\\
& &&&&& &&&&& &&&&& &&&&& &&&1&& &&&&& &&&1\\
& &&&&& &&&&& &&&&& &1&&&& &&&&1& &&&&& &&&\\
& &&&&& &&&&& &&&&& &1&&&& &&&1&& &&1&&& &&&1\\
& &&&&& &&&&& &&&&& &&1&&& 1&&&&& &&1&&& &&&1\\
& &&&&& &&&&& &&&&& &&&&& &&&&& &&&1&& &&&1\\
& &&&&& &&&&& &&&&& &&&&& &&&&& &&&&& &&&\\
& &&&&& &&&&& &&&&& &&&&& &&&&& &&&&& &&&\\
& &&&&& &&&&& &&&&& &&&&& &&&&& &&&&& &&&\\
& 1&1&&1&& 1&1&&1&& 1&1&&1&& 1&1&&1&& 1&1&&1&& 1&1&&1&& &&&\\
};
\draw[solid,lightgray](M-1-1.south west)--(M-1-32.south east);
\draw[solid,lightgray](M-1-35.south west)--(M-1-35.south east);
\draw[dotted,lightgray](M-2-1.south west)--(M-2-32.south east);
\draw[dotted,lightgray](M-2-35.south west)--(M-2-35.south east);
\draw[dotted,lightgray](M-3-1.south west)--(M-3-32.south east);
\draw[dotted,lightgray](M-3-35.south west)--(M-3-35.south east);
\draw[dotted,lightgray](M-4-1.south west)--(M-4-32.south east);
\draw[dotted,lightgray](M-4-35.south west)--(M-4-35.south east);
\draw[-latex,thick,gray](M-5-1.south west)--(M-5-32.south east);
\node (jinf) at (M-1-35.north) {$j=\infty$};
\draw[thick,gray](M-5-35.south west)--(M-5-35.south east);
\draw[dotted,lightgray](M-6-1.south west)--(M-6-32.south east);
\draw[dotted,lightgray](M-6-35.south west)--(M-6-35.south east);
\draw[dotted,lightgray](M-7-1.south west)--(M-7-32.south east);
\draw[dotted,lightgray](M-7-35.south west)--(M-7-35.south east);
\draw[dotted,lightgray](M-8-1.south west)--(M-8-32.south east);
\draw[dotted,lightgray](M-8-35.south west)--(M-8-35.south east);
\draw[solid,lightgray](M-9-1.south west)--(M-9-32.south east);
\draw[solid,lightgray](M-9-35.south west)--(M-9-35.south east);
\draw[dotted,lightgray](M-10-1.south west)--(M-10-32.south east);
\draw[dotted,lightgray](M-10-35.south west)--(M-10-35.south east);
\draw[dotted,lightgray](M-11-1.south west)--(M-11-32.south east);
\draw[dotted,lightgray](M-11-35.south west)--(M-11-35.south east);
\draw[dotted,lightgray](M-12-1.south west)--(M-12-32.south east);
\draw[dotted,lightgray](M-12-35.south west)--(M-12-35.south east);
\draw[solid,lightgray](M-13-1.south west)--(M-13-32.south east);
\draw[solid,lightgray](M-13-35.south west)--(M-13-35.south east);
\draw[dotted,lightgray](M-14-1.south west)--(M-14-32.south east);
\draw[dotted,lightgray](M-14-35.south west)--(M-14-35.south east);
\draw[dotted,lightgray](M-15-1.south west)--(M-15-32.south east);
\draw[dotted,lightgray](M-15-35.south west)--(M-15-35.south east);
\draw[dotted,lightgray](M-16-1.south west)--(M-16-32.south east);
\draw[dotted,lightgray](M-16-35.south west)--(M-16-35.south east);
\draw[solid,lightgray](M-17-1.south west)--(M-17-32.south east);
\draw[solid,lightgray](M-17-35.south west)--(M-17-35.south east);
\draw[solid,lightgray](M-1-1.north east)--(M-18-1.south east);
\draw[solid,lightgray](M-21-1.north east)--(M-21-1.south east);
\draw[dotted,lightgray](M-1-2.north east)--(M-18-2.south east);
\draw[dotted,lightgray](M-21-2.north east)--(M-21-2.south east);
\draw[dotted,lightgray](M-1-3.north east)--(M-18-3.south east);
\draw[dotted,lightgray](M-21-3.north east)--(M-21-3.south east);
\draw[dotted,lightgray](M-1-4.north east)--(M-18-4.south east);
\draw[dotted,lightgray](M-21-4.north east)--(M-21-4.south east);
\draw[dotted,lightgray](M-1-5.north east)--(M-18-5.south east);
\draw[dotted,lightgray](M-21-5.north east)--(M-21-5.south east);
\draw[solid,lightgray](M-1-6.north east)--(M-18-6.south east);
\draw[solid,lightgray](M-21-6.north east)--(M-21-6.south east);
\draw[dotted,lightgray](M-1-7.north east)--(M-18-7.south east);
\draw[dotted,lightgray](M-21-7.north east)--(M-21-7.south east);
\draw[dotted,lightgray](M-1-8.north east)--(M-18-8.south east);
\draw[dotted,lightgray](M-21-8.north east)--(M-21-8.south east);
\draw[dotted,lightgray](M-1-9.north east)--(M-18-9.south east);
\draw[dotted,lightgray](M-21-9.north east)--(M-21-9.south east);
\draw[dotted,lightgray](M-1-10.north east)--(M-18-10.south east);
\draw[dotted,lightgray](M-21-10.north east)--(M-21-10.south east);
\draw[-latex,thick,gray](M-1-11.north east)--(M-18-11.south east);
\node [anchor=east] (iinf) at (M-21-1.west) {$i=\infty$};
\draw[thick,gray](M-21-11.north east)--(M-21-11.south east);
\draw[dotted,lightgray](M-1-12.north east)--(M-18-12.south east);
\draw[dotted,lightgray](M-21-12.north east)--(M-21-12.south east);
\draw[dotted,lightgray](M-1-13.north east)--(M-18-13.south east);
\draw[dotted,lightgray](M-21-13.north east)--(M-21-13.south east);
\draw[dotted,lightgray](M-1-14.north east)--(M-18-14.south east);
\draw[dotted,lightgray](M-21-14.north east)--(M-21-14.south east);
\draw[dotted,lightgray](M-1-15.north east)--(M-18-15.south east);
\draw[dotted,lightgray](M-21-15.north east)--(M-21-15.south east);
\draw[solid,lightgray](M-1-16.north east)--(M-18-16.south east);
\draw[solid,lightgray](M-21-16.north east)--(M-21-16.south east);
\draw[dotted,lightgray](M-1-17.north east)--(M-18-17.south east);
\draw[dotted,lightgray](M-21-17.north east)--(M-21-17.south east);
\draw[dotted,lightgray](M-1-18.north east)--(M-18-18.south east);
\draw[dotted,lightgray](M-21-18.north east)--(M-21-18.south east);
\draw[dotted,lightgray](M-1-19.north east)--(M-18-19.south east);
\draw[dotted,lightgray](M-21-19.north east)--(M-21-19.south east);
\draw[dotted,lightgray](M-1-20.north east)--(M-18-20.south east);
\draw[dotted,lightgray](M-21-20.north east)--(M-21-20.south east);
\draw[solid,lightgray](M-1-21.north east)--(M-18-21.south east);
\draw[solid,lightgray](M-21-21.north east)--(M-21-21.south east);
\draw[dotted,lightgray](M-1-22.north east)--(M-18-22.south east);
\draw[dotted,lightgray](M-21-22.north east)--(M-21-22.south east);
\draw[dotted,lightgray](M-1-23.north east)--(M-18-23.south east);
\draw[dotted,lightgray](M-21-23.north east)--(M-21-23.south east);
\draw[dotted,lightgray](M-1-24.north east)--(M-18-24.south east);
\draw[dotted,lightgray](M-21-24.north east)--(M-21-24.south east);
\draw[dotted,lightgray](M-1-25.north east)--(M-18-25.south east);
\draw[dotted,lightgray](M-21-25.north east)--(M-21-25.south east);
\draw[solid,lightgray](M-1-26.north east)--(M-18-26.south east);
\draw[solid,lightgray](M-21-26.north east)--(M-21-26.south east);
\draw[dotted,lightgray](M-1-27.north east)--(M-18-27.south east);
\draw[dotted,lightgray](M-21-27.north east)--(M-21-27.south east);
\draw[dotted,lightgray](M-1-28.north east)--(M-18-28.south east);
\draw[dotted,lightgray](M-21-28.north east)--(M-21-28.south east);
\draw[dotted,lightgray](M-1-29.north east)--(M-18-29.south east);
\draw[dotted,lightgray](M-21-29.north east)--(M-21-29.south east);
\draw[dotted,lightgray](M-1-30.north east)--(M-18-30.south east);
\draw[dotted,lightgray](M-21-30.north east)--(M-21-30.south east);
\draw[solid,lightgray](M-1-31.north east)--(M-18-31.south east);
\draw[solid,lightgray](M-21-31.north east)--(M-21-31.south east);
\begin{scope}[on background layer]
\draw[densely dashed,red](M-21-15.center) -- (M-12-15.center) -| (M-8-17.center) -| (M-7-19.center)  -- (M-7-35.center);
\draw[densely dashed,red](M-21-17.center) -- (M-15-17.center) -| (M-9-18.center) -| (M-8-20.center)  -- (M-8-35.center);
\draw[densely dashed,red](M-21-18.center) -- (M-16-18.center) -| (M-11-19.center) -| (M-10-20.center) -| (M-9-21.center)  -- (M-9-35.center);
\end{scope}
\node [anchor=east] at (M.west) {$A^*=$};
\end{tikzpicture}
\end{center}

}
\end{ex}

\subsection{Tensor product rule}\label{subsec:tensor product rule}
From now on, we fix $A\in \wh{\mc M}_{m\times n}$ and $j\in [n]$.
If there is no confusion, let us write $\tf_j$, $\varepsilon_j$, and $\varphi_j$ instead of $\tf_j^t$, $\varepsilon_j^t$, and $\varphi_j^t$ for simplicity.
In the {remainder} of this section, we will focus on the proof of
\begin{equation}\label{eq:main equation}
 \tf_j \Psi(A) =\Psi (\tf_jA).
\end{equation}

Let  
\begin{equation}\label{eq:sigma ast}
\begin{split}
\sigma & = ( \cdots,\ \underbrace{-\dots -}_{a_{i-1 j+1}},\ \underbrace{+\dots +}_{a_{i-1 j}},\ \underbrace{-\dots -}_{a_{i j+1}},\ \underbrace{+\dots +}_{a_{i j}},\ \cdots ),\\
{\sigma}^* & = 
( \cdots,\ \underbrace{-\dots -}_{{a}^*_{i-1 j+1}},\ \underbrace{+\dots +}_{{a}^*_{i-1 j}},\ \underbrace{-\dots -}_{{a}^*_{i j+1}},\ \underbrace{+\dots +}_{{a}^*_{i j}},\ \cdots) \cdot (\underbrace{+}_{{a}^*_{\infty, j}},\ \underbrace{-}_{{a}^*_{\infty j+1}}), 
\end{split}
\end{equation}
where $\sigma^*$ is a concatenation of two sequences. 
By tensor product rule \eqref{eq:tensor_product_rule}, we see that
\begin{equation}\label{eq:Psi and f_j}
\td{f}_j\left(A^{\flat}\otimes{\bf s}\right) = 
\begin{cases}
A^{\flat} \otimes \td{f}_j{\bf s} & \text{if the leftmost $+$ in $\td{{\sigma}^*}$ corresponds to $(\infty, j)$}, \\
\left(\td{f}_j A^{\flat}\right) \otimes {\bf s} & \text{if the leftmost $+$ in $\td{{\sigma}^*}$ corresponds to $(i, j)$ for some $i < \infty $}, \\
{\bf 0} & \text{if $\td{{\sigma}^*}$ has no $+$}.
\end{cases}
\end{equation}
In terms of $A^*$, this can be simplified as
\begin{equation}\label{eq:f on A*}
\tf_j {A}^* = 
\begin{cases}
{A}^* - \wh{E}_{ij} + \wh{E}_{i\, j+1} & \text{if the leftmost $+$ in $\td{{\sigma}^*}$ corresponds to $(i, j)$}, \\
{\bf 0} & \text{if $\td{{\sigma}^*}$ has no $+$},
\end{cases}
\end{equation}
where we {define} $\wh{E}_{\infty j} = \sum_{k \in \Z} E_{\infty\, j+kn}$.
 
\begin{lem}\label{lem:Psi preserves signature}
We have $\td{f}_jA \neq {\bf 0}$ if and only if $\td{f}_j\Psi(A) \neq {\bf 0}$.
\end{lem}
\pf We may assume that there exists a non-zero cell in the $j$-th column. Otherwise, we have $\td{f}_jA = \td{f}_j\Psi(A) = {\bf 0}$. 

Let $\{{\bf z}_k\}_{k \in \Z}$ be the set of zig-zags associated to $d_A^{\tt sw}$. Let $k_0$ (resp. $k_1$) be the minimal (resp. maximal) value of $d_A^{\tt sw}$ in the $j$-th column. 
For $k_0 \le k \le k_1$, let
$i_k$ be the minimal row index with $(i_k, j) \in {\bf z}_k$. Note that $(i_k, j + 1) \in {\bf z}_k$. 

Put
\begin{equation*}
\begin{split}
\sigma_k & = (\underbrace{+\dots +}_{a_{i_{k}\, j}},\ \underbrace{-\dots -}_{a_{i_{k}+1\, j+1}},\ \cdots,\ \underbrace{+\dots +}_{a_{i_{k+1}-1\, j}},\ \underbrace{-\dots -}_{a_{i_{k+1}\, j+1}}),\\
{\sigma}^*_k & = (\underbrace{+\dots +}_{{a}^*_{i_{k}\, j}},\ \underbrace{-\dots -}_{{a}^*_{i_{k}+1\, j+1}},\ \cdots,\ \underbrace{+\dots +}_{{a}^*_{i_{k+1}-1 \, j}},\ \underbrace{-\dots -}_{{a}^*_{i_{k+1}\, j+1}}),\\
\end{split}
\end{equation*}
for $k_0 \le k < k_1$, and 
\begin{align*}
\sigma_{-\infty} & = (\cdots,\ \underbrace{-\dots -}_{a_{i_{k_{0}}-1\, j+1}},\ \underbrace{-\dots -}_{a_{i_{k_{0}}\, j+1}}), & \sigma_{\infty} & = (\underbrace{+\dots +}_{a_{i_{k_{1}}\, j}},\ \underbrace{-\dots -}_{a_{i_{k_{1}}+1\, j+1}},\ \cdots), \\
{\sigma}^*_{-\infty} & = (\cdots,\ \underbrace{-\dots -}_{{a}^*_{i_{k_{0}}-1\, j+1}},\ \underbrace{-\dots -}_{{a}^*_{i_{k_{0}}\, j+1}}), & {\sigma}^*_{\infty} & = (\underbrace{+\dots +}_{{a}^*_{i_{k_{1}}\, j}},\ \underbrace{-\dots -}_{{a}^*_{i_{k_{1}}+1\, j+1}},\ \cdots) \cdot (\underbrace{+}_{{{a}^*_{\infty\, j}}},\ \underbrace{-}_{{{a}^*_{\infty\, j+1}}}).
\end{align*}
Then $\sigma$ and ${\sigma}^*$ in \eqref{eq:sigma ast} decompose as follows:
\begin{equation*}\label{eq:decomposition of signature}
\begin{split}
\sigma & = \sigma_{-\infty} \cdot \sigma_{k_{0}} \cdot \dots \cdot \sigma_{k_{1}-1} \cdot \sigma_{\infty},\\
{\sigma}^* & = {\sigma}^*_{-\infty} \cdot {\sigma}^*_{k_{0}} \cdot \dots \cdot {\sigma}^*_{k_{1}-1} \cdot {\sigma}^*_{\infty}.
\end{split}
\end{equation*}

Suppose that $k$ is given with $k_{0} \le k< k_{1}$. Let $u$ be the maximal row index with $(u, j) \in {\bf z}_k$, and let $v$ be the minimal row index with $(v, j+1) \in {\bf z}_{k+1}$. Note that $i_k \le u$ and $v \le i_{k+1}$. 
Suppose first that $u < v$. Then we have
\begin{equation*}
\begin{split}
\sigma_k & = (\underbrace{+\dots +}_{a_{i_{k}\, j}},\ +\dots+,\  \underbrace{+\dots +}_{a_{u\, j}},\ \underbrace{-\dots -}_{a_{v\, j+1}},\ -\dots-,\ \underbrace{-\dots -}_{a_{i_{k+1}\, j+1}}),\\
{\sigma}^*_k & = (\underbrace{+\dots +}_{a_{i_{k}\, j} - 1},\ +\dots+,\ \underbrace{+\dots +}_{a_{u\, j} + 1},\ \underbrace{-\dots -}_{a_{v\, j+1} - 1},\ -\dots-,\ \underbrace{-\dots -}_{a_{i_{k+1}\, j+1} + 1}) = \sigma_k,
\end{split}
\end{equation*}
and hence $\td{\sigma_k} = \td{{\sigma}^*_k}$.
Next, suppose that $u \ge v$. 
In this case, we have $i_k < v \le u < i_{k+1}$. 
Therefore, $(i_k, j)$, $(v, j+1)$ are inner corners, and $(u, j)$, $(i_{k+1}, j+1)$ are outer corners. Then we have
\begin{equation*}
\begin{split}
\sigma_k & = (\underbrace{+\dots +}_{a_{i_{k}\, j}},\ +\dots+,\ \underbrace{-\dots -}_{a_{v\, j+1}},\ \cdots,\  \underbrace{+\dots +}_{a_{u\, j}},\ -\dots-,\ \underbrace{-\dots -}_{a_{i_{k+1}\, j+1}}),\\
{\sigma}^*_k & = (\underbrace{+\dots +}_{a_{i_{k}\, j} - 1},\ +\dots+,\ \underbrace{-\dots -}_{a_{v\,j+1} - 1},\ \cdots,\ \underbrace{+\dots +}_{a_{u\, j} + 1},\ -\dots-,\ \underbrace{-\dots -}_{a_{i_{k+1}\, j+1} + 1}).
\end{split}
\end{equation*} 
Since one cancelling pair $(+,-)$ of $\sigma_k$ in 
$$(\underbrace{+\dots +}_{a_{i_{k}\, j}},\ +\dots+,\ \underbrace{-\dots -}_{a_{v\, j+1}})$$ is moved to a pair $(+,-)$ of ${\sigma}^*_k$ in
$$(\underbrace{+\dots +}_{a_{u\, j} + 1},\ -\dots-,\ \underbrace{-\dots -}_{a_{i_{k+1}\, j+1} + 1}),$$ 
we conclude that $\td{\sigma_k} = \td{{\sigma}^*_k}$. 
By similar argument, we see that $(\sigma_{\pm \infty})^\sim = ({\sigma}^*_{\pm \infty})^\sim$, where $(\,\cdot\,)^\sim$ means $\td{(\,\cdot\,)}$.
Hence we have
\begin{equation}\label{eq:non-red = red for k}
 \td{\sigma_k} = \td{{\sigma}^*_k}\quad (-\infty \le k \le \infty).
\end{equation}
Since reducing a sequence does not depend on the order of cancelling $(+,-)$, we have
\begin{equation*}
\begin{split}
 \td{\sigma} 
 &= (\sigma_{-\infty} \cdot \cdots \cdot \sigma_{\infty})^{\sim} 
 = ((\sigma_{-\infty})^\sim \cdot \cdots \cdot (\sigma_{\infty})^\sim)^{\sim} \\
 &= (({\sigma}^*_{-\infty})^\sim \cdot \dots \cdot ({\sigma}^*_{\infty})^\sim)^{\sim} = ({\sigma}^*_{-\infty} \cdot \dots \cdot {\sigma}^*_{\infty})^{\sim} = \td{{\sigma}^*}.
\end{split}
\end{equation*}
This shows that $\varphi_j(A) = \varphi_j(\Psi(A))$ and $\varepsilon_j(A) = \varepsilon_j(\Psi(A))$, and hence the lemma follows. 
\qed\newline

From now on, we assume that $\td{f}_jA \neq {\bf 0}$ and $\td{f}_j\Psi(A) \neq {\bf 0}$. 
We also {adopt} the following notations:
\begin{itemize}
 \item[$\bullet$] $u$ : the row index corresponding to the leftmost $+$ in $\td{\sigma}$, 
 
 \item[$\bullet$] $u^*$ : the row index corresponding to the leftmost $+$ in $\td{\sigma^*}$,
 
 \item[$\bullet$] $s= d_A^{\tt sw}(u,j)$,
 
 \item[$\bullet$] $\td{A} = \tf_jA = (\td{a}_{ij})_{i,j \in \Z^\ast}$. 
\end{itemize}
We have 
\begin{equation}\label{eq:A*-1}
 \td{A}=\tf_j{A} = A - \wh{E}_{u j} + \wh{E}_{u\, j+1},\quad  
 \tf_j {A}^* = A^* - \wh{E}_{u^* j} + \wh{E}_{u^*\, j+1}.
\end{equation}
Note that the leftmost $+$ in $\td{{\sigma}^*}$ also appears in $({\sigma}^*_{s})^\sim$ by \eqref{eq:non-red = red for k}. More explicitly, $u^*$ is the minimal row index such that $u \le u^*$ and ${a}^*_{u^*j} \neq 0$. 
We have
\begin{equation*}
\begin{cases}
 u < u^* & \text{if $(u, j)$ is an inner corner of ${\bf z}_{s}$  with $a_{uj}=1$},\\
 u = u^* & \text{otherwise}.
\end{cases}
\end{equation*}

\subsection{Southwest channel numberings on $A$ and $\td{A}$}\label{subsec:sw numbering of A and tdA}
In this subsection, we discuss the relation between the southwest channel numberings on $A$ and $\td{A}$.

For $(x,y)\in \Z^2$, let ${(x, y)}^{\wedge} = \{\,\tau^k(x,y)\,|\,k\in\Z\,\}$.
We have 
\begin{equation*}
\begin{split}
{\rm supp}(\td{A}) - {\rm supp}(A) & \subseteq {(u, j+1)}^\wedge, \quad
{\rm supp}(A) - {\rm supp}(\td{A})  \subseteq {(u, j)}^\wedge,
\end{split}
\end{equation*}
where the equalities hold when $a_{uj+1}=0$ and $a_{uj}=1$ respectively.

We first give an example, where a proper numbering on $\td{A}$ is induced by a proper numbering on $A$.

\begin{ex}{\rm
{Let $A$ be the generalized affine permutation given in Example \ref{ex:A}. It is easily checked that $\tf_2^tA = A - \wh{E}_{5\,2} + \wh{E}_{5\,3}$. Consider the set of zig-zags $Z = \{ {\bf z}_k \}_{k \in \Z}$ correspodning to the southwest channel numbering $d^{\tt st}_A$ on $A$. If we draw the zig-zgas ${\bf z}_1$, ${\bf z}_2$ and ${\bf z}_3$ over $\tf_2^tA$ as follows;
}
\begin{center}
\begin{tikzpicture}[every node/.style={font=\footnotesize,scale=1}]
\matrix (M)[matrix of math nodes,nodes in empty cells,nodes={rectangle,minimum height=1.0em,minimum width=1.0em,inner sep=0pt,anchor=center,align=center}]
{
& &&&&& &&&&& &&&&& &&&&& &&&&& &&&&& &\\
& &&1&&& &&&1&& &&&&& &&&&& &&&&& &&&&& &\\
& 1&&&&& &&1&&& &&2&&& &&&&& &&&&& &&&&& &\\
& &1&&1&& &&&&& 1&&&1&& &&&&& &&&&& &&&&& &\\
& &&&&& &&&&& &1&&&1& &&&&& &&&&& &&&&& &\\
& &&&&& &&1&&& &&&1&& &&&&& &&&&& &&&&& &\\
& &&&&& 1&&&&& &&1&&& &&2&&& &&&&& &&&&& &\\
& &&&&& &1&&1&& &&&&& 1&&&1&& &&&&& &&&&& &\\
& &&&&& &&&&& &&&&& &1&&&1& &&&&& &&&&& &\\
& &&&&& &&&&& &&1&&& &&&1&& &&&&& &&&&& &\\
& &&&&& &&&&& 1&&&&& &&1&&& &&2&&& &&&&& &\\
& &&&&& &&&&& &1&&1&& &&&&& 1&&&1&& &&&&& &\\
& &&&&& &&&&& &&&&& &&&&& &1&&&1& &&&&& &\\
& &&&&& &&&&& &&&&& &&1&&& &&&1&& &&&&& &\\
& &&&&& &&&&& &&&&& 1&&&&& &&1&&& &&2&&& &\\
& &&&&& &&&&& &&&&& &1&&1&& &&&&& 1&&&1&& &\\
& &&&&& &&&&& &&&&& &&&&& &&&&& &1&&&1& &\\
& &&&&& &&&&& &&&&& &&&&& &&&&& &&&&& &\\
& &&&&& &&&&& &&&&& &&&&& &&&&& &&&&& &\\
};
\draw[solid,lightgray](M-1-1.south west)--(M-1-32.south east);
\draw[dotted,lightgray](M-2-1.south west)--(M-2-32.south east);
\draw[dotted,lightgray](M-3-1.south west)--(M-3-32.south east);
\draw[dotted,lightgray](M-4-1.south west)--(M-4-32.south east);
\draw[-latex,thick,gray](M-5-1.south west)--(M-5-32.south east);
\node (j) at (M-6-33) {$j$};
\node (j) at (M-7-33) {${\bf z}_1$};
\node (j) at (M-8-33) {${\bf z}_2$};
\node (j) at (M-9-33) {${\bf z}_3$};
\draw[dotted,lightgray](M-6-1.south west)--(M-6-32.south east);
\draw[dotted,lightgray](M-7-1.south west)--(M-7-32.south east);
\draw[dotted,lightgray](M-8-1.south west)--(M-8-32.south east);
\draw[solid,lightgray](M-9-1.south west)--(M-9-32.south east);
\draw[dotted,lightgray](M-10-1.south west)--(M-10-32.south east);
\draw[dotted,lightgray](M-11-1.south west)--(M-11-32.south east);
\draw[dotted,lightgray](M-12-1.south west)--(M-12-32.south east);
\draw[solid,lightgray](M-13-1.south west)--(M-13-32.south east);
\draw[dotted,lightgray](M-14-1.south west)--(M-14-32.south east);
\draw[dotted,lightgray](M-15-1.south west)--(M-15-32.south east);
\draw[dotted,lightgray](M-16-1.south west)--(M-16-32.south east);
\draw[solid,lightgray](M-17-1.south west)--(M-17-32.south east);
\draw[solid,lightgray](M-1-1.north east)--(M-18-1.south east);
\draw[dotted,lightgray](M-1-2.north east)--(M-18-2.south east);
\draw[dotted,lightgray](M-1-3.north east)--(M-18-3.south east);
\draw[dotted,lightgray](M-1-4.north east)--(M-18-4.south east);
\draw[dotted,lightgray](M-1-5.north east)--(M-18-5.south east);
\draw[solid,lightgray](M-1-6.north east)--(M-18-6.south east);
\draw[dotted,lightgray](M-1-7.north east)--(M-18-7.south east);
\draw[dotted,lightgray](M-1-8.north east)--(M-18-8.south east);
\draw[dotted,lightgray](M-1-9.north east)--(M-18-9.south east);
\draw[dotted,lightgray](M-1-10.north east)--(M-18-10.south east);
\draw[-latex,thick,gray](M-1-11.north east)--(M-18-11.south east);
\node (i) at (M-19-12) {$i$};
\draw[dotted,lightgray](M-1-12.north east)--(M-18-12.south east);
\draw[dotted,lightgray](M-1-13.north east)--(M-18-13.south east);
\draw[dotted,lightgray](M-1-14.north east)--(M-18-14.south east);
\draw[dotted,lightgray](M-1-15.north east)--(M-18-15.south east);
\draw[solid,lightgray](M-1-16.north east)--(M-18-16.south east);
\draw[dotted,lightgray](M-1-17.north east)--(M-18-17.south east);
\draw[dotted,lightgray](M-1-18.north east)--(M-18-18.south east);
\draw[dotted,lightgray](M-1-19.north east)--(M-18-19.south east);
\draw[dotted,lightgray](M-1-20.north east)--(M-18-20.south east);
\draw[solid,lightgray](M-1-21.north east)--(M-18-21.south east);
\draw[dotted,lightgray](M-1-22.north east)--(M-18-22.south east);
\draw[dotted,lightgray](M-1-23.north east)--(M-18-23.south east);
\draw[dotted,lightgray](M-1-24.north east)--(M-18-24.south east);
\draw[dotted,lightgray](M-1-25.north east)--(M-18-25.south east);
\draw[solid,lightgray](M-1-26.north east)--(M-18-26.south east);
\draw[dotted,lightgray](M-1-27.north east)--(M-18-27.south east);
\draw[dotted,lightgray](M-1-28.north east)--(M-18-28.south east);
\draw[dotted,lightgray](M-1-29.north east)--(M-18-29.south east);
\draw[dotted,lightgray](M-1-30.north east)--(M-18-30.south east);
\draw[solid,lightgray](M-1-31.north east)--(M-18-31.south east);
\draw[densely dashed, gray](M-14-18)circle(0.4em);
\draw[gray](M-14-19)circle(0.4em);
\begin{scope}[on background layer]
\draw[thin,red](M-18-15.south) -- (M-12-15.center) -| (M-8-17.center) -| (M-7-19.center)  -- (M-7-32.east);
\draw[thin,red](M-18-17.south) -- (M-15-17.center) -| (M-9-18.center) -| (M-8-20.center)  -- (M-8-32.east);
\draw[thin,red](M-18-18.south) -- (M-16-18.center) -| (M-11-19.center) -| (M-10-20.center) -| (M-9-21.center)  -- (M-9-32.east);
\end{scope}
\end{tikzpicture}
\end{center}
{then one can see that $Z$ satisfies the conditions (z.1)-(z.3) in Section \ref{subsec:mbc} with respect to $\tf_2^tA$. Here, the dashed and solid circles are the positions where $\tf_2^tA$ differs from $A$. We conclude that $Z$ induces a proper numbering on $\tf_2^tA$.

It is also easily checked that $\tf_3^t\tf_2^tA = \tf_2^tA - \wh{E}_{2\, 3} + \wh{E}_{2\, 4}$. However, $Z$ does not give a proper numbering on $\tf_3^t\tf_2^tA$, since an inner corner $(2, 3)$ of ${\bf z}_3$ is not a non-zero cell of $\tf_3^t\tf_2^tA$. If we modify a segment $(9,8), (8,8), (7.8). (6.8). (6.9)$ of ${\bf z}_3$ by $(9, 8), (9, 9), (8, 9), (7, 9), (6, 9)$ (see below), it remedies the failure of the condition (z.1). The modified zig-zag, which is denoted by ${\bf z}'_3$, does not intersect with ${\bf z}_4$ as follow,
}

\begin{center}
\begin{tikzpicture}[every node/.style={font=\footnotesize,scale=1}]
\matrix (M)[matrix of math nodes,nodes in empty cells,nodes={rectangle,minimum height=1.0em,minimum width=1.0em,inner sep=0pt,anchor=center,align=center}]
{
& &&&&& &&&&& &&&&& &&&&& &&&&& &&&&& &\\
& &&1&&& &&&1&& &&&&& &&&&& &&&&& &&&&& &\\
& 1&&&&& &&&1&& &&2&&& &&&&& &&&&& &&&&& &\\
& &1&&1&& &&&&& 1&&&1&& &&&&& &&&&& &&&&& &\\
& &&&&& &&&&& &1&&&1& &&&&& &&&&& &&&&& &\\
& &&&&& &&1&&& &&&1&& &&&&& &&&&& &&&&& &\\
& &&&&& 1&&&&& &&&1&& &&2&&& &&&&& &&&&& &\\
& &&&&& &1&&1&& &&&&& 1&&&1&& &&&&& &&&&& &\\
& &&&&& &&&&& &&&&& &1&&&1& &&&&& &&&&& &\\
& &&&&& &&&&& &&1&&& &&&1&& &&&&& &&&&& &\\
& &&&&& &&&&& 1&&&&& &&&1&& &&2&&& &&&&& &\\
& &&&&& &&&&& &1&&1&& &&&&& 1&&&1&& &&&&& &\\
& &&&&& &&&&& &&&&& &&&&& &1&&&1& &&&&& &\\
& &&&&& &&&&& &&&&& &&1&&& &&&1&& &&&&& &\\
& &&&&& &&&&& &&&&& 1&&&&& &&&1&& &&2&&& &\\
& &&&&& &&&&& &&&&& &1&&1&& &&&&& 1&&&1&& &\\
& &&&&& &&&&& &&&&& &&&&& &&&&& &1&&&1& &\\
& &&&&& &&&&& &&&&& &&&&& &&&&& &&&&& &\\
& &&&&& &&&&& &&&&& &&&&& &&&&& &&&&& &\\
};
\draw[solid,lightgray](M-1-1.south west)--(M-1-32.south east);
\draw[dotted,lightgray](M-2-1.south west)--(M-2-32.south east);
\draw[dotted,lightgray](M-3-1.south west)--(M-3-32.south east);
\draw[dotted,lightgray](M-4-1.south west)--(M-4-32.south east);
\draw[-latex,thick,gray](M-5-1.south west)--(M-5-32.south east);
\node (j) at (M-6-33) {$j$};
\node (j) at (M-7-33) {${\bf z}_1$};
\node (j) at (M-8-33) {${\bf z}_2$};
\node (j) at (M-9-33) {${\bf z}'_3$};
\node (j) at (M-11-33) {${\bf z}_4$};
\draw[dotted,lightgray](M-6-1.south west)--(M-6-32.south east);
\draw[dotted,lightgray](M-7-1.south west)--(M-7-32.south east);
\draw[dotted,lightgray](M-8-1.south west)--(M-8-32.south east);
\draw[solid,lightgray](M-9-1.south west)--(M-9-32.south east);
\draw[dotted,lightgray](M-10-1.south west)--(M-10-32.south east);
\draw[dotted,lightgray](M-11-1.south west)--(M-11-32.south east);
\draw[dotted,lightgray](M-12-1.south west)--(M-12-32.south east);
\draw[solid,lightgray](M-13-1.south west)--(M-13-32.south east);
\draw[dotted,lightgray](M-14-1.south west)--(M-14-32.south east);
\draw[dotted,lightgray](M-15-1.south west)--(M-15-32.south east);
\draw[dotted,lightgray](M-16-1.south west)--(M-16-32.south east);
\draw[solid,lightgray](M-17-1.south west)--(M-17-32.south east);
\draw[solid,lightgray](M-1-1.north east)--(M-18-1.south east);
\draw[dotted,lightgray](M-1-2.north east)--(M-18-2.south east);
\draw[dotted,lightgray](M-1-3.north east)--(M-18-3.south east);
\draw[dotted,lightgray](M-1-4.north east)--(M-18-4.south east);
\draw[dotted,lightgray](M-1-5.north east)--(M-18-5.south east);
\draw[solid,lightgray](M-1-6.north east)--(M-18-6.south east);
\draw[dotted,lightgray](M-1-7.north east)--(M-18-7.south east);
\draw[dotted,lightgray](M-1-8.north east)--(M-18-8.south east);
\draw[dotted,lightgray](M-1-9.north east)--(M-18-9.south east);
\draw[dotted,lightgray](M-1-10.north east)--(M-18-10.south east);
\draw[-latex,thick,gray](M-1-11.north east)--(M-18-11.south east);
\node (i) at (M-19-12) {$i$};
\draw[dotted,lightgray](M-1-12.north east)--(M-18-12.south east);
\draw[dotted,lightgray](M-1-13.north east)--(M-18-13.south east);
\draw[dotted,lightgray](M-1-14.north east)--(M-18-14.south east);
\draw[dotted,lightgray](M-1-15.north east)--(M-18-15.south east);
\draw[solid,lightgray](M-1-16.north east)--(M-18-16.south east);
\draw[dotted,lightgray](M-1-17.north east)--(M-18-17.south east);
\draw[dotted,lightgray](M-1-18.north east)--(M-18-18.south east);
\draw[dotted,lightgray](M-1-19.north east)--(M-18-19.south east);
\draw[dotted,lightgray](M-1-20.north east)--(M-18-20.south east);
\draw[solid,lightgray](M-1-21.north east)--(M-18-21.south east);
\draw[dotted,lightgray](M-1-22.north east)--(M-18-22.south east);
\draw[dotted,lightgray](M-1-23.north east)--(M-18-23.south east);
\draw[dotted,lightgray](M-1-24.north east)--(M-18-24.south east);
\draw[dotted,lightgray](M-1-25.north east)--(M-18-25.south east);
\draw[solid,lightgray](M-1-26.north east)--(M-18-26.south east);
\draw[dotted,lightgray](M-1-27.north east)--(M-18-27.south east);
\draw[dotted,lightgray](M-1-28.north east)--(M-18-28.south east);
\draw[dotted,lightgray](M-1-29.north east)--(M-18-29.south east);
\draw[dotted,lightgray](M-1-30.north east)--(M-18-30.south east);
\draw[solid,lightgray](M-1-31.north east)--(M-18-31.south east);
\draw[dashed, gray](M-11-19)circle(0.4em);
\draw[gray](M-11-20)circle(0.4em);
\begin{scope}[on background layer]
\draw[thin,red](M-18-15.south) -- (M-12-15.center) -| (M-8-17.center) -| (M-7-19.center)  -- (M-7-32.east);
\draw[thin,red](M-18-17.south) -- (M-15-17.center) -| (M-9-18.center) -| (M-8-20.center)  -- (M-8-32.east);
\draw[thin,red](M-18-18.south) -- (M-16-18.center) -| (M-14-19.center) -| (M-10-20.center) -| (M-9-21.center)  -- (M-9-32.east);
\draw[thin,red](M-18-20.south) -- (M-16-20.center) -| (M-12-22.center) -| (M-11-24.center)  -- (M-11-32.east);
\draw[densely dashed,red](M-14-19.center) |- (M-11-20.center);

\end{scope}
\end{tikzpicture}
\end{center}
{
where dashed red line is the original part of ${\bf z}_3$, and the dashed circle and solid circle are the positions where $\tf_3^t\tf_2^tA$ differs from $\tf_2^tA$. Hence the modified set of zig-zags $Z'=\{\, \cdots, {\bf z}'_0, {\bf z}_1, {\bf z}_2, {\bf z}'_3, {\bf z}_4, \cdots\, \}$ give a proper numbering on $\tf_3^t\tf_2^tA$.

In the remainder of this section, we will see that the induced numberings on $\tf_2^tA$ and $\tf_3^t\tf_2^tA$ are, in fact, the southwest channel numberings.}
}
\end{ex}\vskip 2mm

{The following lemma describes how to construct a proper numbering on $\td{A}$ from a given numbering $d$ on $A$ in general.}

\begin{lem}\label{lem:tilde of numbering}
We have the following.
\begin{itemize}
    \item[(1)] Let $d$ be a proper numbering on $A$ with the associated zig-zags $Z = \{ {\bf z}_k \}_{k \in \Z}$ and $d(u, j) = s$. Then there exists a proper numbering ${d}^-$ on $\td{A}$ satisfying
    \begin{equation*}\label{eq:d1 tilde}
    \begin{split}
     {d}^-(c) & = d(c)  \quad\text{ if $c \in {\rm supp}(A) \cap {\rm supp}(\td{A})$}, \\
     {d}^-(u, j+1) & =  
      \begin{cases}
       d(u, j)  & \text{if $u$ is minimal such that $(u, j)\in {\bf z}_{s}$} ,\\
       d(u, j) + 1 & \text{otherwise}.
      \end{cases}
    \end{split}
    \end{equation*}
    \item[(2)] Let $d$ be a proper numbering on $\td{A}$ with the associated zig-zags $Z = \{ {\bf z}_k \}_{k \in \Z}$ and $d(u, j+1) = t$. Then there exists a proper numbering $d^+$ on $A$ satisfying
    \begin{equation*}\label{eq:d2 tilde}
    \begin{split}
        d^+(c) & = d(c)  \quad\text{ if $c \in {\rm supp}(A) \cap {\rm supp}(\td{A})$}, \\
        d^+(u, j) & =  
        \begin{cases}
         d(u, j) - 1 & \text{if $(u, j) \in {\bf z}_{t-1}$} ,\\
         d(u, j) & \text{otherwise}.
        \end{cases}
    \end{split}
    \end{equation*}
\end{itemize}
In particular, the widths of $A$ and $\td{A}$ are the same.
\end{lem}
\pf 
(1) We construct a set of zig-zags $Z^-$ (by adjusting $Z$) which {satisfies} the conditions (z.1)-(z.3) in Section \ref{subsec:mbc} with respect to $\td{A}$ and hence gives a proper numbering $d^-$ on $\td{A}$.

Let $i_{s}$ be the minimal row index with $(i_{s}, j) \in {\bf z}_{s}$. 
If $i_{s}=u$, then $(u, j+1) \in {\bf z}_{s}$ by definition of zig-zag. 
Suppose that $i_{s}< u$. 
Then $(i_{s}, j)$ is an inner corner of ${\bf z}_{s}$, and $a_{i_{s}\,j} >0$. Consider a subsequence of $\sigma$ in \eqref{eq:sigma ast}
\begin{equation}\label{eq:sigma1}
\begin{split}
 & (\underbrace{+\dots+}_{a_{i_{s} j}},\ \underbrace{-\dots -}_{a_{i_{s}+1\, j+1}},\ \cdots,\ \underbrace{+\dots +}_{a_{u-1\, j}},\ \underbrace{-\dots -}_{a_{u\, j+1}}).
\end{split}
\end{equation}
Since $(u, j)$ is the cell corresponding to the leftmost $+$ in $\td{\sigma}$, there exists no $+$ in the reduced form of \eqref{eq:sigma1}. 
This implies that there exists some $a_{v\, j+1} > 0$ for some $v$ with $i_{s} < v \le u$ so that $+$ in the cell $(i_{s}, j)$ is paired with $-$ in $(v, j+1)$. 
It is easy to see that $d(v, j+1) = s + 1$, and hence $(u, j+1) \in {\bf z}_{s+1}$. 

Hence we see that ${\rm supp}(\td{A}) \subseteq \bigsqcup_{k \in \Z} {\bf z}_k$ and $Z$ satisfies the conditions in (z.2) and (z.3) for $\td{A}$.
Note that the condition (z.1) fails if and only if $(u,j)$ is an inner corner of ${\bf z}_{s}$ with $a_{uj} = 1$. \vskip 2mm

\noindent{\em Case 1}. Suppose that $(u, j)$ is not an inner corner of ${\bf z}_{s}$ or $(u, j)$ is an inner corner of ${\bf z}_{s}$  with $a_{uj}>1$.
Then $Z^-:=Z$ satisfies the condition (z.1), and induces a proper numbering ${d}^-$ on $\td{A}$ as given in (\ref{eq:numbering defined by zig-zags}).   
Hence ${d}^-$ satisfies
\begin{equation*}
    {d}^-(u, j+1) = 
    \begin{cases}
     s  & \text{if $u = i_{s}$} ,\\
     s + 1 & \text{if $i_{s}<u$},
    \end{cases}
\end{equation*} 
and ${d}^-(c) = d(c)$ for $c \in {\rm supp}(A) \cap {\rm supp}(\td{A})$.\vskip 2mm

\noindent{\em Case 2}. Suppose that $(u, j)$ is an inner corner of ${\bf z}_{s}$ with $a_{uj} = 1$. 
In this case, the condition (z.1) fails since the inner corner $(u, j)$ of ${\bf z}_{s}$ does not lie in ${\rm supp}(\td{A})$. 
Now let us modify ${\bf z}_{s}$ as follows: Let $(v, j)$ be an outer corner of ${\bf z}_{s}$ and let 
\begin{equation*}
    w = {\rm min}\{\, i \in {\Z} \,\vert \, \text{$u < i \le v$ and $a_{ij}>0$}\, \}.
\end{equation*} 
Note that $w = u^*$ when $d = d^{\tt sw}_A$.
Consider a subsequence of $\sigma$ 
\begin{equation}\label{eq:sigma2}
\begin{split}
& (\underbrace{+}_{a_{u j}},\ \underbrace{-\dots -}_{a_{u+1\, j+1}},\ \cdots\,\ \underbrace{+\dots +}_{a_{w-1\, j}},\ \underbrace{-\dots -}_{a_{w\, j+1}}).
\end{split}
\end{equation}
We see that $a_{ij} = 0$ for $u < i \le w-1$ by definition of $w$, and moreover $a_{ij+1} = 0$ for $u < i \le w$ since the reduced form of \eqref{eq:sigma2} is $(+)$. Indeed the sequence in \eqref{eq:sigma2} is $(+)$. According to this observation, we define a zig-zag ${\bf z}^-_{s}$ by replacing the cells
\begin{equation*}
    (w - 1, j), \, \cdots, \, (u+1, j), \,(u, j)
\end{equation*}
in ${\bf z}_{s}$ with the following cells
\begin{equation*}
    (w, j+1), \, (w - 1, j+1), \, \cdots, \, (u+1, j+1).
\end{equation*} Then each inner corner of ${\bf z}^-_{s}$ lies in ${\rm supp}(\td{A})$.

Let $Z^-$ be the set of zig-zags obtained from $Z$ by replacing $\{\,\tau^k{\bf z}_{s}\,|\,k\in\Z\,\}$ with $\{\,\tau^k{\bf z}^-_{s}\,|\,k\in\Z\,\}$. 
Then $Z^-$ satisfies the conditions (z.1)-(z.3) for $\td{A}$, and hence induces a proper numbering ${d}^-$ on $\td{A}$. 
It is easy to see that ${d}^-(u, j+1) = s$ since $u= i_{s}$, and that ${d}^-(c) = d(c)$ for $c \in {\rm supp}(A) \cap {\rm supp}(\td{A})$. 

By definition, the proper numberings $d$ and $d^-$ have the same flow, which implies that $A$ and $\td{A}$ have the same width.
This proves (1).

(2) As in (1), we construct a set of zig-zags from $Z$ to which a proper numbering $d^+$ on $A$ is associated.\vskip 2mm

\noindent{\em Case 1}. Suppose first that $(u, j) \in {\bf z}_{t-1}$.  
Then there exists an inner corner $(v, j)$ of ${\bf z}_{t-1}$ and an inner corner $(w, j+1)$ of ${\bf z}_{t}$ with $v < w \le u$. In particular, we have $\td{a}_{v j} > 0$. Consider a subsequence of $\sigma$
\begin{equation*}
\begin{split}
(\underbrace{+\dots+}_{\td{a}_{v j}},\,  \underbrace{-\dots-}_{\td{a}_{v - 1\, j+1}},\, \cdots,\, \underbrace{+\dots +}_{\td{a}_{u-1\, j}},\, \underbrace{-\dots -}_{\td{a}_{u\, j+1}}).
\end{split}
\end{equation*}
Since $\td{a}_{v j}>0$ and $(u, j+1)$ is the cell corresponding to the rightmost $-$, we have 
\begin{equation*}
    \td{a}_{w\, j+1} + \td{a}_{w+1\, j+1} + \cdots + \td{a}_{u\, j+1} \ge 2.
\end{equation*} 
This shows that $w<u$ and the inner corner $(w, j+1)$ of ${\bf z}_{t}$ lies in ${\rm supp}(A)$. 
Hence $Z^+:=Z$ satisfies the conditions (z.1)-(z.3), and it induces a proper numbering $d^+$ on $A$. 
It is obvious that $d^+(u, j) = t-1$ and $d^+(c) = d(c)$ for $c \in {\rm supp}(A) \cap {\rm supp}(\td{A})$.\vskip 2mm

\noindent{\em Case 2}.  Suppose that $(u, j) \notin {\bf z}_{t-1}$. 
Let $i_{t}$ be the maximal row index with $(i_{t}, j+1) \in {\bf z}_{t}$. 
Since $(u, j) \notin {\bf z}_{t-1}$, we have $\td{a}_{ij}=0$ for $u \le i < i_{t}$. Consider a subsequence of $\sigma$
\begin{equation}\label{eq:sigma3}
\begin{split}
(\underbrace{-\dots-}_{\td{a}_{u\,j+1}},\ \underbrace{+\dots+}_{\td{a}_{u j}},\ \cdots ,\ \underbrace{-\dots -}_{\td{a}_{i_{t}\, j+1}}) = 
(\underbrace{-\dots -}_{\td{a}_{u\, j+1}},\  \cdots ,\ \underbrace{-\dots -}_{\td{a}_{i_{t} \,j+1}}).
\end{split}
\end{equation} 
Since the rightmost $-$ in $\td{\sigma}$ corresponding to position $(u, j+1)$ is the one in \eqref{eq:sigma3}, we see that $\td{a}_{i\, j+1}=0$ for $u < i \le i_{t}$. 
We define ${\bf z}^+_{t}$ to be a zig-zag by replacing the cells
\begin{equation*}
    (i_{t}, j+1), \, (i_{t} - 1, j+1), \, \cdots, \, (u+1, j+1)
\end{equation*}
in ${\bf z}_{t}$ with the following cells
\begin{equation*}
    (i_{t} - 1, j), \, \cdots, \, (u-1, j),\, (u, j).
\end{equation*} 
Then $(u, j)$ is an inner corner of ${\bf z}^+_{t}$, and $(u, j) \in {\rm supp}(A)$.
Let $Z^+$ be the set of zig-zags obtained from $Z$ by replacing $\{\,\tau^k{\bf z}_{t}\,|\,k\in\Z\,\}$ with $\{\,\tau^k{\bf z}^+_{t}\,|\,k\in\Z\,\}$. 
Then $Z^+$ satisfies the conditions (z.1)-(z.3) for $A$, and it induce a proper numbering $d^+$ on $A$. We have $d^+(u, j) = t$ and $d^+(c) = d(c)$ for $c \in {\rm supp}(A) \cap {\rm supp}(\td{A})$. This proves (2). 
\qed 

\begin{rem}\label{rem:properties of d tilde}
{\rm
Let $d_1$ be a proper numbering on $A$. 
If follows from the construction of $d^\pm$ in the proof of Lemma \ref{lem:tilde of numbering} that
\begin{itemize}
    \item[(1)] $(d_1^-)^+ = d_1$,
    \item[(2)] if $d_2$ is another proper numbering on $A$ such that $d_1(c) \le d_2(c)$ for $c \in {\rm supp}(A)$, then $d_1^-(c) \le d_2^-(c)$ for $c \in {\rm supp}(\td{A})$.
\end{itemize}
The similar properties hold for a proper numbering on $\td{A}$.
}
\end{rem}

\begin{lem}\label{lem:tilde of channel}
We have the following.
\begin{itemize}
\item[(1)] Let $C$ be a channel of $A$. Then there exists a channel $C^-$ of $\td{A}$ given by
\begin{equation*}
\begin{split}
\begin{cases}
\big(C - (u,j)^\wedge \big) \cup (u,j+1)^\wedge & 
\text{if $(u,j) \in C$, $a_{uj}=1$ and $(i,j+1) \not\in C$ for all $i$},\\
\big(C - (u,j)^\wedge \big) \cup (v,j)^\wedge & 
\text{if $(u,j)\in C$, $a_{uj}=1$ and $(i,j+1) \in C$ for some $i$},\\
\ C & \text{otherwise},
\end{cases}
\end{split}
\end{equation*}
where $v$ is the minimal row index such that $u<v$ and $a_{v j} >0$ if it exists.
    
\item[(2)]
Let $C$ be a channel of $\td{A}$. Then there exists a channel $C^+$ of $A$ given by
\begin{equation*}
\begin{split}
\begin{cases}
\big(C - (u,j+1)^\wedge \big) \cup (u,j)^\wedge & 
\text{if $(u,j+1) \in C$, $a_{u\,j+1}=1$ and $(i,j) \not\in C$ for all $i$},\\
\big(C - (u,j+1)^\wedge \big) \cup (w,j+1)^\wedge & 
\text{if $(u,j+1)\in C$, $a_{u\,j+1}=1$ and $(i,j) \in C$ for some $i$},\\
\ C & \text{otherwise},
\end{cases}
\end{split}
\end{equation*}
where $w$ is the maximal row index such that $w<u$ and $a_{w\, j+1} >0$ if it exists.
\end{itemize}
\end{lem}
\pf Let us prove (1) only, since the proof of (2) is similar. \vskip 2mm

First, suppose $(u, j) \notin C$ or $a_{uj} > 1$. 
Then $C \subseteq {\rm supp}(\td{A})$ and $C$ is a channel of $\td{A}$. 
We put $C^- = C$ in this case.

Now suppose that $(u, j) \in C$ and $a_{uj} = 1$. 
Then we have $C \nsubseteq {\rm supp}(\td{A})$. Let us write
\begin{equation*}
 C = \{ \cdots >_{\tt NW} c_{s-1} >_{\tt NW} c_{s} >_{\tt NW} c_{s+1} >_{\tt NW} \cdots\}
\end{equation*}
with $c_{s} = (u, j)$. {We have two cases.}\vskip 2mm

\noindent{\em Case 1.} Suppose that ${(i,j+1)} \not\in C$ for all $i$.
Let $C^-$ be a set obtained from $C$ by replacing $c_s^\wedge=(u,j)^\wedge\subset C$ with $(u,j+1)^\wedge$, which is clearly a stream of $\td{A}$ by assumption.\vskip 2mm

\noindent{\em Case 2.} Suppose that $(i,j+1) \in C$ for some $i$.
Then we have $c_{s+1} = (u', j+1)$ for some $u'>u$. 
Consider a subsequence of $\sigma$ in \eqref{eq:sigma ast}
\begin{equation}\label{eq:sigma4}
\begin{split}
(\underbrace{+}_{a_{uj}},\, \underbrace{-\dots-}_{a_{u+1\, j+1}},\, \cdots,\, \underbrace{+\dots+}_{a_{u'-1\, j}},\, \underbrace{-\dots -}_{a_{u'\, j+1}}).
\end{split}
\end{equation} 
Since $-$ in $c_{s+1}$ is paired with $+$ in \eqref{eq:sigma4} other than $+$ in $c_s=(u, j)$, we have $a_{ij}>0$ for some  $u<i<u'$. Let $v$ be the minimal such one. Then we have
\begin{equation}\label{eq:flow C^-}
c_{s-1} >_{\tt NW} (v, j+1) >_{\tt NW} c_{s+1}.
\end{equation} 
Let $C^-$ be a set obtained from $C$ by replacing $c_s^\wedge=(u,j)^\wedge\subset C$ with $(v,j+1)^\wedge$, which is a stream of $\td{A}$ by \eqref{eq:flow C^-}. 

By definition, $C^-$ has the same flow as $C$. Since $A$ and $\td{A}$ have the same width by Lemma \ref{lem:tilde of numbering}, $C^-$ is a stream of maximal flow, and hence a channel of $\td{A}$ (cf.~Remark \ref{rem:period=flow=width}).
\qed

\begin{rem}\label{rem:properties of c tilde}
{\rm
Let $C$ be channel of $A$. Under the above hypothesis, we have
\begin{equation*}
(C^-)^+ - C =
\begin{cases}
(v,j) & \text{ if $(u,j) \in C$, $a_{uj}=1$ and $(i,j+1) \in C$ for some $i$},\\
(u,j+1) & \text{ if $(u,j) \in C$, $a_{uj}=1$, $a_{u\,j+1}> 0$ and $(i,j+1) \notin C$ for all $i$},\\
\emptyset & \text{ otherwise}.
\end{cases}
\end{equation*} 
Note that there exists no cell in ${\rm supp}(A)$ between $(u, j)$ and $(v, j)$ and between $(u, j)$ and $(u, j+1)$. 
Hence it follows that
\begin{equation*}
C \succcurlyeq_{\tt sw} (C^-)^+\quad \text{ or }\quad (C^-)^+ \succcurlyeq_{\tt sw} C,
\end{equation*}
and there exists no other channel between $C$ and $(C^-)^+$. 
It is also easy to check that if $C'$ is another channel of $A$ with $C{\succcurlyeq_{\tt sw}} C'$, then we have 
\begin{equation}\label{eq:ineq sw and pm}
C^- \succcurlyeq_{\tt sw} (C')^-.
\end{equation}

The similar properties also hold with respect to channels of $\td{A}$.
}
\end{rem}

\begin{lem}\label{lem:dtilde for channel numbering}
We have the following.
\begin{itemize}
 \item[(1)] Let $C$ be a channel of $A$. If $d$ is the channel numbering on $A$ associated to $C$, then $d^-$ is the channel numbering on $\td{A}$ associated to $C^-$.
 
 \item[(2)] Let $C$ be a channel of $\td{A}$. If $d$ is the channel numbering on $\td{A}$ associated to $C$, then $d^+$ is the channel numbering on $A$ associated to $C^+$.
\end{itemize}

\end{lem}
\pf Let us prove (1) only since the proof of (2) is similar. 

Let $d'$ be the channel numbering on $\td{A}$ associated to $C^-$. Since the widths of $A$ and $\td{A}$ coincide, we may assume that ${d}^-$ coincides with $d'$ on the channel $C^-$ by adding a constant to $d^-$. 
Hence we have $d' \le {d}^-$ by Lemma \ref{lem:characterization of ch numbering}, and $(d')^+ \le ({d}^-)^+=d$ by Remark \ref{rem:properties of d tilde}.

Let $\ell$ be the common width of $A$ and $\td{A}$. If $\ell>1$, then we see from Lemma \ref{lem:tilde of channel}(1) that there exists $c \in C \cap {C}^-$ such that $c \in {\rm supp}(A) \cap {\rm supp}(\td{A})$. 
We have 
$$d(c)={d}^-(c)=d'(c)=(d')^+(c).$$ 
So $(d')^+$ and $d$ also coincide on $C$ (cf.~Remark \ref{rem:properties of c tilde}), and $d \le (d')^+$ by Lemma \ref{lem:characterization of ch numbering}.
Therefore, we have $d = (d')^+$, and $d^- = d'$ by Remark \ref{rem:properties of d tilde}.

If $\ell=1$, then it is not possible to have ${d}^-(u, j+1) = d(u, j) + 1$ or $C^-=\big(C - (u,j)^\wedge \big) \cup (v,j)^\wedge$ since we must have another cell $(i, j+1)$ with $(u, j) >_{\tt NW} (i, j+1) >_{\tt NW} (u,j) + (m,n)$. 
Hence we see directly that $(d')^+(u,j)=d(u,j)={d}^-(u,j+1)=d'(u,j+1)$.
By similar arguments as in the above case, we conclude that $d^-=d'$.
\qed\newline

Now we can describe the southwest channel numbering on $\td{A}$ in terms of the one on $A$.

\begin{prop}\label{prop:swn of A tilde}
Let $d$ be the southwest channel numbering on $A$. Then $d^-$ is the southwest channel numbering on $\td{A}$. Equivalently, let $d$ be the southwest channel numbering on $\td{A}$. Then $d^+$ is the southwest channel numbering on $A$. 
\end{prop}
\pf
Let $C_1=C^{\tt sw}_A$ and $C_2 = C^{\tt sw}_{\td A}$. Let $d'$ be the southwest channel numbering on $\td{A}$.
By Lemma \ref{lem:dtilde for channel numbering}, we have 
$$
{d}^- = d^{{C}_1^-}_{\td{A}},\quad (d')^+ = d^{{C}_2^+}_A.
$$
Thus it suffices to show that either $C_1 = C_2^+$ or $C_1^- = C_2$, which implies that $d=(d')^+$ or ${d}^- = d'$, respectively (see also Remark \ref{rem:properties of d tilde}).

Since $C_1$ and $C_2$ are the southwest channels, we have
\begin{equation}\label{eq:sw order on two swc}
    C_1 \succcurlyeq_{\tt sw} C^+_2,\quad  C_2 \succcurlyeq_{\tt sw} C^-_1.
\end{equation}
By Remark \ref{rem:properties of c tilde} (cf.~\eqref{eq:ineq sw and pm}), we get from \eqref{eq:sw order on two swc}
\begin{equation}\label{eq:ineq for channel pm}
    C_1 \succcurlyeq_{\tt sw} C_2^+ \succcurlyeq_{\tt sw} (C_1^-)^+,\quad 
    C_2 \succcurlyeq_{\tt sw} C_1^- \succcurlyeq_{\tt sw}  (C_2^+)^-.
\end{equation}

We claim that $C_1 = C^+_2$ if $C_2 \succ_{\tt sw} C_1^-$.
By \eqref{eq:ineq for channel pm}, we have $C_2 \succ_{\tt sw} (C^+_2)^-$.
By Remark \ref{rem:properties of c tilde}, we see that $C_2 \succ_{\tt sw} (C^+_2)^-$ occurs only when 
$$C_2^+=\big(C_2 - (u,j+1)^\wedge \big) \cup (w,j+1)^\wedge.$$ 
Then we have $C_2^+ = (C_2^+)^-$. 

On the other hand, we have $C^-_1 = (C_2^+)^-$ since there is no other channel between $C_2$ and $(C_2^+)^-$. 
Hence we get $C^-_1 = (C_2^+)^- = C^+_2$, and in particular $(w, j+1) \in C^-_1$. Since we have $(u, j) \notin C_1$ by Lemma \ref{lem:tilde of channel}(1), it follows that $C_1 = C^-_1 = C^+_2$. This proves the claim. \qed

\subsection{Proof of \eqref{eq:main equation}}\label{subsec:proof of main commutativity}
Now we are in a position to prove \eqref{eq:main equation}.
Let $d=d_A^{\tt sw}$ and let $Z = \{ {\bf z}_k \}_{k \in \Z}$ be the set of zig-zags associated to $d$. 
Let $Z^-$ be the set of zig-zags associated to $d^-$ (see the proof of Lemma \ref{lem:tilde of numbering}(1)).
Note that $d^- = d_{\td A}^{\tt sw}$ by Proposition \ref{prop:swn of A tilde}.
\vskip 2mm

\noindent{\em Case 1.} Suppose that $(u, j)$ is not an inner corner of ${\bf z}_{s}$ or $(u, j)$ is an inner corner of ${\bf z}_{s}$  with $a_{uj}>1$.
Since $Z=Z^-$ and the cells corresponding the leftmost $+$ in $A$ and ${A}^*$ coincide in this case, we have $\td{{A}^*}= (\td{A})^*$. \vskip 2mm

\noindent{\em Case 2.} Suppose that $(u, j)$ is an inner corner of ${\bf z}_{s}$  with $a_{uj}=1$.

Let us first compare ${\bf z}_{s}$ with the modified zig-zag ${\bf z}^-_{s}$. 
Let $(v_0, j)$ be an outer corner of ${\bf z}_{s}$ and let $v_1$ be the minimal row index with $(v_1, j+1) \in {\bf z}_{s}$. 
Note that the inner and outer corners of $Z$ and $Z^-$ always coincide other than the following cells (more precisely, their orbits under $\tau^{\pm 1}$):
\begin{equation}\label{eq:cells important}
(u, j),\, (u^*, j),\, (v_0, j),\, (v_1, j+1),\, (u, j+1),\, (u^*, j+1),
\end{equation} 
where $v_1 \le u < u^* \le v_0 \le \infty$. 
For the reader's convenience, we summarize the positions of the cells in \eqref{eq:cells important} as follows:
\begin{center}
\begin{tabular}{c|c|c}
    & ${\bf z}_{s}$ & ${\bf z}^-_{s}$ \\
     \hline
    $(u, j)$   & inner corner & \\
    $(u^*, j)$ & outer corner if $u^* = v_0$ & inner corner if $u^*< v_0$\\
    $(v_0, j)$ & outer corner& outer corner if $u^* < v_0$\\ \hline
    $(v_1, j+1)$ & inner corner if $v_1 < u$ & inner corner\\ 
    $(u, j+1)$ & outer corner if $v_1 < u$ & inner corner if $v_1 = u$\\
    $(u^*, j+1)$ &  & outer corner
\end{tabular}
\end{center}
Hence we may write
\begin{equation}\label{eq:A*-2}
\begin{split}
    {A}^* & = A - \wh{E}_{uj} + \wh{E}_{v_0 j} - \wh{E}_{v_1\, j+1} + \wh{E}_{u\, j+1} +B, \\
    (\td{A})^* & = \td{A} - \wh{E}_{u^*j} + \wh{E}_{v_0 j} - \wh{E}_{v_1\, j+1} + \wh{E}_{u^*\, j+1}  +B,
\end{split}
\end{equation}
where $B$ is a finite linear combination of $\wh{E}_{kl}$'s over the cells $(k,l)$ not belonging to \eqref{eq:cells important}.
Combining \eqref{eq:A*-1} and \eqref{eq:A*-2}, we have
\begin{equation*}
\begin{split}
    \td{{A}^*} 
    &= {A}^* - \wh{E}_{u^*j} + \wh{E}_{u^*\,j+1} \\
    &= \left(A - \wh{E}_{uj} + \wh{E}_{v_0 j} - \wh{E}_{v_1\, j+1} + \wh{E}_{u\, j+1} +B\right) - \wh{E}_{u^*j} + \wh{E}_{u^*\,j+1} \\
    &= \left(A - \wh{E}_{uj} + \wh{E}_{u\, j+1}\right) - \wh{E}_{u^*j} + \wh{E}_{v_0 j} - \wh{E}_{v_1\, j+1} + \wh{E}_{u^*\, j+1} + B \\
     & = \td{A} - \wh{E}_{u^*j} + \wh{E}_{v_0 j} - \wh{E}_{v_1\, j+1} + \wh{E}_{u^*\, j+1} + B\\
     & = (\td{A})^*.
\end{split}
\end{equation*} 
By {\em Case 1} and {\em Case 2}, we have $\td{{A}^*}=(\td{A})^*$. Let us write $\Psi(\td{f}_jA) = (\td{f}_jA)^{\flat} \otimes {\bf s}'$. From $\td{A^*} = (\td{A})^*$ and Proposition \ref{prop:swn of A tilde}, we see that
\begin{equation*}
(\td{f}_jA)^{\flat} \otimes {\bf s}' =
\begin{cases}
A^\flat \otimes \td{f}_j{\bf s} & \text{ if $u^* = \infty$,} \\
(\td{f}_jA^{\flat}) \otimes {\bf s} & \text{ if $u^* < \infty.$}
\end{cases}
\end{equation*}
Comparing this with \eqref{eq:Psi and f_j}, we have \eqref{eq:main equation}.

By \eqref{eq:Psi composition}, $\kappa$ commutes with $\tf_j$, and hence $\kappa$ commutes with $\te_j$ for  $j\in \{\,0,1,\dots,n-1\,\}$.

\subsection{Proof of Theorem \ref{thm:main-2}}
We have proved that $\kappa$ commutes with $\te^t_j$ and $\tf^t_j$ for $j\in \{\,0,1,\dots,n-1\,\}$.
Let us finish the proof of Theorem \ref{thm:main-2} by showing that $\kappa$ commutes with $\te_i$ and $\tf_i$ for $i\in \{\,1,\dots,m-1\,\}$.

First, it is not difficult to see that Proposition \ref{prop:swn of A tilde} still holds if we replace the southwest channel numberings with the northeast channel numberings (cf. Remark \ref{rem:ne ch}). Hence by the same arguments as Section \ref{subsec:proof of main commutativity} we have 
\begin{equation}\label{eq:x_i commutes with *}
 \td{x}_i{A}^*=(\td{x}_iA)^*
\end{equation}
for $i\in \{\,0,1,\dots,m-1\,\}$ and $x\in \{\, e,f\, \}$.
This implies that $\Psi$ commutes with $\te_i$ and $\tf_i$ for $i\in \{\,1,\dots,m-1\,\}$ (see Remark \ref{rem:e_0 and f_0} for $i=0$). Hence $\kappa$ commutes with $\te_i$ and $\tf_i$ for $i\in \{\,1,\dots,m-1\,\}$.
This completes the proof of Theorem \ref{thm:main-2}. \qed

\begin{rem}\label{rem:e_0 and f_0}
{\rm 
We should remark that $\te_0$ and $\tf_0$ may not commute with $\Psi$. 
Let $A\in \wh{\mc M}_{m\times n}$ be given such that $\td{x}_0A\neq {\bf 0}$ for $x\in\{e,f\}$. Suppose that $\Psi(A)=A^\flat\ot {\bf s}$ and ${\bf s}=(\ba,\bb, r)$.

If $\td{x}_0(A^\flat\ot {\bf s})=A^\flat\ot \td{x}_0{\bf s}$, then it follows from \eqref{eq:x_i commutes with *} that $\Psi(\td{x}_0A)=A^\flat \ot {\bf s}'$ and ${\bf s}'=(\td{x}_0\ba,\bb, r')$, where $r'=r+1$ (resp. $r-1$) if $x=e$ (resp. $x=f$). Since $\td{x}_0{\bf s}=(\td{x}_0\ba,\bb, r)\neq {\bf s}'$, we have $\Psi(\td{x}_0A)\neq \td{x}_0\Psi(A)$. 
If $\kappa(A)=(P_0,Q)$, then by applying $\Psi$ repeatedly we have 
\begin{equation*}
\kappa(\td{x}_0A)=(\td{x}_0P_0,Q'),
\end{equation*}
for some $Q'\in \mc{B}_n(\la)$ with $Q'\neq Q$.
} 
\end{rem}

\section{Dual affine RSK correspondence}\label{sec:dual aff RSK}
In this section, we construct a dual analogue of Theorem \ref{thm:main-1}.

\subsection{Generalized dual affine permutations}\label{subsec:notations for dual}
Let
\begin{equation*}
\wh{\mc{N}}_{m\times n}=
\left\{\,
A = (a_{ij})_{i,j \in \Z}\,\,\Bigg\vert\,\, 
\begin{array}{l}
(1)\ \text{$a_{ij} \in \{0, 1\}$ and $a_{i+m\, j+n} = a_{ij}$ for all $i,j\in\Z$},\\
(2)\ \text{for each $j$, $a_{ij}=0$ except for finitely many $i$'s}.
\end{array}
\,\right\}.
\end{equation*}
For $A \in \wh{\mathcal{N}}_{m \times n }$,
we define a {\em (dual) standardization of $A$} to be
\begin{equation*}
 A^{\tt st'} = (A^\bullet)^{\circ'}, 
\end{equation*}
(see Remark \ref{rem:st other versions} for the definition of $A^\bullet$).
We may adopt the same notations in Section \ref{subsec:standardization} for $A^{\tt st'}$.
Note that a cell $c$ in ${\rm supp}(A)$ corresponds to a unique cell in ${\rm supp}(A^{\tt st'}_c)$.

\begin{ex}\label{ex:Adual}{\rm
Let  $A \in \wh{\mc{N}}_{3 \times 4}$ be given as follows.\vskip 2mm
\begin{center}
\begin{tikzpicture}[every node/.style={font=\footnotesize,scale=1}]
\matrix (M)[matrix of math nodes,nodes in empty cells,nodes={rectangle,minimum height=1.0em,minimum width=1.0em,inner sep=0pt,anchor=center,align=center}]
{
& &&&& &&&& &&&& &&&& &\\
& 1&&1&& 1&&&& &&&& &&&& &\\
& &1&&& &1&&& &&&& &&&& &\\
& &1&&1& &&1&& &&&& &&&& &\\
& &&&& 1&&1&& 1&&&& &&&& &\\
& &&&& &1&&& &1&&& &&&& &\\
& &&&& &1&&1& &&1&& &&&& &\\
& &&&& &&&& 1&&1&& 1&&&& &\\
& &&&& &&&& &1&&& &1&&& &\\
& &&&& &&&& &1&&1& &&1&& &\\
& &&&& &&&& &&&& &&&& &\\
& &&&& &&&& &&&& &&&& &\\
};
\draw[solid,lightgray](M-1-1.south west)--(M-1-18.south east);
\draw[dotted,lightgray](M-2-1.south west)--(M-2-18.south east);
\draw[dotted,lightgray](M-3-1.south west)--(M-3-18.south east);
\draw[-latex,thick,gray](M-4-1.south west)--(M-4-18.south east);
\node (j) at (M-5-19) {$j$};
\draw[dotted,lightgray](M-5-1.south west)--(M-5-18.south east);
\draw[dotted,lightgray](M-6-1.south west)--(M-6-18.south east);
\draw[solid,lightgray](M-7-1.south west)--(M-7-18.south east);
\draw[dotted,lightgray](M-8-1.south west)--(M-8-18.south east);
\draw[dotted,lightgray](M-9-1.south west)--(M-9-18.south east);
\draw[solid,lightgray](M-10-1.south west)--(M-10-18.south east);
\draw[solid,lightgray](M-1-1.north east)--(M-11-1.south east);
\draw[dotted,lightgray](M-1-2.north east)--(M-11-2.south east);
\draw[dotted,lightgray](M-1-3.north east)--(M-11-3.south east);
\draw[dotted,lightgray](M-1-4.north east)--(M-11-4.south east);
\draw[-latex,thick,gray](M-1-5.north east)--(M-11-5.south east);
\node (i) at (M-12-6) {$i$};
\draw[dotted,lightgray](M-1-6.north east)--(M-11-6.south east);
\draw[dotted,lightgray](M-1-7.north east)--(M-11-7.south east);
\draw[dotted,lightgray](M-1-8.north east)--(M-11-8.south east);
\draw[solid,lightgray](M-1-9.north east)--(M-11-9.south east);
\draw[dotted,lightgray](M-1-10.north east)--(M-11-10.south east);
\draw[dotted,lightgray](M-1-11.north east)--(M-11-11.south east);
\draw[dotted,lightgray](M-1-12.north east)--(M-11-12.south east);
\draw[solid,lightgray](M-1-13.north east)--(M-11-13.south east);
\draw[dotted,lightgray](M-1-14.north east)--(M-11-14.south east);
\draw[dotted,lightgray](M-1-15.north east)--(M-11-15.south east);
\draw[dotted,lightgray](M-1-16.north east)--(M-11-16.south east);
\draw[solid,lightgray](M-1-17.north east)--(M-11-17.south east);
\end{tikzpicture}
\end{center}
Then $A^{\tt st'}$ is  
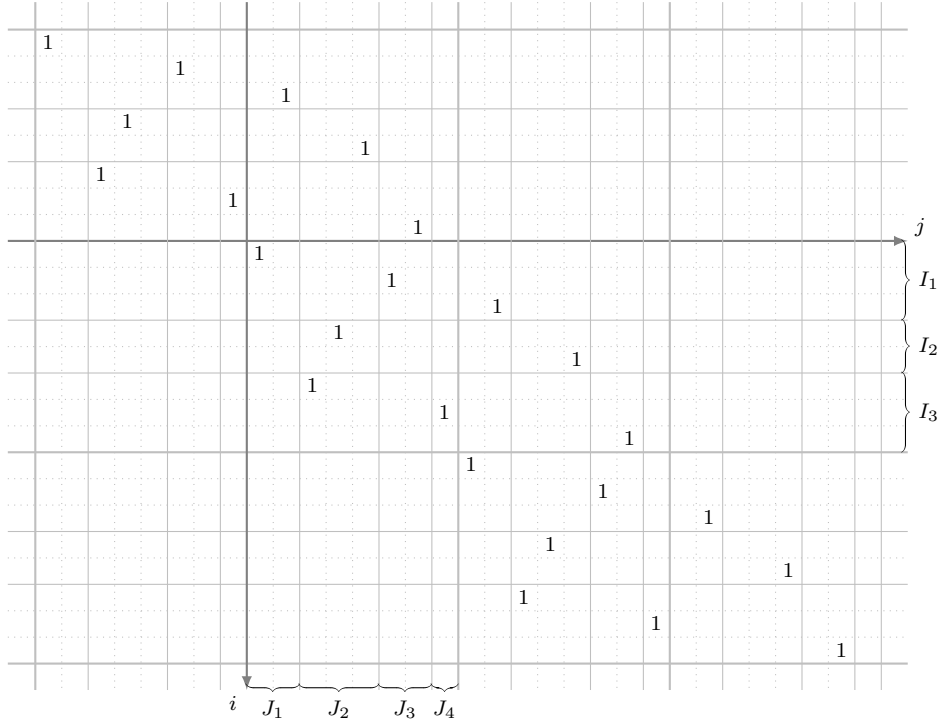
\begin{center}
\begin{tikzpicture}[every node/.style={font=\footnotesize,scale=1},BC/.style = {decorate, decoration={calligraphic brace, amplitude=3pt, raise=1mm},
        pen colour={black}
            }]
\matrix (M)[matrix of math nodes,nodes in empty cells,nodes={rectangle,minimum height=1.0em,minimum width=1.0em,inner sep=0pt,anchor=center,align=center}]
{
& &&&&&&&& &&&&&&&& &&&&&&&& &&&&&&&& &\\
& 1&&&&&&&& &&&&&&&& &&&&&&&& &&&&&&&& &\\
& &&&&&1&&& &&&&&&&& &&&&&&&& &&&&&&&& &\\
& &&&&&&&& &1&&&&&&& &&&&&&&& &&&&&&&& &\\
& &&&1&&&&& &&&&&&&& &&&&&&&& &&&&&&&& &\\
& &&&&&&&& &&&&1&&&& &&&&&&&& &&&&&&&& &\\
& &&1&&&&&& &&&&&&&& &&&&&&&& &&&&&&&& &\\
& &&&&&&&1& &&&&&&&& &&&&&&&& &&&&&&&& &\\
& &&&&&&&& &&&&&&1&& &&&&&&&& &&&&&&&& &\\
& &&&&&&&& 1&&&&&&&& &&&&&&&& &&&&&&&& &\\
& &&&&&&&& &&&&&1&&& &&&&&&&& &&&&&&&& &\\
& &&&&&&&& &&&&&&&& &1&&&&&&& &&&&&&&& &\\
& &&&&&&&& &&&1&&&&& &&&&&&&& &&&&&&&& &\\
& &&&&&&&& &&&&&&&& &&&&1&&&& &&&&&&&& &\\
& &&&&&&&& &&1&&&&&& &&&&&&&& &&&&&&&& &\\
& &&&&&&&& &&&&&&&1& &&&&&&&& &&&&&&&& &\\
& &&&&&&&& &&&&&&&& &&&&&&1&& &&&&&&&& &\\
& &&&&&&&& &&&&&&&& 1&&&&&&&& &&&&&&&& &\\
& &&&&&&&& &&&&&&&& &&&&&1&&& &&&&&&&& &\\
& &&&&&&&& &&&&&&&& &&&&&&&& &1&&&&&&& &\\
& &&&&&&&& &&&&&&&& &&&1&&&&& &&&&&&&& &\\
& &&&&&&&& &&&&&&&& &&&&&&&& &&&&1&&&& &\\
& &&&&&&&& &&&&&&&& &&1&&&&&& &&&&&&&& &\\
& &&&&&&&& &&&&&&&& &&&&&&&1& &&&&&&&& &\\
& &&&&&&&& &&&&&&&& &&&&&&&& &&&&&&1&& &\\
& &&&&&&&& &&&&&&&& &&&&&&&& &&&&&&&& &\\
& &&&&&&&& &&&&&&&& &&&&&&&& &&&&&&&& &\\
};
\draw[thick,lightgray](M-1-1.south west)--(M-1-34.south east);
\draw[dotted,lightgray](M-2-1.south west)--(M-2-34.south east);
\draw[dotted,lightgray](M-3-1.south west)--(M-3-34.south east);
\draw[solid,lightgray](M-4-1.south west)--(M-4-34.south east);
\draw[dotted,lightgray](M-5-1.south west)--(M-5-34.south east);
\draw[solid,lightgray](M-6-1.south west)--(M-6-34.south east);
\draw[dotted,lightgray](M-7-1.south west)--(M-7-34.south east);
\draw[dotted,lightgray](M-8-1.south west)--(M-8-34.south east);
\draw[-latex,thick,gray](M-9-1.south west)--(M-9-34.south east);
\node (j) at (M-9-35) {$j$};
\draw[dotted,lightgray](M-10-1.south west)--(M-10-34.south east);
\draw[dotted,lightgray](M-11-1.south west)--(M-11-34.south east);
\draw[solid,lightgray](M-12-1.south west)--(M-12-34.south east);
\draw[dotted,lightgray](M-13-1.south west)--(M-13-34.south east);
\draw[solid,lightgray](M-14-1.south west)--(M-14-34.south east);
\draw[dotted,lightgray](M-15-1.south west)--(M-15-34.south east);
\draw[dotted,lightgray](M-16-1.south west)--(M-16-34.south east);
\draw[thick,lightgray](M-17-1.south west)--(M-17-34.south east);
\draw[dotted,lightgray](M-18-1.south west)--(M-18-34.south east);
\draw[dotted,lightgray](M-19-1.south west)--(M-19-34.south east);
\draw[solid,lightgray](M-20-1.south west)--(M-20-34.south east);
\draw[dotted,lightgray](M-21-1.south west)--(M-21-34.south east);
\draw[solid,lightgray](M-22-1.south west)--(M-22-34.south east);
\draw[dotted,lightgray](M-23-1.south west)--(M-23-34.south east);
\draw[dotted,lightgray](M-24-1.south west)--(M-24-34.south east);
\draw[thick,lightgray](M-25-1.south west)--(M-25-34.south east);
\draw[thick,lightgray](M-1-1.north east)--(M-26-1.south east);
\draw[dotted,lightgray](M-1-2.north east)--(M-26-2.south east);
\draw[solid,lightgray](M-1-3.north east)--(M-26-3.south east);
\draw[dotted,lightgray](M-1-4.north east)--(M-26-4.south east);
\draw[dotted,lightgray](M-1-5.north east)--(M-26-5.south east);
\draw[solid,lightgray](M-1-6.north east)--(M-26-6.south east);
\draw[dotted,lightgray](M-1-7.north east)--(M-26-7.south east);
\draw[solid,lightgray](M-1-8.north east)--(M-26-8.south east);
\draw[-latex,thick,gray](M-1-9.north east)--(M-26-9.south east);
\draw[dotted,lightgray](M-1-10.north east)--(M-26-10.south east);
\node (i) at (M-27-9) {$i$};
\draw[solid,lightgray](M-1-11.north east)--(M-26-11.south east);
\draw[dotted,lightgray](M-1-12.north east)--(M-26-12.south east);
\draw[dotted,lightgray](M-1-13.north east)--(M-26-13.south east);
\draw[solid,lightgray](M-1-14.north east)--(M-26-14.south east);
\draw[dotted,lightgray](M-1-15.north east)--(M-26-15.south east);
\draw[solid,lightgray](M-1-16.north east)--(M-26-16.south east);
\draw[thick,lightgray](M-1-17.north east)--(M-26-17.south east);
\draw[dotted,lightgray](M-1-18.north east)--(M-26-18.south east);
\draw[solid,lightgray](M-1-19.north east)--(M-26-19.south east);
\draw[dotted,lightgray](M-1-20.north east)--(M-26-20.south east);
\draw[dotted,lightgray](M-1-21.north east)--(M-26-21.south east);
\draw[solid,lightgray](M-1-22.north east)--(M-26-22.south east);
\draw[dotted,lightgray](M-1-23.north east)--(M-26-23.south east);
\draw[solid,lightgray](M-1-24.north east)--(M-26-24.south east);
\draw[thick,lightgray](M-1-25.north east)--(M-26-25.south east);
\draw[dotted,lightgray](M-1-26.north east)--(M-26-26.south east);
\draw[solid,lightgray](M-1-27.north east)--(M-26-27.south east);
\draw[dotted,lightgray](M-1-28.north east)--(M-26-28.south east);
\draw[dotted,lightgray](M-1-29.north east)--(M-26-29.south east);
\draw[solid,lightgray](M-1-30.north east)--(M-26-30.south east);
\draw[dotted,lightgray](M-1-31.north east)--(M-26-31.south east);
\draw[solid,lightgray](M-1-32.north east)--(M-26-32.south east);
\draw[solid,lightgray](M-1-33.north east)--(M-26-33.south east);
\draw[BC] (M-10-34.north) -- node[right=2mm] {$I_1$} (M-12-34.south);
\draw[BC] (M-13-34.north) -- node[right=2mm] {$I_2$} (M-14-34.south);
\draw[BC] (M-15-34.north) -- node[right=2mm] {$I_3$} (M-17-34.south);
\draw[BC] (M-26-11.east) -- node[below=2mm] {$J_1$} (M-26-10.west);
\draw[BC] (M-26-14.east) -- node[below=2mm] {$J_2$} (M-26-12.west);
\draw[BC] (M-26-16.east) -- node[below=2mm] {$J_3$} (M-26-15.west);
\draw[BC] (M-26-17.east) -- node[below=2mm] {$J_4$} (M-26-17.west);

\end{tikzpicture}
\end{center}
}
\end{ex}

We define partial orders $>_{\tt nW}$ and $<_{\tt Ne}$ on $\Z^2$ as follows: 
\begin{itemize}
    \item[(1)] $c_1 >_{\tt nW} c_2$ if and only if $i_1 \le i_2$ and $j_1 < j_2$,
    \item[(2)] $c_1 <_{\tt Ne} c_2$ if and only if $i_1 > i_2$ and $j_1 \le j_2$,
\end{itemize}
for $c_1=(i_1, j_1)$ and $c_2=(i_2, j_2)$ in $\Z^2$. 

With respect to these partial orders, we have natural dual analogues of the notions and their properties in Section \ref{sec:prel}, the proofs of which are almost parallel to those in the case of $\wh{\mc M}_{m\times n}$.
Let us summarize them as follows: suppose that $A\in \wh{\mc{N}}_{m\times n}$ is given.

\begin{itemize}
 \item[$\bullet$] A proper numbering $d$ on $A$ is defined as in Definition \ref{def:proper} with respect to $>_{\tt nW}$ instead of $>_{\tt NW}$. Let $d^{\tt st'}$ denote the proper numbering on $A^{\tt st'}$ which {corresponds} to $d$. Lemmas \ref{lem:d and dst} and \ref{lem:period is unique} hold for proper numberings on $A$.
 
 \item[$\bullet$] A stream ${\bf s}$ and a channel $C$ of $A$ are defined in the same way as in Definition  \ref{def:channel} with respect to $>_{\tt nW}$. Let ${\bf s}^{\tt st'}$ and $C^{\tt st'}$ denote the stream and channel of $A^{\tt st'}$ which naturally correspond to ${\bf s}$ and $C$, respectively.
 
 \item[$\bullet$] Let $\mathcal{C}'_A$ be the set of channels of $A$. 
 We define a partial order $\succcurlyeq_{\tt sw}$ on $\mathcal{C}'_A$ as in \eqref{eq:SW order on channels} by using $\le_{\tt Ne}$ instead of $\le_{\tt ne}$ (see also \eqref{eq:C_ne C_sw}). 
Then $\mathcal{C}'_A$ has the greatest and the smallest elements with respect to $\succcurlyeq_{\tt sw}$. It can be proved by the same arguments as in Proposition \ref{prop:lattice property} using $A^{\tt st'}$.
 
 \item[$\bullet$] Let $C^{\tt sw'}_A$ be the greatest element in $\mathcal{C}'_A$, which we call the southwest channel of $A$. We have 
 $(C_A^{\tt sw'})^{\tt st'}=C_{A^{\tt st'}}^{\tt sw'}$.
 
 \item[$\bullet$]  For $C\in \mc{C}'_A$, we define a channel numbering $d^C_A$ as in \eqref{eq:channel numbering} with respect to $>_{\tt nW}$. Then $d^C_A$ is a well-defined proper numbering on $A$ (cf.~Proposition \ref{prop:ch numbering}), and
\begin{equation*}\label{eq:sw proper A and Ast'}
\left(d^C_A\right)^{\tt st'}=d^{C^{\tt st'}}_{A^{\tt st'}}.
\end{equation*}

 \item[$\bullet$] Let $d^{\tt sw'}_A=d^{C^{\tt sw'}_A}_A$ be the southwest channel numbering on $A$. Lemma \ref{lem:characterization of ch numbering} also holds.

\end{itemize}

\subsection{Row semistandard tableaux and offset vectors}

For $a \ge 1$, let $\mc{B}'_n((a))=\mc{B}'((a))$ be the set of $T\in SST_\Z((a))$ such that $T(a)-T(1) \le n$, where $T(j)$ is the $j$-th entry of $T$ from the left for $1 \le j \le a$.
Let $\tau_n=\tau$ be a bijection on $\mc{B}'((a))$ where $\tau(T)$ is the row semistandard tableau obtained from $T\in \mc{B}'((a))$ by replacing its entries $T(1)\le T(2)\le \dots\le T(a)$ with $T(2)\le T(3)\le \dots \le T(a)\le T(1)+n$.

Let $b\ge 1$. For $\alpha=(\alpha_1,\dots,\alpha_b)\in \Z^b$ and $(T_1,\dots,T_b)\in \mc{B}'((a))\times \dots \times \mc{B}'((a))$, we define
\begin{equation*}
\tau^{\alpha}(T_1,\dots,T_b) = \left(\tau^{\alpha_1}(T_1),\dots,\tau^{\alpha_a}(T_b)\right).
\end{equation*}

Let $R=(a^b)$.
Let us regard $T\in RSST_{[n]}(R)$ as $T=(T^b,\dots,T^1)\in \mc{B}'((a))\times \dots \times \mc{B}'((a))$, where $T^i$ is the $i$-th row of $T$ from the bottom. For $1\le i\le b-1$, let $r_i$ be the minimal non-negative integer such that 
$$
(T^{i+1},\tau^{r_i}(T^i))
$$ is a $\Z$-semistandard tableau of shape $(a^2)$,
and put $\eta_i = r_i + r_{i+1} +\dots + r_{b-1}$. We call $r=(r_1, \dots, r_{b-1})$ an {\em offset vector} and $\eta = (\eta_1, \dots, \eta_{b-1})$ a {\em symmetrized offset vector} of $T$. Note that $\tau^{\eta_{\tt rev}}(T) \in SST_\Z(R)$.

\subsection{Dual affine RSK correspondence}\label{subsec:dual aff RSK}

Let $d$ be a proper numbering on $A\in \wh{\mc N}_{m\times n}$. 
Note that each level set $d^{-1}(k)$ forms a chain with respect to $<_{\tt Ne}$. 

Let $\{{\bf z}_k\}_{k \in \Z}$ be the set of zig-zags associated to $d$, where the inner corners of ${\bf z}_k$ are the set of elements in $d^{-1}(k)$  maximal with respect to $>_{\tt nW}$ (cf.~Section \ref{subsec:mbc}). 
Then $\{ {\bf z}_k \}_{k \in \Z}$ satisfies
\begin{itemize}
    \item[(${\rm z'}$.1)] the inner corners of each ${\bf z}_k$ are contained in ${\rm supp}(A)$,
    \item[(${\rm z'}$.2)] ${\rm supp}(A) \subseteq \bigcup_{k \in \Z} {\bf z}_k$, 
    \item[(${\rm z'}$.3)] ${\bf z}_k$ is located to the southeast of ${\bf z}_{k-1}$ for $k\in\Z$ in the sense of \eqref{eq:NW chain on zig-zag} with respect to \!\! $>_{\tt nW}$.
\end{itemize}

\begin{rem}\label{rem:zig-zag for proper numbering dual}
{\rm
We should remark that no outer cell of ${\bf z}_k$ belongs to ${\rm supp}(A)$, and ${\bf z}_k$'s are not always mutually disjoint. More precisely, two horizontal lines (or line segments) in ${\bf z}_k$ and ${\bf z}_l$ ($k<l$) may have non-trivial intersection, while vertical lines (or line segmenents) in ${\bf z}_k$ and ${\bf z}_l$ ($k<l$) are always disjoint.
} 
\end{rem}

Suppose that a non-zero $A\in \wh{\mc N}_{m\times n}$ is given and let
\begin{itemize}

\item[$\bullet$] $\{{\bf z}_k\}_{k\in\Z}$ : the set of zig-zags associated to $d_A^{\tt sw'}$,

\item[$\bullet$] $A^{\flat'}$ : the matrix obtained from $A$ by the same rule as in $A^\flat$ with respect to $\{{\bf z}_k\}_{k\in\Z}$, 

\item[$\bullet$] $A^{(t)}$ : the matrices in $\wh{\mc N}_{m\times n}$ defined inductively by
\begin{equation*}
A^{(0)}=A,\quad A^{(t)}=\left( A^{(t-1)} \right)^{\flat'} \quad (t\ge 1).
\end{equation*}
\end{itemize}
Note that $A^{(s-1)}\neq {\mathbb O}$ and $A^{(s)} ={\mathbb O}$ for some $s\ge 1$.
For $1\le t\le s$, we let
\begin{itemize}
\item[$\bullet$] $\{{\bf z}^{(t)}_k\}_{k\in \Z}$ : the set of zig-zags associated to $d^{\tt sw'}_{A^{(t-1)}}$,

\item[$\bullet$] ${\bf s}^{(t)}=(\ba_t,\bb_t,\rho_t)$ : the stream of the back-post corners of $\{{\bf z}^{(t)}_k\}_{k\in \Z}$ with flow $\mu_t$,

\end{itemize}
where we can check that $\mu=(\mu_1,\dots,\mu_s)\in \cP_{s}$ as in Lemma \ref{lem:width decreasing}.
Now let
\begin{itemize}
\item[$\bullet$] $P_0$, $Q_0$ : the tableau of shape $\la=\mu'$ defined as in Section \ref{subsec:mbc},

\item[$\bullet$] $P_0^t$: the tableau of shape $\mu$ obtained by flipping $P_0$ with respect to the main diagonal,

\item[$\bullet$] $\rho=(\rho_1,\dots,\rho_s)\in \Z^s$.
 
\end{itemize}
\begin{lem}\label{lem:P_0 Q-0 for dual}
Under the above hypothesis, $P_0^t$ is row semistandard, while $Q_0$ is column semistandard. 
\end{lem}
\pf It follows immediately from Remark \ref{rem:zig-zag for proper numbering dual}.
\qed\newline

Let $l = \ell(\la)$ and let $m_i$ be the number of occurrences of $i$ in $\la'$ for $1\le i \le l$.
Let
\begin{itemize}
\item[$\bullet$] $R_i$ : the Young diagram given in \eqref{eq:rec decomp for lambda},

\item[$\bullet$] $P^{(i)}_0$, $Q^{(i)}_0$ : the subtableaux of $P_0$ and $Q_0$ corresponding to $R_i$, respectively,

\item[$\bullet$] $\rho^{(i)}\in \Z^{m_i}$ : the subsequence of $\rho$ corresponding to the columns of $R_i$, 

\item[$\bullet$] $\eta^{(i)}\in \cP_{m_i-1}$ : the symmetrized offset vector of $\left(P^{(i)}_0\right)^t$.

\end{itemize}
for $1\le i\le l$ with $m_i \ge 1$, and
\begin{itemize}
\item[$\bullet$] $Q$ : the tableau defined as in \eqref{eq:(P,Q)}.
\end{itemize}
Then we have the following, which is a dual analogue of Theorem \ref{thm:main-1}.
\begin{thm}\label{thm:main-1'}
We have a bijection
\begin{equation*}
\xymatrixcolsep{3pc}\xymatrixrowsep{0.5pc}\xymatrix{
\kappa':\ \ \wh{\mc N}_{m\times n} \ \ar@{->}[r]  & \ \
\displaystyle{\bigsqcup_{\la\in \cP_n} RSST_{[m]}(\la') \times \mc{B}_n(\la)} \\
\ \ A \  \ar@{|->}[r]  &\ \ (P_0^t,Q)},
\end{equation*}
where $\kappa'({\mathbb O}) = (\emptyset, \emptyset)$.
\end{thm}

\begin{ex}{\rm
Let $A$ be the generalized permutation in Example \ref{ex:Adual}. 
Then \vskip 2mm

\begin{center}
\begin{tikzpicture}[every node/.style={font=\footnotesize,scale=0.75}]
\matrix (M)[matrix of math nodes,nodes in empty cells,nodes={rectangle,minimum height=1.2em,minimum width=1.2em,inner sep=0pt,anchor=center,align=center}]
{
& &&&& &&&& &&&& &&&& &\\
& 1&&1&& 1&&&& &&&& &&&& &\\
& &1&&& &1&&& &&&& &&&& &\\
& &1&&1& &&1&& &&&& &&&& &\\
& &&&& 1&&1&& 1&&&& &&&& &\\
& &&&& &1&&& &1&&& &&&& &\\
& &&&& &1&&1& &&1&& &&&& &\\
& &&&& &&&& 1&&1&& 1&&&& &\\
& &&&& &&&& &1&&& &1&&& &\\
& &&&& &&&& &1&&1& &&1&& &\\
& &&&& &&&& &&&& &&&& &\\
& &&&& &&&& &&&& &&&& &\\
};
\draw[solid,lightgray](M-1-1.south west)--(M-1-18.south east);
\draw[dotted,ultra thin](M-2-1.south west)--(M-2-18.south east);
\draw[dotted,ultra thin](M-3-1.south west)--(M-3-18.south east);
\draw[-latex,thick,gray](M-4-1.south west)--(M-4-18.south east);
\node (j) at (M-5-19) {$j$};
\draw[dotted,ultra thin](M-5-1.south west)--(M-5-18.south east);
\draw[dotted,ultra thin](M-6-1.south west)--(M-6-18.south east);
\draw[solid,lightgray](M-7-1.south west)--(M-7-18.south east);
\draw[dotted,ultra thin](M-8-1.south west)--(M-8-18.south east);
\draw[dotted,ultra thin](M-9-1.south west)--(M-9-18.south east);
\draw[solid,lightgray](M-10-1.south west)--(M-10-18.south east);
\draw[solid,lightgray](M-1-1.north east)--(M-11-1.south east);
\draw[dotted,ultra thin](M-1-2.north east)--(M-11-2.south east);
\draw[dotted,ultra thin](M-1-3.north east)--(M-11-3.south east);
\draw[dotted,ultra thin](M-1-4.north east)--(M-11-4.south east);
\draw[-latex,thick,gray](M-1-5.north east)--(M-11-5.south east);
\node (i) at (M-12-6) {$i$};
\draw[dotted,ultra thin](M-1-6.north east)--(M-11-6.south east);
\draw[dotted,ultra thin](M-1-7.north east)--(M-11-7.south east);
\draw[dotted,ultra thin](M-1-8.north east)--(M-11-8.south east);
\draw[solid,lightgray](M-1-9.north east)--(M-11-9.south east);
\draw[dotted,ultra thin](M-1-10.north east)--(M-11-10.south east);
\draw[dotted,ultra thin](M-1-11.north east)--(M-11-11.south east);
\draw[dotted,ultra thin](M-1-12.north east)--(M-11-12.south east);
\draw[solid,lightgray](M-1-13.north east)--(M-11-13.south east);
\draw[dotted,ultra thin](M-1-14.north east)--(M-11-14.south east);
\draw[dotted,ultra thin](M-1-15.north east)--(M-11-15.south east);
\draw[dotted,ultra thin](M-1-16.north east)--(M-11-16.south east);
\draw[solid,lightgray](M-1-17.north east)--(M-11-17.south east);
\begin{scope}[on background layer]
\draw[thin,red](M-11-9.south) -- (M-7-9.center) -| (M-5-10.center)  -- (M-5-18.east);
\draw[thin,red](M-11-10.south) -- (M-8-10.center) -| (M-6-11.center)  -- (M-6-18.east);
\draw[thin,red](M-11-11.south) -- (M-9-11.center) -| (M-7-12.center)  -- (M-7-18.east);
\end{scope}
\node[anchor=east,scale=1.5] at (M.west) {$A^{(0)}=$};
\end{tikzpicture}\quad
\begin{tikzpicture}[every node/.style={font=\footnotesize,scale=0.75}]
\matrix (M)[matrix of math nodes,nodes in empty cells,nodes={rectangle,minimum height=1.2em,minimum width=1.2em,inner sep=0pt,anchor=center,align=center}]
{
& &&&& &&&& &&&& &&&& &\\
& &1&1&& &&&& &&&& &&&& &\\
& &&1&& &&&& &&&& &&&& &\\
& &1&&& 1&&&& &&&& &&&& &\\
& &&&& &1&1&& &&&& &&&& &\\
& &&&& &&1&& &&&& &&&& &\\
& &&&& &1&&& 1&&&& &&&& &\\
& &&&& &&&& &1&1&& &&&& &\\
& &&&& &&&& &&1&& &&&& &\\
& &&&& &&&& &1&&& 1&&&& &\\
& &&&& &&&& &&&& &&&& &\\
& &&&& &&&& &&&& &&&& &\\
};
\draw[solid,lightgray](M-1-1.south west)--(M-1-18.south east);
\draw[dotted,ultra thin](M-2-1.south west)--(M-2-18.south east);
\draw[dotted,ultra thin](M-3-1.south west)--(M-3-18.south east);
\draw[-latex,thick,gray](M-4-1.south west)--(M-4-18.south east);
\node (j) at (M-5-19) {$j$};
\draw[dotted,ultra thin](M-5-1.south west)--(M-5-18.south east);
\draw[dotted,ultra thin](M-6-1.south west)--(M-6-18.south east);
\draw[solid,lightgray](M-7-1.south west)--(M-7-18.south east);
\draw[dotted,ultra thin](M-8-1.south west)--(M-8-18.south east);
\draw[dotted,ultra thin](M-9-1.south west)--(M-9-18.south east);
\draw[solid,lightgray](M-10-1.south west)--(M-10-18.south east);
\draw[solid,lightgray](M-1-1.north east)--(M-11-1.south east);
\draw[dotted,ultra thin](M-1-2.north east)--(M-11-2.south east);
\draw[dotted,ultra thin](M-1-3.north east)--(M-11-3.south east);
\draw[dotted,ultra thin](M-1-4.north east)--(M-11-4.south east);
\draw[-latex,thick,gray](M-1-5.north east)--(M-11-5.south east);
\node (i) at (M-12-6) {$i$};
\draw[dotted,ultra thin](M-1-6.north east)--(M-11-6.south east);
\draw[dotted,ultra thin](M-1-7.north east)--(M-11-7.south east);
\draw[dotted,ultra thin](M-1-8.north east)--(M-11-8.south east);
\draw[solid,lightgray](M-1-9.north east)--(M-11-9.south east);
\draw[dotted,ultra thin](M-1-10.north east)--(M-11-10.south east);
\draw[dotted,ultra thin](M-1-11.north east)--(M-11-11.south east);
\draw[dotted,ultra thin](M-1-12.north east)--(M-11-12.south east);
\draw[solid,lightgray](M-1-13.north east)--(M-11-13.south east);
\draw[dotted,ultra thin](M-1-14.north east)--(M-11-14.south east);
\draw[dotted,ultra thin](M-1-15.north east)--(M-11-15.south east);
\draw[dotted,ultra thin](M-1-16.north east)--(M-11-16.south east);
\draw[solid,lightgray](M-1-17.north east)--(M-11-17.south east);
\begin{scope}[on background layer]
\draw[thin,red](M-11-7.south) -- (M-5-7.center)  -- (M-5-18.east);
\draw[thin,red](M-11-8.south) -- (M-5-8.center)  -- (M-5-18.east);
\draw[thin,red](M-11-10.south) -- (M-7-10.center)  -- (M-7-18.east);
\end{scope}
\node[anchor=east,scale=1.5] at (M.west) {$A^{(1)}=$};
\end{tikzpicture} 
\end{center}

\begin{center}
\begin{tikzpicture}[every node/.style={font=\footnotesize,scale=0.75}]
\matrix (M)[matrix of math nodes,nodes in empty cells,nodes={rectangle,minimum height=1.2em,minimum width=1.2em,inner sep=0pt,anchor=center,align=center}]
{
& &&&& &&&& &&&& &&&& &\\
& &&&& &&&& &&&& &&&& &\\
& &&1&& &&&& &&&& &&&& &\\
& &1&&& &&&& &&&& &&&& &\\
& &&&& &&&& &&&& &&&& &\\
& &&&& &&1&& &&&& &&&& &\\
& &&&& &1&&& &&&& &&&& &\\
& &&&& &&&& &&&& &&&& &\\
& &&&& &&&& &&1&& &&&& &\\
& &&&& &&&& &1&&& &&&& &\\
& &&&& &&&& &&&& &&&& &\\
& &&&& &&&& &&&& &&&& &\\
};
\draw[solid,lightgray](M-1-1.south west)--(M-1-18.south east);
\draw[dotted,ultra thin](M-2-1.south west)--(M-2-18.south east);
\draw[dotted,ultra thin](M-3-1.south west)--(M-3-18.south east);
\draw[-latex,thick,gray](M-4-1.south west)--(M-4-18.south east);
\node (j) at (M-5-19) {$j$};
\draw[dotted,ultra thin](M-5-1.south west)--(M-5-18.south east);
\draw[dotted,ultra thin](M-6-1.south west)--(M-6-18.south east);
\draw[solid,lightgray](M-7-1.south west)--(M-7-18.south east);
\draw[dotted,ultra thin](M-8-1.south west)--(M-8-18.south east);
\draw[dotted,ultra thin](M-9-1.south west)--(M-9-18.south east);
\draw[solid,lightgray](M-10-1.south west)--(M-10-18.south east);
\draw[solid,lightgray](M-1-1.north east)--(M-11-1.south east);
\draw[dotted,ultra thin](M-1-2.north east)--(M-11-2.south east);
\draw[dotted,ultra thin](M-1-3.north east)--(M-11-3.south east);
\draw[dotted,ultra thin](M-1-4.north east)--(M-11-4.south east);
\draw[-latex,thick,gray](M-1-5.north east)--(M-11-5.south east);
\node (i) at (M-12-6) {$i$};
\draw[dotted,ultra thin](M-1-6.north east)--(M-11-6.south east);
\draw[dotted,ultra thin](M-1-7.north east)--(M-11-7.south east);
\draw[dotted,ultra thin](M-1-8.north east)--(M-11-8.south east);
\draw[solid,lightgray](M-1-9.north east)--(M-11-9.south east);
\draw[dotted,ultra thin](M-1-10.north east)--(M-11-10.south east);
\draw[dotted,ultra thin](M-1-11.north east)--(M-11-11.south east);
\draw[dotted,ultra thin](M-1-12.north east)--(M-11-12.south east);
\draw[solid,lightgray](M-1-13.north east)--(M-11-13.south east);
\draw[dotted,ultra thin](M-1-14.north east)--(M-11-14.south east);
\draw[dotted,ultra thin](M-1-15.north east)--(M-11-15.south east);
\draw[dotted,ultra thin](M-1-16.north east)--(M-11-16.south east);
\draw[solid,lightgray](M-1-17.north east)--(M-11-17.south east);
\begin{scope}[on background layer]
\draw[thin,red](M-11-7.south) -- (M-7-7.center) -| (M-6-8.center)  -- (M-6-18.east);
\end{scope}
\node[anchor=east,scale=1.5] at (M.west) {$A^{(2)}=$};
\end{tikzpicture}
\quad
\begin{tikzpicture}[every node/.style={font=\footnotesize,scale=0.75}]
\matrix (M)[matrix of math nodes,nodes in empty cells,nodes={rectangle,minimum height=1.2em,minimum width=1.2em,inner sep=0pt,anchor=center,align=center}]
{
& &&&& &&&& &&&& &&&& &\\
& &&&& &&&& &&&& &&&& &\\
& &&&& &&&& &&&& &&&& &\\
& &&1&& &&&& &&&& &&&& &\\
& &&&& &&&& &&&& &&&& &\\
& &&&& &&&& &&&& &&&& &\\
& &&&& &&1&& &&&& &&&& &\\
& &&&& &&&& &&&& &&&& &\\
& &&&& &&&& &&&& &&&& &\\
& &&&& &&&& &&1&& &&&& &\\
& &&&& &&&& &&&& &&&& &\\
& &&&& &&&& &&&& &&&& &\\
};
\draw[solid,lightgray](M-1-1.south west)--(M-1-18.south east);
\draw[dotted,ultra thin](M-2-1.south west)--(M-2-18.south east);
\draw[dotted,ultra thin](M-3-1.south west)--(M-3-18.south east);
\draw[-latex,thick,gray](M-4-1.south west)--(M-4-18.south east);
\node (j) at (M-5-19) {$j$};
\draw[dotted,ultra thin](M-5-1.south west)--(M-5-18.south east);
\draw[dotted,ultra thin](M-6-1.south west)--(M-6-18.south east);
\draw[solid,lightgray](M-7-1.south west)--(M-7-18.south east);
\draw[dotted,ultra thin](M-8-1.south west)--(M-8-18.south east);
\draw[dotted,ultra thin](M-9-1.south west)--(M-9-18.south east);
\draw[solid,lightgray](M-10-1.south west)--(M-10-18.south east);
\draw[solid,lightgray](M-1-1.north east)--(M-11-1.south east);
\draw[dotted,ultra thin](M-1-2.north east)--(M-11-2.south east);
\draw[dotted,ultra thin](M-1-3.north east)--(M-11-3.south east);
\draw[dotted,ultra thin](M-1-4.north east)--(M-11-4.south east);
\draw[-latex,thick,gray](M-1-5.north east)--(M-11-5.south east);
\node (i) at (M-12-6) {$i$};
\draw[dotted,ultra thin](M-1-6.north east)--(M-11-6.south east);
\draw[dotted,ultra thin](M-1-7.north east)--(M-11-7.south east);
\draw[dotted,ultra thin](M-1-8.north east)--(M-11-8.south east);
\draw[solid,lightgray](M-1-9.north east)--(M-11-9.south east);
\draw[dotted,ultra thin](M-1-10.north east)--(M-11-10.south east);
\draw[dotted,ultra thin](M-1-11.north east)--(M-11-11.south east);
\draw[dotted,ultra thin](M-1-12.north east)--(M-11-12.south east);
\draw[solid,lightgray](M-1-13.north east)--(M-11-13.south east);
\draw[dotted,ultra thin](M-1-14.north east)--(M-11-14.south east);
\draw[dotted,ultra thin](M-1-15.north east)--(M-11-15.south east);
\draw[dotted,ultra thin](M-1-16.north east)--(M-11-16.south east);
\draw[solid,lightgray](M-1-17.north east)--(M-11-17.south east);
\begin{scope}[on background layer]
\draw[thin,red](M-11-8.south) -- (M-7-8.center)  -- (M-7-18.east);
\end{scope}
\node[anchor=east,scale=1.5] at (M.west) {$A^{(3)}=$};
\end{tikzpicture}
\end{center}
where the red lines denote the {zig-zags} associated to $d^{\tt sw'}_{A^{(t-1)}}$ for $1 \le t \le 4$, and 
\begin{center}
$P^t_0 = 
\begin{ytableau}
1&2&3\\
1&1&3\\
2\\
3
\end{ytableau}
$\,,\quad $Q_0 = 
\begin{ytableau}
1&1&2&3\\
2&2\\
4&3
\end{ytableau}$\,,\quad $\rho = (2, 1, 0, 0)$.
\end{center}
\vskip 2mm

In this case, $R_1$ and $R_3$ are the only non-trivial rectangles in the decomposition of the shape of $P_0$ and $Q_0$.
It is easy to see that $\eta^{(3)}_{\tt rev} = (0, 2)$, $\eta^{(1)}_{\tt rev} = (0, 0)$, and hence\vskip 2mm

$$Q = 
\begin{ytableau}
4&5&2&3\\
5&6\\
6&7 
\end{ytableau}\,.$$\vskip 2mm
}
\end{ex}

\subsection{Proof of Theorem \ref{thm:main-1'}}
The proof of Theorem \ref{thm:main-1'} is almost parallel to that of Theorem \ref{thm:main-1} in Section \ref{sec:Proof of affine RSK}. 
So we give an outline of its proof and leave the details to the reader. 

Let $\alpha=(\alpha_1,\dots,\alpha_m)\in \Z_{\geq 0}^m$ and $\beta=(\beta_1,\dots,\beta_n)\in \Z_{\geq 0}^n$ be given. Let
\begin{equation*}
 \wh{\mc N}_{m \times n}(\alpha,\beta)=\{\,A\,|\,A\in \wh{\mc N}_{m\times n},\ {\rm row}(A)=\alpha, \ {\rm col}(A)=\beta\,\},
\end{equation*}
where ${\rm row}(A)$ and ${\rm col}(A)$ are given as in \eqref{eq:row and col of A}. 

Let $K\ge 1$ be given. 
We say that $i \in [K]$ is an {\em ascent} of $B\in \td{W}_K$ if $(i, j), (i+1, j') \in {\rm supp}(B)$ and $j < j'$. 
 we say that $B \in \td{W}_K$ is {\em $\alpha$-ascending} if for any $k \in {\rm supp}(\alpha)$ and $i$ with 
\begin{equation*}
\alpha_1 + \dots + \alpha_{k-1} < i < \alpha_1 + \dots + \alpha_{k-1} +\alpha_k,
\end{equation*}
$i$ is an ascent of $B$.
Let $\td{W}_{K,[\alpha,\beta]}$ denote the set of $B\in \td{W}_K$ such that $B$ is $\alpha$-ascending and $B^t$ is $\beta$-descending. 
As in Lemma \ref{lem:inverse std}, we have a bijection
\begin{equation*}
\xymatrixcolsep{3pc}\xymatrixrowsep{0.5pc}\xymatrix{
\wh{\mc N}_{m \times n }(\alpha,\beta) \ \ar@{->}[r]  & \ \
 \td{W}_{K,[\alpha,\beta]} \\
\ \ A \  \ar@{|->}[r]  &\ \ A^{\tt st'} }. 
\end{equation*}

Let $\la\in \cP$ be given with $|\la|=K$. 
Let $RSST_{[m]}(\la)_\alpha$ be the set of $T\in RSST_{[m]}(\la)$ with content $\alpha$.
For $T\in RSST_{[m]}(\la)_\alpha$, we define $T^{\tt st'}$ to be a tableau obtained from $T$ by replacing each $k \in [m]$ in $T$ with $\alpha_k\neq 0$ by the consecutive numbers
\begin{equation*}
\alpha_1 + \dots + \alpha_{k-1} + 1 < \dots < \alpha_1 + \dots + \alpha_{k-1} + \alpha_{k} 
\end{equation*}
from the bottom row to the top one, and from left to right in each row.
Then we have 
\begin{equation*}
 T^{\tt st'}\in RST_{[K]}(\la),
\end{equation*} 
We call $T^{\tt st'}$ the {\em standardization of $T$}.
Suppose that $R=(a^b) \in \cP$ for some $a, b \ge 1$.
As in Lemma \ref{lem:std and msstd}, we have 
\begin{equation}\label{eq:off vec equal}
\tau_m^{\nu_{\tt rev}}(T) \in SST_\Z(R) \quad\text{if and only if}\quad \tau_K^{\nu_{\tt rev}}(T^{\tt st'}) \in SST_\Z(R),
\end{equation} 
for $T \in RSST_{[m]}(R)$ and $\nu \in \cP_{a-1}$.
In particular, the symmetrized offset vector of $T$ is the same as that of $T^{\tt st'}$.

Let $S \in RST_{[K]}(\la)$ be given. 
We say that $i \in [K]$ is a {\em (row) ascent} of $S$ if the letter $i+1$ does not appear below $i$ in $T$.
We define $S \in RST_{[K]}(\la)$ to be {\em $\alpha$-ascending} in the same way as in the case of $CST_{[K]}(\la)$ with respect to ascent.
Let $RST_{[K],\alpha}(\la)$ be the set of $\alpha$-ascending tableaux in $RST_{[K]}(\la)$. As in Lemma \ref{lem:std of tableaux}, we have a bijection
\begin{equation*}
\xymatrixcolsep{3pc}\xymatrixrowsep{0.5pc}\xymatrix{
RSST_{[m]}(\la)_\alpha \ \ar@{->}[r]  & \ \
 RST_{[K],\alpha}(\la) \\
\ \ T \  \ar@{|->}[r]  &\ \ T^{\tt st'} }. 
\end{equation*}

Let $B \in \td{W}_K$ be given with $\kappa_0(B) = (P_0, Q_0, \rho)$, where $\kappa_0$ is the map in \eqref{eq:Aff RSK 0}.
Suppose that the shape of $P_0$ and $Q_0$ is $\la$ so that $P_0^t\in RST_{[K]}(\la')$, $Q_0 \in CST_{[K]}(\la)$. 
By Lemma \ref{lem:RSK and descent}, we have
as in Corollary \ref{cor:affine RS preserves descents} 
\begin{equation}\label{eq:aux-1}
 \text{$B\in \td{W}_{K,[\alpha,\beta]}$ if and only if $P^t_0\in RST_{[K],\alpha}(\la')$ and $Q_0 \in CST_{[K],\beta}(\la)$}.
\end{equation}

We also have the following analogue of Proposition \ref{prop:affine RSK and standardization}.
For $A \in \wh{\mc N}_{m \times n}$, let $\kappa'_0(A) = (P_0, Q_0, \rho)$ be given in Section \ref{subsec:dual aff RSK}. Then
\begin{equation}\label{eq:aux-2}
 \kappa_0(A^{\tt st'}) = (P_0^{\tt st'}, Q^{\tt st}_0, \rho),
\end{equation}
where we understand $P_0^{\tt st'}$ as $((P_0^t)^{\tt st'})^t$.

Let $(P_0, Q_0, \rho)$ be a triple such that $(P_0^t, Q_0, \rho) \in RSST_{[m]}(\la')\times CSST_{[n]}(\la)\times \Z^{\la_1}$ for some $\la\in \cP_n$.
Keeping the notations in Section \ref{subsec:dual aff RSK}, we let 
\begin{itemize}
\item[$\bullet$] $\theta^{(i)}\in \cP_{m_i-1}$ : the symmetrized offset vector of $Q^{(i)}_0$,

\item[$\bullet$] $\zeta^{(i)} = \theta^{(i)} - \eta^{(i)} \in \Z^{m_i}$ with the last component being zero,

\end{itemize}
for $1 \le i \le l$ with $m_i \ge 1$. 
We say that $(P_0, Q_0, \rho)$ is {\em dominant} if 
$\rho_{\tt rev} - \zeta \in \mc{P}(\la)$, where $\zeta=(\zeta^{(l)},\dots,\zeta^{(1)})  \in \Z^{m_l} \times \dots \times \Z^{m_1}$.
Then as in Lemma \ref{lem:std preserves dominance}, it follows from \eqref{eq:off vec equal} that 
\begin{equation}\label{eq:aux-3}
 \text{$(P_0, Q_0, \rho)$ is dominant if and only if $(P^{\tt st'}_0, Q^{\tt st}_0, \rho)$ is dominant.}
\end{equation}
Furthermore, by Remark \ref{rem:dominance CPY}, \eqref{eq:aux-3} is equivalent to saying that $\rho$ is dominant in the sense of \cite[Definition 5.8]{CPY}. 

Let $\Omega'_{\rm dom}$ be the set of $(P_0, Q_0, \rho)$ which are dominant. 
Then it follows from \eqref{eq:aux-1}, \eqref{eq:aux-2}, \eqref{eq:aux-3} and \cite[Theorem 5.1, Theorem 5.11]{CPY} that
we have a bijection 
\begin{equation*}
\xymatrixcolsep{2.5pc}\xymatrixrowsep{0.5pc}
\xymatrix{
\kappa'_0: \wh{\mc N}_{m\times n}  \ \ar@{->}[r]  & \ \ \Omega'_{\rm dom}
}.
\end{equation*} 
Finally, by the same arguments as in \ref{lem:kapp_0 to kappa}, we conclude that
\begin{equation*}
\xymatrixcolsep{3pc}\xymatrixrowsep{0.5pc}\xymatrix{
\kappa':\ \ \wh{\mc N}_{m\times n} \ \ar@{->}[r]  & \ \
\displaystyle{\bigsqcup_{\la\in \cP_n} RSST_{[m]}(\la') \times \mc{B}_n(\la)} \\
\ \ A \  \ar@{|->}[r]  &\ \ (P_0^t,Q)},
\end{equation*}
is a well-defined bijection. This proves Theorem \ref{thm:main-1'}.\qed

\subsection{Affine crystal morphisms}

We first remark that $\wh{\mc N}_{m \times n}$ has a structure of normal $A_{m-1}^{(1)}$-crystal with $P^0$ and $P^0_{\rm cl}$-weights with respect to ${\rm wt}, \te_i, \tf_i$ for $i\in \{\,0,1,\dots,m-1\,\}$ for $m \ge 2$, where ${\rm wt}={\rm wt}^0$ in \eqref{eq:wt A} and ${\rm wt}^0_{\rm cl}$ in \eqref{eq:wt_0 A} respectively. {For $A \in \wh{\mc N}_{m \times n}$, $\te_iA$ and $\tf_i A$ are given as follows:

\begin{itemize}
\item[(1)] Let
\begin{equation*}
\sigma = (\ \cdots,\ \underbrace{+}_{a_{i\, j+1}},\ \underbrace{-}_{a_{i+1\, j+1}},\ \underbrace{+}_{a_{i\, j}},\ \underbrace{-}_{a_{i+1\, j}},\ \cdots \ ), 
\end{equation*} and let $\td{\sigma}$ be the reduced one.
\item[(2)] If $\td{\sigma}$ has at least one $-$, then we have
\begin{equation*}
 \te_i A =A + \wh{E}_{i j_0} - \wh{E}_{i+1\, j_0},
\end{equation*}
where $j_0$ is the column index of $A$ corresponding the rightmost $-$ in $\td{\sigma}$. If $\td{\sigma}$ has no $-$, then we have $\te_iA={\bf 0}$. Similarly, if $\td{\sigma}$ has at least one $+$, then we have
\begin{equation*}
 \tf_i A =A - \wh{E}_{i j_1} + \wh{E}_{i+1\, j_1},
\end{equation*}
where $j_1$ is the column index of $A$ corresponding the leftmost $+$ in $\td{\sigma}$. If $\td{\sigma}$ has no $+$, then we have $\tf_iA={\bf 0}$.
\end{itemize}
Similarly, $\wh{\mc N}_{m \times n}$ has a structure of normal $A_{n-1}^{(1)}$-crystal with $P^0$ and $P^0_{\rm cl}$-weights with respect to ${\rm wt}^t, \te^t_j, \tf^t_j$ for $j\in \{\,0,1,\dots,n-1\,\}$ for $n \ge 2$, where ${\rm wt}^t(A)={\rm wt}(A^t)$, and  $\te^t_j A$, $\tf^t_j A$ are given as follows:

\begin{itemize}
\item[(1)] Let
\begin{equation*}
\sigma' = (\ \cdots,\ \underbrace{+}_{a_{i\, j}},\ \underbrace{-}_{a_{i\, j+1}},\ \underbrace{+}_{a_{i+1\, j}},\ \underbrace{-}_{a_{i+1\, j+1}},\ \cdots \ ).
\end{equation*} and let $\td{\sigma'}$ be the reduced one.
\item[(2)] If $\td{\sigma'}$ has at least one $-$, then we have
\begin{equation*}
 \te_j^t A =A + \wh{E}_{i_0\, j} - \wh{E}_{i_0\, j+1},
\end{equation*}
where $i_0$ is the row index of $A$ corresponding the rightmost $-$ in $\td{\sigma'}$. If $\td{\sigma'}$ has no $-$, then we have $\te_j^tA={\bf 0}$. Similarly, if $\td{\sigma'}$ has at least one $+$, then we have
\begin{equation*}
 \tf_j^t A =A - \wh{E}_{i_1\, j} + \wh{E}_{i_1\, j + 1},
\end{equation*}
where $i_1$ is the row index of $A$ corresponding the leftmost $+$ in $\td{\sigma'}$. If $\td{\sigma'}$ has no $+$, then we have $\tf_j^tA={\bf 0}$.
\end{itemize}
 }

The following is an analogue of Proposition \ref{prop:bocrystal}, which is straightforward to check. 
\begin{prop}
The operators $\te_i$ and $\tf_i$ for  $i\in \{\,0,1,\dots,m-1\,\}$ commute with $\te_j^t$ and $\tf_j^t$ for $j\in \{\,0,1,\dots,n-1\,\}$ on $\wh{\mc{N}}_{m\times n}\cup\{{\bf 0}\}$.
\end{prop}

For $r\ge 1$, the set $SST_{[m]}((r))$ has a structure of normal $A_{m-1}^{(1)}$-crystal with $P^0_{\rm cl}$-weights with respect to ${\rm wt}^0_{\rm cl}$, $\te_i$, $\tf_i$ for $i\in \{\,0,1,\dots,m-1\,\}$, where for $T\in SST_{[m]}((r))$
\begin{itemize}
 \item[(1)] ${\rm wt}^0_{\rm cl}(T)=\sum_{j=1}^r {\rm cl}(\epsilon_{T(j)})$,
 
 \item[(2)] $\tf_i T$ is the row semistandard tableau obtained by replacing an entry $i$ in $T$ with $i+1\pmod{n}$ if exists, and ${\bf 0}$ otherwise.
 
\end{itemize}
Recall that $SST_{[m]}((r))$ is isomorphic to the crystal of the Kirillov-Reshetikhin module corresponding to $r\,{\rm cl}(\varpi_1)$ \cite{KMN2}.

In general, for $\mu=(\mu_1,\dots,\mu_l)\in \cP$, we regard $RSST_{[m]}(\mu)$ as an $A_{m-1}^{(1)}$-crystal with $P^0_{\rm cl}$ by identifying $T\in RSST_{[m]}(\mu)$ with $T^{\ell}\ot \dots \ot T^1\in SST_{[m]}(\mu_1)\ot \dots \ot SST_{[m]}(\mu_l)$, where $T^i$ is the $i$-th row of $T$ from the bottom.

Let
\begin{equation*}
\mc{S}_{m\times n}= \bigsqcup_{\la\in \cP_n} 
RSST_{[m]}(\la')\times \mc{B}_n(\la).
\end{equation*}
We assume the normal $A_{m-1}^{(1)}$-crystal structure with $P^0_{\rm cl}$-weights for $m \ge 2$ and $A_{n-1}^{(1)}$-crystal structure with $P^0$-weights for $n \ge 2$ on $\mc{S}_{m\times n}$ in the same manner as in $\mc{T}_{m\times n}$.

\begin{thm}\label{thm:main-2'}
The bijection 
\begin{equation*}
\xymatrixcolsep{3pc}\xymatrixrowsep{0.5pc}\xymatrix{
\kappa' :\ \ \wh{\mc N}_{m\times n} \ \ar@{->}[r]  & \ \ \mc{S}_{m\times n}}
\end{equation*}
commutes with $\te_i$, $\tf_i$ for $i\in \{\,1,\dots,m-1\,\}$ and $\te^t_j$, $\tf^t_j$ for $j\in \{\,0,1,\dots,n-1\,\}$.
\end{thm}
\pf The proof is almost parallel to that of Theorem \ref{thm:main-2}, where we need to modify the arguments in Section \ref{sec:proof of main-2} with respect to the notions in Section \ref{subsec:notations for dual}. We leave the details to the reader.
\qed\newline

Let $\kappa'_1=\pi_1\circ\kappa'$, where $\pi_1$ is the projection of $\mc{S}_{m\times n}$ along the first component. Then $\kappa'_1$ commutes with $\te_i$ and $\tf_i$ for $i\in \{\,0,1,\dots,m-1\,\}$, and preserves ${\rm wt}^0_{\rm cl}$. Hence we have the following. 

\begin{cor}
A generalized dual affine permutation $A\in \wh{\mc N}_{m\times n}$ is $A_{m-1}^{(1)}$-crystal equivalent to $P^t_0$, where $\kappa'(A)=(P^t_0,Q)$. 
\end{cor}

Moreover, if we define ${\rm\bf wt}^t$ on $\wh{\mc{N}}_{m \times n}$ in the same way as in \eqref{eq:new aff wt} with respect to $\kappa'$, then we have the following analogue of Corollary \ref{cor:decomp M into An crystal}. 

\begin{cor}\label{eq:decomp N into An crystal}
If we regard $\wh{\mc N}_{m\times n}$ as an $A_{m-1}$-crystal with $P^0_{\rm cl}$-weights and as an $A_{n-1}^{(1)}$-crystal with $P^0$-weights with respect to ${\rm\bf wt}^t$, then it is an $(A_{m-1}, A_{n-1}^{(1)})$-bicrystal and $\kappa'$ is an isomorphism of $(A_{m-1}, A_{n-1}^{(1)})$-bicrystals. In particular, a generalized dual affine permutation $A\in \wh{\mc N}_{m\times n}$ is $A_{n-1}^{(1)}$-crystal equivalent to $Q$, where $\kappa'(A)=(P^t_0,Q)$. 
\end{cor}
We remark that {both $m$ and $n$ do not need to be greater} than 1 in Theorem \ref{thm:main-2'} and its corollaries. In particular, when $m=1$ we have the following multiplicity-free decomposition
\begin{equation*}
 \wh{\mc N}_{1\times n} \cong \bigoplus_{\la\in \cP_n} 
\mc{B}_n(\la),
\end{equation*}
since $RSST_{[1]}(\la')$ consists of single element for all $\la\in \cP_n$.

\section*{Declarations}
\subsection*{Conflict of Interest} The authors declare that they have no conflict of interest.

\end{document}